\documentclass[11pt]{article}

\usepackage[margin=1in]{geometry}

\usepackage{amsmath,amssymb,amsfonts}%
\usepackage{amsthm}%
\usepackage{bbm}
\usepackage{mathtools}
\usepackage{dsfont}
\usepackage{mathrsfs}
\usepackage{xcolor}

\usepackage{booktabs}       
\usepackage{nicefrac}       
\usepackage{microtype}      
\usepackage{pifont} 
\usepackage{bbm}
\usepackage{multirow}

\usepackage{natbib}
\setcitestyle{authoryear,open={(},close={)}}

\usepackage{bm}

\usepackage[hyperfootnotes=false]{hyperref}
\hypersetup{colorlinks,citecolor=blue,linkcolor=blue}

\usepackage[nameinlink, capitalize]{cleveref}

\theoremstyle{plain}
\newtheorem{theorem}{Theorem}
\newtheorem{proposition}[theorem]{Proposition}%
\newtheorem{lemma}[theorem]{Lemma}%
\newtheorem{corollary}[theorem]{Corollary}

\theoremstyle{remark}
\newtheorem{definition}{Definition}
\newtheorem{example}{Example}%
\newtheorem{remark}{Remark}%
\newcommand{\NN}{\mathbb{N}}
\newcommand{\HH}{\mathbb{H}}
\newcommand{\XX}{\mathbb{X}}
\newcommand{\YY}{\mathbb{Y}}
\newcommand{\RR}{\mathbb{R}}
\newcommand{\corr}{\mathbb{C}\mathrm{orr}}
\newcommand{\cov}{\mathbb{C}\mathrm{ov}}
\newcommand{\var}{\mathbb{V}\mathrm{ar}}
\newcommand{\diff}{\mathrm{d}}
\newcommand{\E}{\mathbb{E}}

\def\eqd{\stackrel{\mbox{\scriptsize{\textup{d}}}}{=}}

\DeclareMathOperator{\iid}{i.i.d.}
\DeclareMathOperator{\as}{a.s.}

\title{Measuring Partial Exchangeability with \\ Reproducing Kernel Hilbert Spaces}
\author{Marta Catalano\footnote{Luiss University, Italy. Email: \texttt{mcatalano@luiss.it}},~ Hugo Lavenant\footnote{Bocconi University, Italy. Email: \texttt{hugo.lavenant@unibocconi.it}}~~and Francesco Mascari\footnote{Bocconi University, Italy. Email: \texttt{francesco.mascari@phd.unibocconi.it}}}

\date{\today}

\begin{document}

\maketitle

\begin{abstract}
In Bayesian multilevel models, the data are structured in interconnected groups, and their posteriors borrow information from one another due to prior dependence between latent parameters. However, little is known about the behaviour of the dependence a posteriori. In this work, we develop a general framework for measuring partial exchangeability for parametric and nonparametric models, both a priori and a posteriori. We define an index that detects exchangeability for common models, is invariant by reparametrization, can be estimated through samples, and, crucially, is well-suited for posteriors. We achieve these properties through the use of Reproducing Kernel Hilbert Spaces, which map any random probability to a random object on a Hilbert space. This leads to many convenient properties and tractable expressions, especially a priori and under mixing. We apply our general framework to i) investigate the dependence a posteriori for the hierarchical Dirichlet process, retrieving a parametric convergence rate under very mild assumptions on the data; ii) eliciting the dependence structure of a parametric model for a principled comparison with a nonparametric alternative.

\medskip

\noindent \emph{Keywords}. Bayesian Nonparametrics, Correlation, Hierarchical Dirichlet Process, Multilevel Model, Random Probability Measure.
\end{abstract}

\maketitle

\section{Introduction} \label{sec:intro}

Bayesian modelling is widely embraced for multilevel data, characterized by distinct yet related group structures, thanks to its flexibility and natural shrinkage effect \citep{lindley_bayes_1972, efron_steins_1973,gelman_bayesian_2003, gelman_data_2007}. We consider models of the general form
\begin{equation} \label{def:multilevel_model}
    X_{i,j}|(\theta_{1},\dots,\theta_{d}) \stackrel{\iid}{\sim} P_{\theta_{i}}, \qquad (\theta_{1},\dots,\theta_{d}) \sim Q,
\end{equation}
where $X_{i,j}$ the $j$-th observation in group $i$, $Q$ is the prior distribution for the parameter vector, and $P_{\theta}$ is the data distribution parametrised by $\theta$. Thanks to de Finetti's theorem \citep{definetti_condition_1938}, this class of models is equivalent to partial exchangeability of the observations in the sense that
\[
    \bigl((X_{1,j})_{j \in \NN}, \dots, (X_{d,j})_{j \in \NN}\bigr) \eqd \bigl((X_{1,\pi_{1}(j)})_{j \in \NN}, \dots, (X_{d,\pi_{d}(j)})_{j \in \NN}\bigr),
\]
for $\pi_{1},\dots,\pi_{d}$ finite permutations of $\NN$, where $\eqd$ denotes equality in distribution.

In the inferential process, the parameters $\theta_i$ are estimated simultaneously, crucially allowing for a \textit{borrowing of information}, first introduced by Tukey as the need ``to borrow strength from either other aspects of the same body of data or from other bodies of data'' \citep{tukey_data_1972}. This feature is strictly related to the dependence among the parameters a priori. As a limiting case, when $\theta_{1} = \dots = \theta_{d}$ almost surely (a.s.), both the dependence and the borrowing are maximal: since the observations are fully exchangeable, the observations in a group carry the same information as observations in the other groups. When $\theta_1,\dots, \theta_d$ are independent, there is no borrowing of information, and the observations of the other groups will not affect the group-specific inference. The situations in between are perhaps the most interesting and require careful prior elicitation through a measure of dependence. We refer to \cite{catalano_measuring_2021, catalano_wasserstein_2024} for an account in the nonparametric setting.

The primary aim of this work is to provide a unifying framework to measure partial exchangeability that allows for a fair comparison between parametric and nonparametric multilevel models, and to investigate the behaviour of the dependence structure after observing the data. We wish to understand whether multilevel models approach or diverge from full exchangeability a posteriori as the number of observations increases, and quantify their speed of convergence. To this end, we need a principled measure of partial exchangeability that can be used a priori and a posteriori, for both finite and infinite-dimensional parameter spaces.

For simplicity, we only consider two groups of observations. If the parameters $\theta_i$ are real-valued, an intuitive measure of partial exchangeability is Pearson's linear correlation between the parameters. When the parameters have the same first and second moment, this measure detects exchangeability, in the sense that $\corr(\theta_1, \theta_2) = 1$ if and only if $(X_{1,j})_{j \in \NN}$ and $(X_{2,j})_{j \in \NN}$ are fully exchangeable. As a measure of partial exchangeability, it is not satisfactory since having the same first two moments is a strong assumption, at least a posteriori; it is not invariant under reparameterization; and clearly, it does not tackle the higher- or infinite-dimensional case. However, the following simple example provides some intuition on the results we aim to achieve for more complex models. Details are provided in the Supplementary Material.
\begin{example} \label{par_ex} 
    Consider $X_{i,j} | \theta_{1},\theta_{2} \stackrel{\iid}{\sim} \mathcal{N}\left(\theta_{i},s^{2}\right)$ for $s>0$, $j \in \NN$,  and $i = 1,2$, with $(\theta_1,\theta_2) \sim  \mathcal{N}\left(\boldsymbol{0},\tau^{2}\Sigma\right),$ where $\tau > 0$, $\Sigma_{11} = \Sigma_{22} = 1$ and $\Sigma_{12} = \Sigma_{21} = \rho \in [-1,1]$. Then, the prior correlation between $\theta_{1}$ and $\theta_{2}$ is $\corr(\theta_{1},\theta_{2}) = \rho$, while a version of the posterior correlation after $n_i$ observations in group $i$ for $i=1,2$ is
     \begin{equation*}
        \corr\bigl(\theta_{1},\theta_{1} \big| \boldsymbol{X}^{(n_1,n_2)} = \boldsymbol{x}^{(n_{1},n_{2})}\bigr) =  \frac{\rho}{\sqrt{1 + n_{1}\frac{\tau^{2}}{s^{2}}(1 - \rho^{2})}\sqrt{1 + n_{2}\frac{\tau^{2}}{s^{2}}(1 - \rho^{2})}},  
    \end{equation*}
    where here and in the following we use the compact notation $\boldsymbol{x}^{(n_{1},n_{2})} = \bigl((x_{1,j})_{j=1}^{n_{1}},(x_{2,j})_{j=1}^{n_{2}}\bigr)$.
\end{example}

We notice that the posterior correlation converges to $0$ as at least one of the sample sizes $n_{1},n_{2}$ diverges, independently of the \emph{true} data-generating process. Interestingly, it differs from the typical assumptions of asymptotic analyses, such as data being independent and identically distributed (i.i.d.) from an unknown distribution or generated by the model. The asymptotic behaviour of the correlation aligns with our intuition of borrowing information a posteriori: as the amount of data from each group increases, we rely on more evidence and thus the borrowing becomes weaker. The convergence rate of the posterior correlation is $\mathcal{O}\bigl((n_{1}n_{2})^{-1/2}\bigr)$ as $\max(n_{1},n_{2}) \rightarrow +\infty$ whenever the prior correlation $\rho$ is different from $\pm 1$, hinting to a fast decay of the borrowing of information. One of the main objectives of this paper is to establish a framework to reproduce this asymptotic investigation in a nonparametric setting as well. Nonparametric models for partially exchangeable observations are widely spread in the Bayesian literature. The first proposal dates back to \cite{cifarelli_problemi_1978}, but they have been later popularised by the seminal work of \cite{maceachern_dependent_1999, maceachern_dependent_2000}. Over the last quarter of a century, there has been a wealth of proposals, nicely reviewed in \cite{quintana_dependent_2022, wade2025bayesian}.

The first crucial idea to avoid the lack of invariance under reparametrization and extend the measure to nonparametric partially exchangeable models is to use de Finetti's Theorem \citep{definetti_condition_1938}. This foundational theorem ensures that, for any partially exchangeable sequence, the law of the parameters becomes unique once it is embedded in the space of probabilities. More precisely, models in \eqref{def:multilevel_model} are uniquely characterized by the law of the vector of random probabilities $(\tilde P_{\theta_1}, \dots, \tilde P_{\theta_d})$, which for simplicity we denote as $(\tilde P_1, \tilde P_2)$ when $d=2$. In other terms, the law of $(\tilde P_1, \tilde P_2)$ does not change under reparametrization of the model. Thus, in the sequel, we work with the model
\begin{equation}  \label{part_ex}
    X_{i,j} \big| \tilde P_{1}, \tilde P_{2} \stackrel{\iid}{\sim} \tilde P_{i}, \qquad \bigl(\tilde P_{1}, \tilde P_{2}\bigr) \sim Q,
\end{equation}
for $j \in \NN$ and $i = 1,2$, rather than with~\eqref{def:multilevel_model}. The following conceptual step is to measure dependence at the level of the random probabilities $\tilde P_{1}$ and $\tilde P_{2}$. Since full exchangeability of the observations is recovered when $\tilde P_{1} = \tilde P_{2}$ a.s., we need an index of dependence $I(Q)$ that detects almost sure equality, i.e., $I(Q) = 1$ if and only if $\tilde{P}_{1} = \tilde{P}_{2}$ a.s..

In the Bayesian nonparametric literature, there are two primary methods for measuring dependence between random measures on general Polish spaces. The first is the set-wise correlation $\corr\bigl(\tilde P_1(A), \tilde P_2(A)\bigr)$, for any measurable set $A$. Although it is based only on the first two moments of the random measures, it is the most used in practice because its expression a priori stands out for tractability and interpretability for the majority of Bayesian nonparametric models (see, e.g., \cite{rodriguez_nested_2008,leisen_vector_2013,griffin_compound_2017, camerlenghi_distribution_2019,beraha_semihierarchical_2021,ascolani_nonparametric_2023,denti_common_2023,lijoi_flexible_2023, horiguchi_tree_2024,colombi_hierarchical_2025}). In most of these settings, the set-wise correlation does not depend on the set $A$. This has been recently shown to be a property of the general class of multivariate species sampling models \citep{franzolini_multivariate_2025}, which includes nearly all priors mentioned above. However, we show that tractability can fail for random probabilities that do not belong to this class, such as those that arise from parametric models or a posteriori. In such cases, not only does the set-wise correlation depend on the choice of the set $A$, but its value for different sets can also change dramatically from $1$ to $-1$.

A second method that has been recently proposed is to use the Wasserstein distance to measure the discrepancy between $Q$ and the law in the same Fréchet class inducing full exchangeability \citep{catalano_measuring_2021,catalano_wasserstein_2024}. This method considers the full distribution of the random measures, detecting both exchangeability and independence, and can be naturally extended to a multi-group scenario. However, its tractability is limited to completely random vectors \citep{catalano_measuring_2021}, which are the natural multivariate extension of completely random measures \citep{kingman_completely_1967}. Though they are commonly used to build nonparametric priors \citep{lijoi_models_2010}, parametric priors and posteriors rarely belong to this class.

Summing up, the current proposals in the literature are not satisfactory for measuring dependence between parametric priors and a posteriori. A major objective for this work is to define a new index of dependence that detects exchangeability for common models, is well-suited for parametric and posterior random probability measures, and maintains the tractability of $\corr\bigl(\tilde P_1(A), \tilde P_2(A)\bigr)$ for Bayesian nonparametric priors. We achieve this by generalizing the set-wise correlation through the theory of Reproducing Kernel Hilbert Spaces (RKHS). These functional spaces, introduced in their general form in the seminal work by \cite{aronszajn_theory_1950}, are widely used in Machine Learning and Statistics to handle high- or infinite-dimensional data via kernels, enabling efficient computation of inner products without explicitly mapping data to higher dimensions. See \cite{scholkopf_learning_2001, berlinet_reproducing_2004, muandet_kernel_2017} for complete overviews. They also play an important role in Bayesian modelling as they allow the specification of smooth priors over function spaces, e.g., through Gaussian processes \citep{RasmussenWilliams2006} and other Bayesian kernel-based models (see, e.g., \cite{tipping_sparse_2001,sollich_bayesian_2002,pillai_characterizing_2007,maclehose_nonparametric_2009,chakraborty_bayesian_2012}), with important applications to nonparametric Bayesian modelling and functional data analysis. An interesting research direction uses RKHS embeddings to facilitate Bayesian computation, e.g., through Stein and maximum mean discrepancy (see, e.g., \cite{fukumizu2013, park2016approximate, qiang2016stein, chen2019stein, legramanti2025concentration}). Most importantly for our work, RKHS have been used to build popular measures of independence between two random variables through different summaries of the covariance operator. Some prominent examples are the kernel canonical correlation \citep{bach_kernel_2002}, the constrained covariance~\citep{gretton_kernel_2005}, 
 and the centred kernel alignment \cite{kornblith2019similarity}, a normalization of the Hilbert-Schmidt independence criterion \cite{gretton_measuring_2005} first introduced to measure similarity between kernels \cite{cristianini_kernel_2001, cortes_algorithms_2012}.
However, their setting is not directly applicable to our problem, since it cannot be used to detect exchangeability of the observations. For instance, the observations in \cref{par_ex} are exchangeable for $\rho =1$. However, none of the indices above evaluated between observations $X_{1,j_1}$ and $X_{2,j_2}$ from distinct groups are equal to $1$, as shown in \cref{tab:corr_meas}. More generally, we cannot expect to detect exchangeability by using standard measures of dependence between observations. 

\begin{table}[b!]
\centering
\caption{Estimated values of Pearson correlation, Kernelised Canonical Correlation (KCC) with regularisation 
$\varepsilon = 10^{-2}$, Constrained Covariance (COCO), and Centred Kernel Alignment (CKA) between $X_{1,j_{1}}$ and $X_{2,j_{2}}$ with samples of size $1000$ from the model in \cref{par_ex} for $s^{2} = t^{2} = 1$ and $\rho = 1$.} \label{tab:corr_meas}
\begin{tabular}{c c c c}
    \toprule
	Pearson correlation & $\mathrm{KCC}$ for $\varepsilon = 10^{-2}$ & $\mathrm{COCO}$ & $\mathrm{CKA}$ \\
	\midrule
    $0.4907$ & $0.5024$ & $0.1504$ & $0.1010$ \\
    \bottomrule
\end{tabular}
\end{table}

\subsection*{Overview and Organization of the Main Results}

In \cref{sec:rkhs}, we fix the notation and recall the main properties of Reproducing Kernel Hilbert Spaces (RKHS). These are instrumental for \cref{sec:kernel_corr}, where we define our kernel correlation starting from any symmetric positive-definite kernel $k$ on the space of observations $\XX$. Any random probability $\tilde{P}$ on $\XX$ can be mapped into a random element of $\HH_{k}$ through the kernel mean embedding $\mu_{k}(\tilde{P})$  \citep{berlinet_reproducing_2004}. Since the Hilbert structure on $\HH_{k}$ determines a natural notion of correlation, which we denote as $\corr_{\HH_{k}}$, we measure partial exchangeability as the Hilbert correlation between the kernel mean embeddings of $\tilde{P}_{1}$ and $\tilde{P}_{2}$. Using the previous terminology, our index $I_{k}(Q)$ is defined as
\[
    \corr_{k}\bigl(\tilde{P}_{1}, \tilde{P}_{2}\bigr)  := \corr_{\HH_{k}}\bigl(\mu_{k}(\tilde{P}_{1}),\mu_{k}(\tilde{P}_{2})\bigr).
\]
Such a construction can be made for several choices of the kernel $k$. Classical examples include linear, Gaussian, and Laplace kernels. A curious note is that by taking the set-wise kernel $k(x,y) = \mathbbm{1}_{A}(x)\mathbbm{1}_{A}(y)$ for some measurable set $A$, we recover the standard notion of set-wise correlation $\corr\bigl(\tilde{P}_{1}(A),\tilde{P}_{2}(A)\bigr)$. This is thoroughly discussed in \cref{sec:setwise}, where it allows us to draw interesting parallels between our proposal and the widely used set-wise correlation.

In \cref{sec:exchangeability}, we identify settings under which the kernel correlation detects exchangeability, i.e., $\corr_{k}\bigl(\tilde{P}_{1}, \tilde{P}_{2}\bigr) = 1$ if and only if $\tilde{P}_{1} = \tilde{P}_{2}$ a.s.. In particular, the kernel must induce an injective kernel mean embedding on the set of bounded signed measures. This is guaranteed for Gaussian and Laplace kernels, but not for linear or set-wise kernels.

In \cref{sec:estimator}, we express the kernel correlation in terms of partially exchangeable observables. Remarkably, we show that it can be determined by two observables for each group. This inspires a natural asymptotically normal estimator for kernel correlation through independent copies.

In \cref{sec:setwise}, we identify a key structural assumption on the random measures for which the kernel correlation does not depend on the choice of kernel, and we underline that this assumption holds for multivariate species sampling models. Since the set-wise correlation is a specific choice of kernel, this implies (i) the known result that the set-wise correlation does not depend on the set $A$ for multivariate species sampling models \citep{franzolini_multivariate_2025}, (ii) that our kernel correlation coincides with the set-wise correlation for this class of models, thus recovering the interpretable and tractable expressions known in the literature.

In \cref{sec:mixtures}, we explore the performance of $\corr_{k}$ on mixture models \citep{ferguson_bayesian_1983, lo_class_1984, escobar_bayesian_1995}, which are widely used in Bayesian inference for density estimation and clustering. We give conditions to express the kernel correlation of the mixture in terms of the mixing measures with respect to an updated version of the kernel. Our results imply that the kernel correlation remains invariant under mixing for multivariate species sampling models and detects exchangeability for mixture models. We extend this property to parametric models by reinterpreting them as a special case of mixture models.

In \cref{sec:hdp}, we investigate the posterior behaviour of $\corr_{k}$ for the archetype of Bayesian nonparametric hierarchical models, the Hierarchical Dirichlet Process \citep{teh_hierarchical_2006}. Under a \emph{non-degeneracy} assumption on the data, which depends on the kernel and is mild only for injective kernels, we prove that our kernel correlation goes to 0 as $\mathcal{O}\bigl((n_{1}n_{2})^{-1/2}\bigr)$. Notably, the infinite dimensionality of the parameters does not slow down the convergence rate of \cref{par_ex}, which holds for basically any data-generating mechanism and does not depend on the dimension of the observation space. Our theoretical results are also confirmed by numerical simulations in \cref{sec:numerics}. We apply a Gibbs sampling scheme to compute the kernel correlation for different choices of the kernel, and we empirically demonstrate how the Gaussian and Laplace kernels are robust to the choice of hyperparameters, whereas the set $A$ has a significant impact on the value of the set-wise correlation.

Finally, in \cref{sec:models}, we apply the kernel correlation to perform a model comparison between the Gaussian model in \cref{par_ex} and the hierarchical Dirichlet Process. Specifically, we show how the conclusions of the models depend greatly on the value of the kernel correlation, and that one should set its value to be the same in both models for a fair model comparison.

\section{Preliminaries on Reproducing Kernel Hilbert Spaces} \label{sec:rkhs}

In this section, we establish the notation and recall the main properties of Reproducing Kernel Hilbert Spaces that will be used throughout the remainder of the work.

Let $\XX$ be a locally compact Polish Space, endowed with its Borel $\sigma$-algebra $\mathcal{X}$. We denote the space of bounded signed measures on $\XX$ as $\mathcal{M}_{b}(\XX)$, and the subspace of probability measures as $\mathcal{P}(\XX)$, endowing both these spaces with the $\sigma$-algebra induced by evaluation maps $\varphi_{A}(\xi) =  \xi(A)$, for any $A \in \mathcal{X}, \xi \in \mathcal{M}_{b}(\XX)$. 

A kernel $k: \XX \times \XX \rightarrow \RR$ is a measurable function that, in the sequel, will always be assumed to be bounded, symmetric, and positive-definite; see the Supplementary Material for details. Any kernel defines a natural mapping of $\XX$ into the space $\RR^{\XX}$ of functions from $\XX$ to $\RR$ through the feature maps $x \mapsto k(x,\cdot)$. The natural notion of inner product between feature maps $\langle k(x,\cdot), k(y,\cdot) \rangle := k(x,y)$ can be extended to all linear combinations of feature maps. The closure of this space in $\RR^{\XX}$ defines the unique Reproducing Kernel Hilbert Space (RKHS) $\HH_{k}$ induced by the kernel $k$ \citep{aronszajn_theory_1950}. We denote its scalar product as $\langle \cdot, \cdot \rangle_{\HH_{k}}$.

For a kernel $k$, every bounded signed measure $\xi \in \mathcal{M}_{b}(\XX)$ can be associated with a unique element $\mu_{k}(\xi) \in \HH_{k}$ through its kernel mean embedding
\begin{equation} \label{def:kernel_mean_embedding}
    \mu_{k}(\xi)(y) := \int_{\XX}{k(x,y)\diff \xi(x)}.
\end{equation}
The reproducing property ensures the following useful identities for $\xi,\xi_1,\xi_2 \in \mathcal{M}_{b}(\XX)$ and $f \in \HH_{k}$,
\begin{equation} \label{reproduce_kernel_int}
    \langle f, \mu_{k}(\xi) \rangle_{\HH_{k}} = \int_{\XX}{f(x)\diff\xi(x)}, \qquad \langle \mu_{k}(\xi_1),\mu_k(\xi_2) \rangle_{\HH_{k}} = \iint_{\XX \times \XX} k(x,y) \diff \xi_1(x) \diff \xi_2(y).
\end{equation}

A kernel $k$ is $c_{0}$ if it vanishes at infinity, that is, if the set $\{x: |k(x,y)| \geq \varepsilon \}$ is compact for every $\varepsilon > 0, y \in \XX$; it is injective if the feature map $x \mapsto k(x,\cdot)$ is injective; it is characteristic if the kernel mean embedding is injective on $\mathcal{P}(\XX)$. Since $k(x,\cdot) = \mu_k(\delta_x)$, a characteristic kernel is always injective. For a $c_{0}$-kernel, the kernel mean embedding is injective on $\mathcal{M}_b(\XX)$ if and only if $\HH_k$ is dense in $C_{0}(\XX)$, the space of continuous functions over $\XX$ vanishing at infinity endowed with the uniform norm \citep{sriperumbudur_universality_2011}. In this case, we call $k$ $c_{0}$-universal, which implies that $k$ is characteristic since $\mathcal{P}(\XX) \subset \mathcal{M}_{b}(\XX)$.

\begin{table}[b!]
\centering
\caption{Useful properties of kernels: injectivity (\textsc{inj}), continuity (\textsc{cont}), characteristic (\textsc{Char}), and $c_{0}$-universality (\textsc{$c_{0}$-univ}). The set-wise kernel is introduced in \cref{sec:setwise}.} \label{tab:kernels}
\resizebox{\textwidth}{!}{
\begin{tabular}{l l c c c c c}
    \toprule
    \multicolumn{3}{c}{Type of kernel} & \multicolumn{4}{c}{Property of kernel} \\
    \cmidrule(rl){1-3} \cmidrule(rl){4-7}
	Name & Expression for $k(x,y)$ & Space $\XX$ & \textsc{cont} & \textsc{inj} & \textsc{char} &  \textsc{$c_{0}$-univ} \\
	\midrule
    \multirow{2}{*}{Linear} & \multirow{2}{*}{$\langle x,y \rangle_{\HH}$} & \raisebox{-2ex}{\parbox[c]{3cm}{  \centering \setlength{\baselineskip}{0.5\baselineskip} Bounded subset \\ of Hilbert $\HH$}} & \multirow{2}{*}{\ding{51}} & \multirow{2}{*}{\ding{51}} & \multirow{2}{*}{\ding{55}} & \multirow{2}{*}{\ding{55}} \vspace{0.1cm} \\ 
	Gaussian &  $\exp(-\|x-y\|_{\XX}^{2}/(2 \sigma^{2}))$, $\sigma > 0$ & Hilbert & \ding{51} & \ding{51} & \ding{51} & \ding{51} \\
	Laplace & $\exp(-\|x-y\|_{1}/\beta)$, $\beta > 0$ & $\RR^{m}$ & \ding{51} & \ding{51} & \ding{51} & \ding{51} \\
	Set-wise & $\mathbbm{1}_{A}(x)\mathbbm{1}_{A}(y)$, $A \in  \mathcal{X}$ & Measure Space
    & \ding{55} & \ding{55} & \ding{55} & \ding{55} \\
    \bottomrule
\end{tabular}}
\end{table}

In this work, we focus on some of the most prominent examples of kernels in the literature: the linear kernel $k(x,y) = \langle x, y \rangle_{\HH}$, defined on a bounded subset $\XX$ of a Hilbert space $\HH$, the Gaussian kernel $k(x,y) = \exp(-\|x-y\|_{\XX}^{2}/(2\sigma^{2}))$ for some $\sigma > 0$, defined on any Hilbert space $\XX$, and the Laplace kernel $k(x,y) = \exp(-\|x-y\|_{1}/\beta)$ for some $\beta > 0$, defined on $\XX = \RR^{m}$, where $\|x-y\|_{1} = |x_{1} - y_{1}| + \dots + |x_{m} - y_{m}|$. We refer to \citet{muandet_kernel_2017} for a complete reference on common kernels on Euclidean spaces and \citet{guella2022gaussian} for an analysis of Gaussian kernels on Hilbert spaces. Both Gaussian and Laplace kernels are injective and $c_{0}$-universal, while the linear kernel is injective and continuous, but not $c_{0}$-universal. \cref{tab:kernels} provides a summary of the main properties of the kernels presented above and the set-wise kernel, which is introduced later in \cref{sec:setwise}.

\section{Kernel Correlation} \label{sec:kernel_corr}

In this section, we define our kernel correlation between random probability measures on $(\XX, \mathcal{X})$ endowed with a kernel $k$.  A random probability measure $\tilde{P}$ on $\XX$ is a random element on $\mathcal{P}(\XX)$. Its mean measure $\E[\tilde{P}]$ is the probability measure that satisfies $\int{f(x)\diff \E[\tilde{P}](x)} = \E[\int{f(x)\diff \tilde{P}(x)}]$ for any bounded and measurable function $f : \XX \rightarrow \mathbb{R}$. Most random probabilities that we mention in this work are defined on $\XX$. We will not repeat this assumption unless needed.

Given two random probabilities $\tilde{P}_{1}$ and $\tilde{P}_{2}$, their kernel mean embeddings $\mu_{k}(\tilde{P}_{1})$ and $\mu_{k}(\tilde{P}_{2})$ defined in \eqref{def:kernel_mean_embedding} are random elements on $\HH_{k}$. We define the kernel covariance between $\tilde{P}_{1}$ and $\tilde{P}_{2}$ as
\[
    \cov_{k}\bigl(\tilde{P}_{1},\tilde{P}_{2}\bigr) := \E\left[\bigl\langle \mu_{k}(\tilde{P}_{1}) - \E\bigl[\mu_{k}(\tilde{P}_{1})\bigr], \mu_{k}(\tilde{P}_{2}) - \E\bigl[\mu_{k}(\tilde{P}_{2})\bigr]\bigr\rangle_{\HH_{k}}\right],
\]
which is well-defined since $k$ is bounded and thus $\E\bigl[\|\mu_{k}(\tilde{P}_{i})\|_{\HH_{k}}^{2}\bigr] < +\infty$ for $i=1,2$.

\begin{remark}
    For any two random variables $X,Y$ on a Hilbert space $\HH$ with $\E\bigl[\|X\|_{\HH}^{2}\bigr], \E\bigl[\|Y\|_{\HH}^{2}\bigr] < +\infty$, the cross-covariance operator between $X$ and $Y$, $C_{X,Y} := \E[(X - \E[X] )\otimes (Y - \E[Y])]:\HH \to \HH$, is defined as $C_{X,Y}(h) = \E\bigl[ (X - \E[X])\langle Y - \E[Y], h \rangle_{\HH}\bigr]$ for $h \in \HH$. It follows that $\cov_{k}\bigl(\tilde{P}_{1},\tilde{P}_{2}\bigr)$ can be interpreted as the trace of the cross-covariance operator between $\mu_k(\tilde{P}_1)$ and $\mu_k(\tilde{P}_2)$ in $\HH_k$.
\end{remark}

Explicit calculations of the kernel covariance are often possible thanks to the following integral representation.

\begin{proposition} \label{int_ker_cov_var}
    Let $\tilde{P}_{1}$, $\tilde{P}_{2}$ be random probabilities with $P_{0,i} = \E\bigl[\tilde{P}_{i}\bigr]$ for $i=1,2$. Then,
    \[
        \cov_{k}\bigl(\tilde{P}_{1},\tilde{P}_{2}\bigr) = \E\left[\iint k(x,y)\diff \tilde{P}_{1}(x) \diff \tilde{P}_{2}(y) \right] - \iint{k(x,y)\diff P_{0,1}(x)\diff P_{0,2}(y)}.
    \]
\end{proposition}

The definition of $\cov_k$ implies properties similar to those of the standard covariance between real-valued random variables. In particular, the kernel covariance is a symmetric bilinear form and it inherits the law of total covariance of Hilbert spaces: for a random variable $Z$ defined on the same probability space of $\tilde{P}_{1}$ and $\tilde{P}_{2}$,
\[
    \cov_{k}\bigl(\tilde{P}_{1},\tilde{P}_{2}\bigr) = 	\E_{Z}\left[\cov_{k}\bigl(\tilde{P}_{1},\tilde{P}_{2} \big| Z \bigr)\right] + \cov_{k}\left(\E\bigl[\tilde{P}_{1} \big| Z\bigr],\E\bigl[\tilde{P}_{2} \big| Z\bigr]\right).
\]

From the definition of kernel covariance, we derive the notion of kernel variance of a random probability measure $\tilde{P}$ as $\var_{k}(\tilde{P}):= \cov_{k}(\tilde{P},\tilde{P})$. As with usual random variables, a zero variance characterizes deterministic objects. We recall that a random probability $\tilde P$ is deterministic if $\tilde{P} = P$ a.s. for some $P \in \mathcal{P}(\XX)$.

\begin{lemma} \label{lemma:var_equal_0}
    For any random probability $\tilde{P}$, $\var_{k}(\tilde{P}) \geq 0$. If the kernel $k$ is characteristic, then $\var_{k}(\tilde{P}) = 0$ if and only if $\tilde{P}$ is deterministic.
\end{lemma}

We now have all the main ingredients for the definition of kernel correlation. 

\begin{definition} \label{ker_corr}
    Let $\tilde{P}_{1}$, $\tilde{P}_{2}$ be random probabilities such that $\var_{k}\bigl(\tilde{P}_{i}\bigr) > 0$ for $i=1,2$. The \emph{kernel correlation} between $\tilde{P}_{1}$ and $\tilde{P}_{2}$ induced by the kernel $k$ is defined as
    \[
        \corr_{k}\bigl(\tilde{P}_{1}, \tilde{P}_{2}\bigr):= \frac{\cov_{k}\bigl(\tilde{P}_{1},\tilde{P}_{2}\bigr)}{\sqrt{\var_{k}\bigl(\tilde{P}_{1}\bigr)}\sqrt{\var_{k}\bigl(\tilde{P}_{2}\bigr)}}.
    \]
\end{definition}

\begin{proposition} \label{extreme_case_indep}
    The kernel correlation $\corr_{k}\bigl(\tilde{P}_{1}, \tilde{P}_{2}\bigr)$ takes values in $[-1,1]$. Moreover, if $\tilde{P}_{1}$ and $\tilde{P}_{2}$ are independent, $\corr_{k}\bigl(\tilde{P}_{1}, \tilde{P}_{2}\bigr) = 0$. 
\end{proposition}

We observe that the kernel correlation can be considered as a generalization of Pearson's correlation between parameters, as considered in \cref{par_ex}. Indeed, by identifying any random real parameter $\theta$ with the random probability $\delta_\theta$, then $\corr(\theta_{1},\theta_{2}) = \corr_{k}\bigl(\delta_{\theta_{1}},\delta_{\theta_{2}}\bigr)$ for $k(x,y) = xy$ the linear kernel on $\XX = \Theta \subset \RR$ bounded. However, as shown in the next section, the linear kernel is unable to detect the full exchangeability of the observations, which corresponds to $\tilde P_1 = \tilde P_2$ a.s.. We devote the next section to showing that other kernels do not suffer from this limitation.

\section{Detecting Full Exchangeability} \label{sec:exchangeability}

We now investigate under which assumptions our kernel correlation can detect the almost sure equality of two random probability measures. We illustrate two different conditions, one that holds for many Bayesian nonparametric priors and the other that typically holds for the corresponding posteriors.

\begin{lemma} \label{linear_corr}
    Let $k$ be a $c_{0}$-universal kernel. Then $\corr_{k}\bigl(\tilde{P}_{1},\tilde{P}_{2}\bigr)=\pm1$ if and only if $\tilde{P}_{1} - \E[\tilde{P}_{1}] = \alpha\bigl(\tilde{P}_{2} - \E\bigl[\tilde{P}_{2}\bigr]\bigr)$ a.s. for some $\alpha \in \RR \setminus \{0\}$. The sign of $\alpha$ and the one of $\corr_{k}\bigl(\tilde{P}_{1},\tilde{P}_{2}\bigr)$ coincide. 
\end{lemma}

This result is pivotal for the kernel correlation to identify full exchangeability both a priori and a posteriori for a.s. discrete random probability measures.

\begin{theorem} \label{index_equal1_prior}
    Let $\tilde{P}_{i}$ be an a.s. discrete random probability such that $\E\bigl[\tilde{P}_{i}\bigr]$ is atomless, for $i=1,2$. If $k$ is $c_0$-universal, then $\corr_{k}\bigl(\tilde{P}_{1},\tilde{P}_{2}\bigr) = 1$ if and only if $\tilde{P}_{1} = \tilde{P}_{2}$ a.s..
\end{theorem}

We observe that \cref{index_equal1_prior} does not put assumptions on the dependence between the random probabilities, but only on their marginal distribution. This will be crucial in Section~\ref{sec:mixtures}, where we will prove that the kernel correlation detects exchangeability also for mixtures and parametric models.

The assumptions of \cref{index_equal1_prior} hold for several classes of nonparametric priors, including the most common specifications of normalized random measures with independent increments \citep{regazzini_distributional_2003} or species sampling processes \citep{pitman_developments_1996}. Example~SM1 in the Supplementary Material shows that the assumption of $c_{0}$-universality for $k$ cannot be removed. On the other hand, the assumption of atomless mean measure can be relaxed, which is extremely useful when dealing with posterior distributions. Indeed, for many Bayesian nonparametric models, the posterior mean measure is a convex combination of an atomless measure and a discrete measure supported on the observed values. Recall that $z \in \XX$ is a fixed atom for the random measure $\tilde{P}$ if $\mathbb{P}\bigl(\tilde{P}(\{z\}) > 0\bigr) > 0$.

\begin{theorem} \label{index_equal1_post}
    Let $\tilde{P}_{i}$ be an a.s. discrete random probability such that for any fixed atom $z_i$ and for any $\varepsilon > 0$, $\mathbb{P}\bigl(\tilde{P}_{i}(\{z_i\}) < \varepsilon\bigr) > 0$, for $i = 1,2$. If $k$ is $c_0$-universal, then $\corr_{k}\bigl(\tilde{P}_{1},\tilde{P}_{2}\bigr) = 1$ if and only if $\tilde{P}_{1} = \tilde{P}_{2}$ a.s..
\end{theorem}

The assumption on the fixed atoms guarantees that the realizations of the corresponding jumps can be arbitrarily small. It holds, for example, for the posterior of normalized random measures with independent increments \citep{james_posterior_2009}. Example~SM2 in the Supplementary Material shows that the hypothesis cannot be removed.

\section{Estimation from Samples} \label{sec:estimator}

In this section, we express the kernel correlation in terms of the partially exchangeable observables in \eqref{part_ex} with an unforeseen take-home message. Whereas recovering the marginal distribution of the random probabilities requires an infinite sample from a partially exchangeable sequence, computing their kernel correlation needs only four observations, two for each group. This characterization also allows us to estimate the kernel correlation from samples generated by the model using a convenient asymptotically normal estimator.

\begin{theorem} \label{prop_cov_obs}
    Let $(X_{1,1},X_{2,1},X_{1,2},X_{2,2})$ be partially exchangeable observations from \eqref{part_ex}. Then, $\cov_{k}\bigl(\tilde{P}_{1},\tilde{P}_{2}\bigr) = \cov_{\HH_{k}}(k(X_{1,1},\cdot),k(X_{2,1},\cdot))$ and $\var_{k}\bigl(\tilde{P}_{i}\bigr) = \cov_{\HH_{k}}(k(X_{i,1},\cdot), k(X_{i,2},\cdot))$ for $i = 1,2$. In particular,
    \begin{equation*} 
        \corr_{k}(\tilde{P}_{1},\tilde{P}_{2}) = \frac{ \cov_{\HH_{k}}(k(X_{1,1},\cdot),k(X_{2,1},\cdot))}{\sqrt{\cov_{\HH_{k}}(k(X_{1,1},\cdot), k(X_{1,2},\cdot))} \sqrt{\cov_{\HH_{k}}(k(X_{2,1},\cdot), k(X_{2,2},\cdot))}}.
    \end{equation*}
\end{theorem}
    
\begin{remark} \label{rem:dependence_observations}
    If $\XX \subset \RR$ bounded and $k(x,y) = xy$ is the linear kernel, then we obtain $\corr_k\bigl(\tilde P_{1},\tilde P_{2}\bigr) = \cov(X_{1,1}, X_{2,1})/ \sqrt{\cov(X_{1,1}, X_{1,2}) \cov(X_{2,1}, X_{2,2})}$. Note that it differs from $\corr(X_{1,1}, X_{2,1})$, Pearson's linear correlation between $X_{1,1}$ and $X_{2,1}$. If $\tilde{P}_1 = \tilde{P}_2 = \tilde P$ a.s. then $\corr(X_{1,1}, X_{2,1}) < 1$ if $\tilde P \neq \delta_X$ for some random variable $X$, as $(X_{1,2},X_{2,1})$ are conditionally independent and not a.s. equal. In contrast, $\corr_k(\tilde P_{1},\tilde P_{2}\bigr) = 1$ as we showed in \cref{sec:exchangeability}. Hence, we need to go beyond $\corr(X_{1,1},X_{2,1})$ to detect exchangeability of the model, as underlined in the introduction.  
\end{remark}

We note that $\cov_{k}\bigl(\tilde{P}_{1},\tilde{P}_{2}\bigr)$ only requires the joint law of one observation in each group, while $\var_{k}\bigl(\tilde{P}_{i}\bigr)$ requires the joint law of two observations in group $i$. Thus, \cref{prop_cov_obs} suggests an estimator of the kernel correlation through independent $2\times2$ samples of a partially exchangeable sequence. In practice, these samples can be easily obtained whenever one can sample from the law of $\bigl(\tilde P_1, \tilde P_2\bigr)$ directly, as with most parametric models. When $\tilde P_1$ and $\tilde P_2$ are infinite-dimensional, one can either find a finite-dimensional approximation of the law of $\bigl(\tilde P_1, \tilde P_2\bigr)$ or find the predictive distribution by integrating out the random probabilities. This can be seen as the partially exchangeable generalization of the Blackwell-MacQueen urn scheme \citep{blackwell_ferguson_1973}. Most partially exchangeable models in the literature have explicit expressions for the predictive distribution, both a priori and a posteriori.

\begin{proposition} \label{prop_cov_unbiased}
    Let $\bigl(X_{1,1}^{(t)},X_{2,1}^{(t)},X_{1,2}^{(t)},X_{2,2}^{(t)}\bigr)$ be independent partially exchangeable observations from \eqref{part_ex}, for $t = 1,\dots,M$. Then,
    \begin{align*}
        \widehat{\cov}_{k,M}\bigl(\tilde{P}_{1},\tilde{P}_{2}\bigr) &:= \frac{1}{M-1}\sum_{t=1}^{M}{k\bigl(X^{(t)}_{1,1},X^{(t)}_{2,1}\bigr)} - \frac{1}{(M-1)M}\sum_{t=1}^{M}{\sum_{s=1}^{M}{k\bigl(X^{(t)}_{1,1},X^{(s)}_{2,1}\bigr)}},\\
        \widehat{\var}_{k,M}\bigl(\tilde{P}_{i}\bigr) &:= \frac{1}{M-1}\sum_{t=1}^{M}{k\bigl(X^{(t)}_{i,1},X^{(t)}_{i,2}\bigr)} - \frac{1}{(M-1)M}\sum_{t=1}^{M}{\sum_{s=1}^{M}{k\bigl(X^{(t)}_{i,1},X^{(s)}_{i,2}\bigr)}}
    \end{align*}
    are unbiased estimators of $\cov_{k}(\tilde{P}_{1},\tilde{P}_{2})$ and $\var_{k}\bigl(\tilde{P}_{i}\bigr)$ respectively, for $i = 1,2$.
\end{proposition}

Combining these two quantities, we obtain an estimator of the correlation, which notably preserves both the rate and the asymptotic normality of the parametric case.

\begin{proposition} \label{prop_cor_asymp_normal}
    With the notations of \cref{prop_cov_unbiased}, if $\var_{k}\bigl(\tilde{P}_i\bigr) > 0$ for $i=1,2$, 
    \[
        \widehat{\corr}_{k,M}\bigl(\tilde{P}_{1},\tilde{P}_{2}\bigr) := \frac{\widehat{\cov}_{k,M}\bigl(\tilde{P}_{1},\tilde{P}_{2}\bigr)}{\sqrt{\widehat{\var}_{k,M}\bigl(\tilde{P}_{1}\bigr)} \sqrt{\widehat{\var}_{k,M}\bigl(\tilde{P}_{2}\bigr)}}
    \]
    is an asymptotically normal estimator of the kernel correlation, i.e., $\sqrt{M}\bigl(\widehat{\corr}_{k,M}\bigl(\tilde{P}_{1},\tilde{P}_{2}\bigr) - \corr_{k}\bigl(\tilde{P}_{1},\tilde{P}_{2}\bigr)\bigr)$ converges in distribution to a centred Gaussian distribution as $M \to + \infty$.
\end{proposition}

\section{From Set-Wise to Kernel Correlation} \label{sec:setwise}

In this section, we interpret the kernel correlation as a generalization of the widely used set-wise correlation $\corr\bigl(\tilde{P}_{1}(A),\tilde{P}_{2}(A)\bigr)$, for a measurable set $A \in \mathcal{X}$. Despite having some undesirable behaviours due to the lack of continuity, we show that its expression is equal to the one for $c_0$-universal kernels for a large class of Bayesian nonparametric priors.

\begin{proposition} \label{setwise_kernel_corr}
    For $A \in  \mathcal{X}$, the set-wise kernel $k_A(x,y):= \mathbbm{1}_{A}(x)\mathbbm{1}_{A}(y)$ defines a kernel such that $\corr_{k_{A}}\bigl(\tilde{P}_{1},\tilde{P}_{2}\bigr) = \corr\bigl(\tilde{P}_{1}(A),\tilde{P}_{2}(A)\bigr)$.
\end{proposition}

As a direct consequence, note that $\var_{k_A}(\tilde{P}) = 0$ if and only if $\tilde{P}(A)$ is deterministic. In contrast to most of the standard kernels, $k_A$ is not continuous, injective, or characteristic (cf. \cref{tab:kernels}). In particular, it is not $c_{0}$-universal, and thus it does not fulfil the conditions of \cref{index_equal1_prior} and \cref{index_equal1_post}. Unsurprisingly, we may find $\tilde{P}_{1}$ and $\tilde{P}_{2}$ that are not a.s. equal such that $\corr\bigl(\tilde{P}_{1}(A),\tilde{P}_{2}(A)\bigr) = 1$ for some $A$, as shown in Example~SM3 in the Supplementary Material. Even more strikingly, the correlation may change from $-1$ to $+1$ for the \emph{same} pair of random measures simply by changing the set $A$.

\begin{example} \label{post_ex}
    For $W \sim \mathrm{Unif}_{[0,1]}$, $P \in \mathcal{P}(\XX)$ an atomless probability, and $x_{1} \neq x_{2} \in \XX$, we define $\tilde{P}_{i} := W\delta_{x_{i}} + (1 - W)P$ for $i = 1,2$. If we take $A \in \mathcal{X}$ such that $x_{1},x_{2} \notin A$ and $P(A) \neq 0$, then $\tilde{P}_1(A) = \tilde{P}_2(A) = (1-W) P(A)$ a.s.. Thus $\corr\bigl(\tilde{P}_{1}(A),\tilde{P}_{2}(A)\bigr) = 1$. If we take $B \in \mathcal{X}$ such that $P(B) \notin \{0,1\}$, $x_{1} \in B$, $x_{2} \notin B$, then $\tilde{P}_1(B) = P(B) + W (1- P(B))$ while $\tilde{P}_2(B) = (1-W) P(B)$. It follows that $\corr\bigl(\tilde{P}_{1}(B),\tilde{P}_{2}(B)\bigr) = -1$.
\end{example}

Moreover, the lack of continuity of the set-wise kernel leads to a lack of continuity of the kernel correlation, which, in turn, compromises the stability in the assessment of the measure of partial exchangeability, as shown by a slight modification of the above example in the Supplementary Material (Example~SM4). These examples strongly advocate for the use of a $c_0$-universal kernel, such as the Gaussian or the Laplace, which provides more stable measures and better detection of full exchangeability. 

Nevertheless, there is a large class of Bayesian nonparametric priors where these advantages are not necessary. We now identify a structural assumption on the random probabilities, which holds for most commonly used priors in Bayesian Nonparametrics, where the kernel correlation does not depend on the choice of the kernel. Thus, for this class of random probabilities, the kernel correlation coincides with the set-wise correlation for any choice of an injective kernel.

\begin{proposition}  \label{cov_set_ker}
    Let $\tilde{P}_{1}$, $\tilde{P}_{2}$ be random probabilities with same mean measure $\E\bigl[\tilde{P}_{1}\bigr] = \E\bigl[\tilde{P}_{2}\bigr] = P_0$. Then, the following conditions are equivalent for any $\eta \in \RR$ and imply that $\eta \in [-1,1]$:
    \begin{itemize}
        \item[(i)] $\cov\bigl(\tilde{P}_{1}(A),\tilde{P}_{2}(A)\bigr) = \eta P_{0}(A)(1 - P_{0}(A))$ for any measurable set $A \in \mathcal{X}$;
        \item[(ii)] $\cov_{k}\bigl(\tilde{P}_{1},\tilde{P}_{2}\bigr) = \eta\bigl(\int{k(x,x)\diff P_{0}(x)} - \iint{k(x,y)\diff P_{0}(x)\diff P_{0}(y)}\bigr)$ for any kernel $k$.
    \end{itemize}
\end{proposition}

Applying \cref{cov_set_ker} to $\tilde{P}_{1} = \tilde{P}_{2} = \tilde{P}$ a.s., we deduce the following corollary for the variance. 

\begin{corollary}  \label{var_set_ker}
    Let $\tilde{P}$ be a random probability with $P_{0} = \E[\tilde{P}]$. Then, the following conditions are equivalent for any $\lambda \in \RR$ and imply that $\lambda \in [0,1]$:
    \begin{itemize}
        \item[(i)] $\var\bigl(\tilde{P}(A)\bigr) = \lambda P_{0}(A)(1 - P_{0}(A))$ for any measurable set $A \in \mathcal{X}$;
        \item[(ii)] $\var_{k}\bigl(\tilde{P}\bigr) = \lambda\bigl(\int{k(x,x)\diff P_{0}(x)} - \iint{k(x,y)\diff P_{0}(x)\diff P_{0}(y)}\bigr)$ for any  kernel $k$.
    \end{itemize}
\end{corollary}

If the conditions in \cref{cov_set_ker} and \cref{var_set_ker} are met, then the kernel correlation does not depend on the kernel as long as it is well defined, that is, if the kernel variances are strictly positive (cf. \cref{lemma:var_equal_0}). As an easy corollary, we deduce the following fundamental theorem.

\begin{theorem}  \label{corr_set_ker}
    Let $\tilde{P}_{1}$, $\tilde{P}_{2}$ be non-deterministic random probabilities that satisfy the conditions in \cref{cov_set_ker} for some $\eta \in [-1,1]$, and the ones in \cref{var_set_ker} for some $\lambda_{i} \in (0,1]$, for $i = 1,2$. Then for any injective kernel $k$ and any set $A \in \mathcal{X}$ such that $\tilde{P}_i(A)$ is not deterministic for $i = 1,2$,
    \[
        \corr_{k}\bigl(\tilde{P}_{1},\tilde{P}_{2}\bigr) = \corr\bigl(\tilde{P}_{1}(A),\tilde{P}_{2}(A)\bigr) = \frac{\eta}{\sqrt{\lambda_{1}\lambda_{2}}}.
    \]
\end{theorem}

\begin{remark}
    In \cref{corr_set_ker}, we only need the kernel $k$ to be injective so that $\var_{k}\bigl(\tilde{P}_{i}\bigr)$ is non-null for $i = 1,2$, whereas \cref{lemma:var_equal_0} requires the stronger assumption that $k$ is characteristic. However, since $\var_{k}\bigl(\tilde{P}_{i}\bigr)$ can be written as in \cref{var_set_ker}, $\tilde{P}_{i}$ is deterministic if and only if there exists $z \in \XX$ such that $\tilde{P}_{i} = \delta_{z}$ a.s.. The details are in the proof of \cref{corr_set_ker} 
    in the Supplementary Material. 
\end{remark}

Remarkably, the results in \cite{franzolini_multivariate_2025} guarantee that \cref{corr_set_ker} holds for multivariate species sampling processes. This class of multivariate priors includes most of the partially exchangeable models studied in the Bayesian nonparametric literature. For example, in the context of hierarchical models, it includes the hierarchical Dirichlet process \citep{teh_hierarchical_2006}, hierarchical normalized completely random measures \citep{camerlenghi_distribution_2019}, the semi-hierarchical Dirichlet process \citep{beraha_semihierarchical_2021}, and the hidden hierarchical Dirichlet process \citep{lijoi_flexible_2023}, to name a few; we refer to \cite{franzolini_multivariate_2025} for the complete list.

However, the framework of \cref{corr_set_ker} is not enough to analyze parametric and posterior random measures of nonparametric models. Indeed, even if $\tilde{P}_1$ and $\tilde{P}_2$ are multivariate species sampling priors, the posteriors do not belong to this class.

\section{Kernel Correlation for Mixture Models} \label{sec:mixtures}

Mixture models are widely used in Bayesian inference for density estimation and clustering. The use of a nonparametric a.s. discrete prior as a mixing distribution allows for a potentially infinite number of components, with an evident gain in flexibility. In this section, we express the kernel correlation between mixture models in terms of a kernel correlation between mixing distributions with an updated kernel. Moreover, we show that for multivariate species sampling models, the kernel correlation between mixture models coincides with the set-wise correlation between the mixing measures. This validates a standard procedure in applied analyses with mixture models, where the prior elicitation of the borrowing of information is performed through the dependence between the mixing measures. Remarkably, our analysis of mixture models can be used to treat the kernel correlation of parametric models, showing that it detects exchangeability also in these settings.

In this section, $(\XX, \mathcal{X})$ is a latent space, and $(\YY, \mathcal{Y})$ is the space of observations endowed with a kernel $k$. For simplicity, we endow $(\YY, \mathcal{Y})$ with a reference $\sigma$-finite measure denoted by $\mathrm{d} y$, and we will consider only measures having a density with respect to this reference measure. The reader can think of $\YY$ being $\RR^d$, endowed with the Lebesgue measure. Let $f:\YY \times \XX \to [0,1]$ be a probability density kernel. We define the partially exchangeable mixture model with probability density kernel $f$ and mixing distribution $Q$ as 
\begin{equation} \label{eq:mix_model}
    Y_{i,j}|X_{i,j} \sim f(\cdot; X_{i,j}), \qquad X_{i,j} \big| \tilde P_{1}, \tilde P_{2} \stackrel{\iid}{\sim} \tilde P_{i}, \qquad \bigl(\tilde P_{1}, \tilde P_{2}\bigr) \sim Q,
\end{equation}
for $j \in \NN$ and $i=1,2$. We observe that $\{Y_{i,j} \}$ are partially exchangeable, and the model can be rewritten as $Y_{i,j} \big| \tilde P_{1}, \tilde P_{2}  \sim f_{\tilde P_i}$, where, for any $P \in \mathcal{P}(\mathbb{X})$, the \emph{mixture density} is defined as
\[
    f_{P}(\cdot):=\int_\XX f(\cdot; x) \diff P(x).
\]

\begin{remark}
    Consider a density $f$ proportional to $k$ on $\YY = \XX$. Then, $f_{P}$ coincides with the kernel mean embedding up to a multiplicative constant. In particular, if $f$ is the density of $\mathcal{N}(\cdot, \sigma^2)$, $f_{P}$ is the kernel mean embedding of the Gaussian kernel $k(x,y) = \exp(-(x-y)^{2}/(2\sigma^2))$. However, it is worth noting that there is not a one-to-one correspondence between mixture densities and kernel mean embedding since not all densities $f$ are symmetric functions, and their parameter space does not always coincide with the observation space, i.e., $\XX \neq \YY$. Vice versa, not all kernel functions can be rewritten as unnormalized densities, since positivity may fail, as is the case with the linear kernel. 
\end{remark}

We say that the parametric family $\{f(\cdot; x)\}_{x}$ is \emph{identifiable} if $f_{P_1} = f_{P_2}$ a.e. implies $P_1 = P_2$ for any $P_1,P_2 \in \mathcal{P}(\XX)$. This definition has been used, e.g., in  \cite{nguyen_convergence_2013}, and it is a slight modification of the original one in \cite{teicher_identifiability_1961}.

\begin{theorem} \label{mix_corr} 
    Let $f$ be a probability density kernel on $\YY \times \XX$. Then, for any kernel $k$ on $\YY \times \YY$,
    \[
        k_f(x_1,x_2) := \iint_{\YY \times \YY} k(y_1,y_2) f(y_1;x_1) f(y_2;x_2)  \diff y_1 \diff y_2, 
    \]
    is a kernel on $\XX \times \XX$ such that $\cov_{k}\bigl(f_{\tilde P_1},f_{\tilde P_2}\bigr) = \cov_{k_f}\bigl(\tilde{P}_{1},\tilde{P}_{2}\bigr)$. Moreover, if $\{f(\cdot; x)\}_{x}$ is identifiable and $k$ is characteristic, then $k_f$ is characteristic.
\end{theorem}

Consequently, under the assumptions of \cref{mix_corr}, we can extend \cref{corr_set_ker} to mixture models. Remarkably, these imply that for multivariate species sampling models, the kernel correlation between mixture models coincides with the set-wise correlation.

\begin{corollary} \label{corr_set_ker_phi}
    Let $\tilde{P}_{1}$, $\tilde{P}_{2}$ satisfy the assumptions of \cref{corr_set_ker} and let $f$ be a probability density kernel on $\YY \times \XX$ such that $\{f(\cdot; x)\}_{x}$ is identifiable. Then, for every characteristic kernel $k$ on $\YY \times \YY$ and set $A$ such that $\tilde{P}_1(A)$ and $\tilde{P}_2(A)$ are not deterministic,
    \[
        \corr_{k}\bigl(f_{\tilde P_1},f_{\tilde P_2}\bigr) = \corr_{k_f}\bigl(\tilde{P}_{1},\tilde{P}_{2}\bigr) = \corr\bigl(\tilde{P}_{1}(A),\tilde{P}_{2}(A)\bigr).
    \]
\end{corollary}

In the general case, $\corr_{k_f}\bigl(\tilde{P}_{1},\tilde{P}_{2}\bigr)$ will depend on the choice of $k$ and some kernels can be more tractable then others. When $\YY = \RR^{m}$, it is convenient to consider a translation invariant kernel $k(y_1,y_2) = \psi(y_1 - y_2)$ for some positive continuous function $\psi: \RR^{m} \rightarrow (0,+\infty)$. Since the kernel correlation is invariant with respect to scalar multiplication of the kernel, without loss of generality, we can assume $\psi(0) = 1$. For a probability distribution $\nu \in \mathcal{P}(\RR^m)$, we denote by $\hat{\nu}(x) := \int{e^{-i \langle x, z \rangle}\diff\nu(z)}$ its Fourier transform, for $x \in \RR^m$. Bochner's Theorem \citep{bochner_lectures_1959} guarantees $\psi = \hat{\nu}$ for some $\nu \in \mathcal{P}(\RR^{m})$.

\begin{proposition} \label{mix_trans_corr}
	Let $k(y_1,y_2) = \psi(y_1 - y_2)$ 
	be a translation invariant kernel for a continuous $\psi$ such that $\psi(0) = 1$, and let $\nu \in \mathcal{P}(\RR^{m})$ such that $\psi = \hat{\nu}$. Then, $k_f$ in \cref{mix_corr} is equal to
	\[
        k_f(x_1,x_2) = \int_{\RR^{m}}{\hat f(z; x_1) \overline{\hat{f}(z;x_2)}\diff\nu(z)} = \bigl\langle \hat f(\cdot; x_1),\hat f(\cdot; x_2) \bigr\rangle_{L^{2}(\nu;\mathbb{C})}.
	\]
\end{proposition}

\begin{example} \label{mix_ex}
    The Gaussian kernel $k(y_1,y_2) = \exp\left(-(y_1 - y_2)^{2}/(2\sigma^2)\right)$ is a translation invariant kernel on $\RR$ with $\psi(z) = \exp(-z^{2}/(2 \sigma^2))$, for $\sigma > 0$. We notice that $\psi(z) = \hat{\nu}(z)$, where $\nu$ is a normal distribution with mean $0$ and variance $\sigma^{-2}$. For a Gaussian mixture model with $f(\cdot;x)$ being the density of $\mathcal{N}(x,\sigma_{0}^{2})$ for some $\sigma_0>0$. In such case, $\hat f(y;x) = \exp (i x y - \sigma_{0}^{2}y^{2}/2)$. Thus 
    \[
        k_f(x_1,x_2) = \sqrt{\frac{\sigma^{2}}{2\pi}}\int_{\RR}{\exp\left(ix_1z - ix_2z - \frac{1}{2}\left(2\sigma_{0}^{2} + \sigma^2 \right)z^{2}\right)\diff z} = \sqrt{\frac{\sigma^{2}}{2\sigma_{0}^{2} + \sigma^{2}}}\exp\left(-\frac{1}{2}\frac{(x_1 - x_2)^{2}}{2 \sigma_0^2 + \sigma^{2}}\right),
    \]
    which is a Gaussian kernel with updated parameters.
\end{example}

The construction for mixture models presented in this section enables us to revisit the parametric case of \cref{par_ex} within this new framework.

\begin{example}[\cref{par_ex} -- Revisited] \label{ex:par_gau_revisited}
    We can revisit \cref{par_ex} as a mixture model as in \cref{eq:mix_model} with $f(\cdot;x)$ being the density of $\mathcal{N}(x,s^{2})$, and $\tilde{P}_{1} = \delta_{\theta_{1}}$, $\tilde{P}_{2} = \delta_{\theta_{2}}$ having joint law determined by $(\theta_{1},\theta_{2}) \sim \mathcal{N}\left(\boldsymbol{0},\tau^{2}\Sigma\right)$, where $s,\tau > 0$, $\Sigma_{11} = \Sigma_{22} = 1$ and $\Sigma_{12} = \Sigma_{21} = \rho \in [-1,1]$.
    
    We take the Gaussian kernel $k(y_1,y_2) = \exp\left(-(y_1 - y_2)^{2}/(2\sigma^2)\right)$ on $\YY \times \YY$. By \cref{mix_ex} the kernel $k_f$ over $\XX \times \XX$ is $k_{f}(x_1,x_2) =$ $\sqrt{\sigma^{2}/(2s^{2} + \sigma^{2}))}\exp(-(x_1-x_2)^{2}/(2(2s^{2} + \sigma^{2}))$. The computation of $\corr_{k_f}\bigr(\tilde{P}_{1},\tilde{P}_{2}\bigl)$ only involves Gaussian integrals and can be done with 
    (SM24) in the Supplementary Material. We obtain, thanks to \cref{corr_set_ker_phi},
\begin{equation} \label{corr_k_par_ex}
        \corr_{k}\bigr(f_{\tilde{P}_{1}},f_{\tilde{P}_{2}}\bigl)=\corr_{k_f}\bigr(\tilde{P}_{1},\tilde{P}_{2}\bigl) = \frac{\sqrt{\frac{\sigma^{2}}{2\tau^{2}(1 - \rho) + 2s^{2} + \sigma^{2}}} - \sqrt{\frac{\sigma^{2}}{2\tau^{2} + 2s^{2} + \sigma^{2}}}}{\sqrt{\frac{\sigma^{2}}{2s^{2} + \sigma^{2}}} - \sqrt{\frac{\sigma^{2}}{2\tau^{2} + 2s^{2} + \sigma^{2}}}}.
    \end{equation}
    
    Since $k_{f}$ is $c_{0}$-universal, and since $\tilde{P}_{i} = \delta_{\theta_i}$ are discrete a.s. with atomless $\E\bigl[\tilde{P}_{i}\bigr] = \mathcal{N}(0,\tau^{2})$, we can apply \cref{index_equal1_prior}. Hence, it holds that $\corr_{k}\bigl(f_{\tilde{P}_{1}},f_{\tilde{P}_{2}}\bigr) = 1$ if and only if $\tilde{P}_{1} = \tilde{P}_{2}$ a.s.. A similar strategy can be extended to other parametric models as well.
\end{example}

The previous example shows that kernel correlation can detect exchangeability for parametric models by reinterpreting them as mixture models over a vector of a.s. discrete random probabilities that do not belong to the class of multivariate species sampling models. The same rewriting shows that it detects exchangeability for mixture models independently of their dependence structure, thus extending \cref{index_equal1_prior} to this setting.

\begin{corollary} \label{exchangeability_mixtures}
    Let $\tilde{P}_{i}$ be an a.s. discrete random probability such that $\E\bigl[\tilde{P}_{i}\bigr]$ is atomless, for $i=1,2$, and let $f$ be a probability density kernel on $\YY \times \XX$. If $k_f$ is $c_0$-universal, then $\corr_{k}\bigr(f_{\tilde{P}_{1}},f_{\tilde{P}_{2}}\bigl) = 1$ if and only if $f_{\tilde{P}_{1}} = f_{\tilde{P}_{2}}$ a.s..
\end{corollary}

\section{Hierarchical Dirichlet Process} \label{sec:hdp}

In this section, we turn to one of the most popular models for partially exchangeable data, the hierarchical Dirichlet Process \citep{teh_hierarchical_2006}, and we study the behaviour of the kernel correlation both a priori and a posteriori. Remarkably, we are able to recover the same rate of convergence of the parametric model in \cref{par_ex} at the cost of adding an extra assumption on the non-degeneracy of the data sequence with respect to the kernel. The set-wise correlation is not the best choice in this context: the verification of the assumption depends heavily on the choice of the set, and there is no sensible way to choose it before looking at the specific dataset.

A random probability $\tilde{P} \sim \mathrm{DP}(c, P_0)$ follows a Dirichlet Process with concentration $c > 0$ and base measure $P_{0} \in \mathcal{P}(\XX)$ if $\bigl(\tilde{P}(A_{1}),\tilde{P}(A_{2}),\dots,\tilde{P}(A_{m})\bigr) \sim \mathrm{Dir}(c P_{0}(A_{1}),c P_{0}(A_{2}),\dots,c P_{0}(A_{m}))$ for any partition $\{A_{1},A_{2},\dots,A_{m}\}$ of $\XX$, where $\mathrm{Dir}$ indicates the finite-dimensional Dirichlet distribution. Starting from the definition of a Dirichlet Process, the hierarchical Dirichlet Process (hDP) of \cite{teh_hierarchical_2006} provides a natural way of building a dependent nonparametric prior for a vector $\bigl(\tilde P_1, \tilde P_2\bigr)$, which can then be used to define partially exchangeable models as in \eqref{part_ex}. Specifically, $\bigl(\tilde{P}_{1},\tilde{P}_{2}\bigr) \sim \mathrm{hDP}(c,c_0,P_0)$ for $c,c_{0}>0$ and $P_{0}$ an atomless probability measure on $\XX$, if
\begin{equation}  \label{hDP}
    \tilde{P}_{1},\tilde{P}_{2} \big| \tilde{P}_{0} \stackrel{\iid}{\sim} \mathrm{DP}\bigl(c,\tilde{P}_{0}\bigr), \qquad \tilde{P}_{0} \sim \mathrm{DP}(c_{0},P_{0}).
\end{equation}
The calculations in Example~SM5 in the Supplementary material show that the the hDP a priori satisfies \cref{var_set_ker} with $\lambda_{1} = \lambda_{2} = (1 + c +c_0)(1+c)^{-1}(1+c_0)^{-1}$ and \cref{cov_set_ker} with $\eta = (1 + c_0)^{-1}$. Thus by \cref{corr_set_ker},
\[
    \corr_{k}\bigl(\tilde{P}_{1},\tilde{P}_{2}\bigr) = \corr\bigl(\tilde{P}_{1}(A),\tilde{P}_{2}(A)\bigr) = \frac{1+c}{1+c+c_{0}},
\]
for any injective kernel $k$ and any measurable set $A$ such that $P_{0}(A) \notin \{0,1\}$. The expression of $\corr\bigl(\tilde{P}_{1}(A),\tilde{P}_{2}(A)\bigr)$ had already been derived in \cite{camerlenghi_distribution_2019}: our contribution is to show that it coincides with the kernel correlation of \emph{any} characteristic kernel.

We now study the rate of convergence of the kernel correlation of the hDP a posteriori to zero as the number of observations diverges. First, we define a notion of non-degeneracy for a sequence, which in our setting will be the sequence of observations within each group.
\begin{definition} \label{def_non_deg}
    We say that a sequence $(x_{j})_{j \in \NN}$ is non-degenerate with respect to a kernel $k$ if
    \begin{equation} \label{eq:non_deg_cond}
        \liminf_{n \rightarrow +\infty} \frac{1}{n^{2}}\sum_{j=1}^{n}{\sum_{h=1}^{n}{d^{2}_{k}(x_{j},x_{h})}} = \liminf_{n \rightarrow +\infty} \int_{\XX}\int_{\XX}{d_{k}^{2}(x,y)\diff \hat P _{n}(x)\diff\hat P_{n}(y)} > 0,
    \end{equation}
    where $\hat P_{n} = n^{-1}\sum_{j=1}^{n}{\delta_{x_{j}}}$ and $d_k^2(x,y)= k(x,x) - 2k(x,y) + k(y,y)$. 
\end{definition}

We note that $d_{k}$ is a pseudo-metric, and it is a distance if $k$ is injective. Its role in the study of the kernel correlation can be understood thanks to (SM2) in the Supplementary Material.

The notion of non-degeneracy depends heavily on the choice of kernel. We gain a better understanding by considering $(x_{j})_{j \in \NN}$ a realisation of an infinitely exchangeable sequence $X_{j}|\tilde{P} \stackrel{\iid}{\sim} \tilde{P}$, with $\tilde{P} \sim Q$ its de Finetti measure. By de Finetti's Representation Theorem \citep{definetti_prevision_1937}, the empirical measure converges weakly to $\tilde{P}$ a.s. as the number of observations $n$ diverges. In this setting, 
\begin{equation} \label{exch_deg}
	\liminf_{n \rightarrow +\infty} \iint{d_{k}^{2}(x,y)\diff\hat P_{n}(x)\diff\hat P_{n}(y)} = \iint{d_{k}^{2}(x,y)\diff\tilde{P}(x)\diff\tilde{P}(y)} \qquad \as{}
\end{equation}
If the kernel is injective, then $d_k$ is a distance. Hence, the right-hand side is zero if and only if $\tilde{P} = \delta_{Z}$ a.s. for some random variable $Z$ on $\XX$. This means that the sequence $(X_{j})_{j \in \NN}$ is a.s. constant, which is arguably an intuitive notion of degeneracy. However, when $k(x,y) = \mathbbm{1}_{A}(x)\mathbbm{1}_{A}(y)$ for some $A$, the right hand side in \cref{exch_deg} is zero whenever $\tilde{P}(A) \in \{0,1\}$ a.s.. Summarising, the non-degeneracy assumption is not restrictive when the kernel $k$ is injective, and in such cases, it does not depend on $k$. For instance, linear, Gaussian, and Laplace kernels induce the same notion of non-degeneracy. In contrast, for the set-wise kernel, the notion of non-degeneracy depends heavily on the choice of the set, which cannot be chosen before seeing the data.

The following theorem states that if we assume that the sequences of observables are non-degenerate for both groups, there is a regular version of the posterior, derived in \cite{camerlenghi_distribution_2019} and reported in (SM9) of the Supplementary Material, for which we can recover the parametric rate of convergence. For observed data $\boldsymbol{x}^{(n_{1},n_{2})} = \bigl((x_{1,j})_{j=1}^{n_{1}},(x_{2,j})_{j=1}^{n_{2}}\bigr)$, we denote by $\mathcal{L}\bigl(\tilde{P}_{1},\tilde{P}_{2} \big| \boldsymbol{X}^{(n_1,n_2)} = \boldsymbol{x}^{(n_{1},n_{2})}\bigr)$ the point-wise evaluation of the corresponding Markov kernel. 

\begin{theorem} \label{convergence_corr_hDP}
    Consider a partially exchangeable model as in \eqref{part_ex} for $\bigl(\tilde{P}_{1},\tilde{P}_{2}\bigr) \sim \mathrm{hDP}(c,c_0,P_0)$ for $c,c_{0}>0$ and $P_{0}$ an atomless probability measure on $\XX$. If $(x_{i,j})_{j \in \NN}$ is a non-degenerate sequence with respect to a kernel $k$ for $i = 1,2$, then as $\max(n_{1}, n_{2}) \rightarrow +\infty$,
	\[
	   \corr_{k}\bigl(\tilde{P}_{1},\tilde{P}_{2} |\boldsymbol{X}^{(n_1,n_2)} = \boldsymbol{x}^{(n_{1},n_{2})}\bigr) = \mathcal{O}\left(\frac{1}{\sqrt{n_{1}n_{2}}}\right).
	\]
\end{theorem}

Our results on the hierarchical Dirichlet process rely on a fine understanding of its posterior structure, which we revise in the Supplementary Material. These results are also helpful to design the simulations in \cref{sec:numerics}.

\section{Numerical Simulations} \label{sec:numerics}

This section aims to numerically validate the convergence rate obtained in \cref{convergence_corr_hDP} and analyze the stability of the kernel correlation with respect to the choice of kernel. We devise two different strategies to approximate the kernel correlation a posteriori. The \emph{sampling-based} method utilizes the estimator presented in \cref{sec:estimator} and can be applied whenever it is possible to generate $2\times 2$ samples from the model a posteriori. The \emph{analytics-based} method is built on the quasi-conjugacy property of the hDP \citep{teh_hierarchical_2006,camerlenghi_distribution_2019}. We empirically show that the \emph{ad hoc} construction of the analytics-based estimator has a lower variance, underlining the importance of analytic calculations when possible.


\subsection{Sampling-based vs Analytics-based Method for the Hierarchical Dirichlet Process}

Both the sampling-based method and the analytic-based method rely on the quasi-conjugacy property of the augmented hDP model we obtain after the introduction of a specific sequence of latent random variables $\boldsymbol{T}^{(n_{1},n_{2})}$, commonly referred to as \emph{tables} in the restaurant franchise metaphor \citep{teh_hierarchical_2006, camerlenghi_distribution_2019,catalano_unified_2023}. To generate a $2 \times 2$ sample from the posterior distribution, as in the sampling-based algorithm, we generate $\boldsymbol{T}^{(n_{1},n_{2})}$ conditionally on $\boldsymbol{X}^{(n_{1},n_{2})}$ and apply the predictive distribution of the quasi-conjugate scheme. Similarly, in the analytics-based algorithm, we find the exact expression of the kernel correlation a posteriori conditionally on $\boldsymbol{T}^{(n_{1},n_{2})}$, and then we generate $R$ copies of $\boldsymbol{T}^{(n_{1},n_{2})}|\boldsymbol{X}^{(n_{1},n_{2})}$ to marginalize the tables out. We refer to Section~SM3 of the Supplementary Material for a detailed description of the two algorithms.

\begin{figure}
    \includegraphics[width=0.495\linewidth]{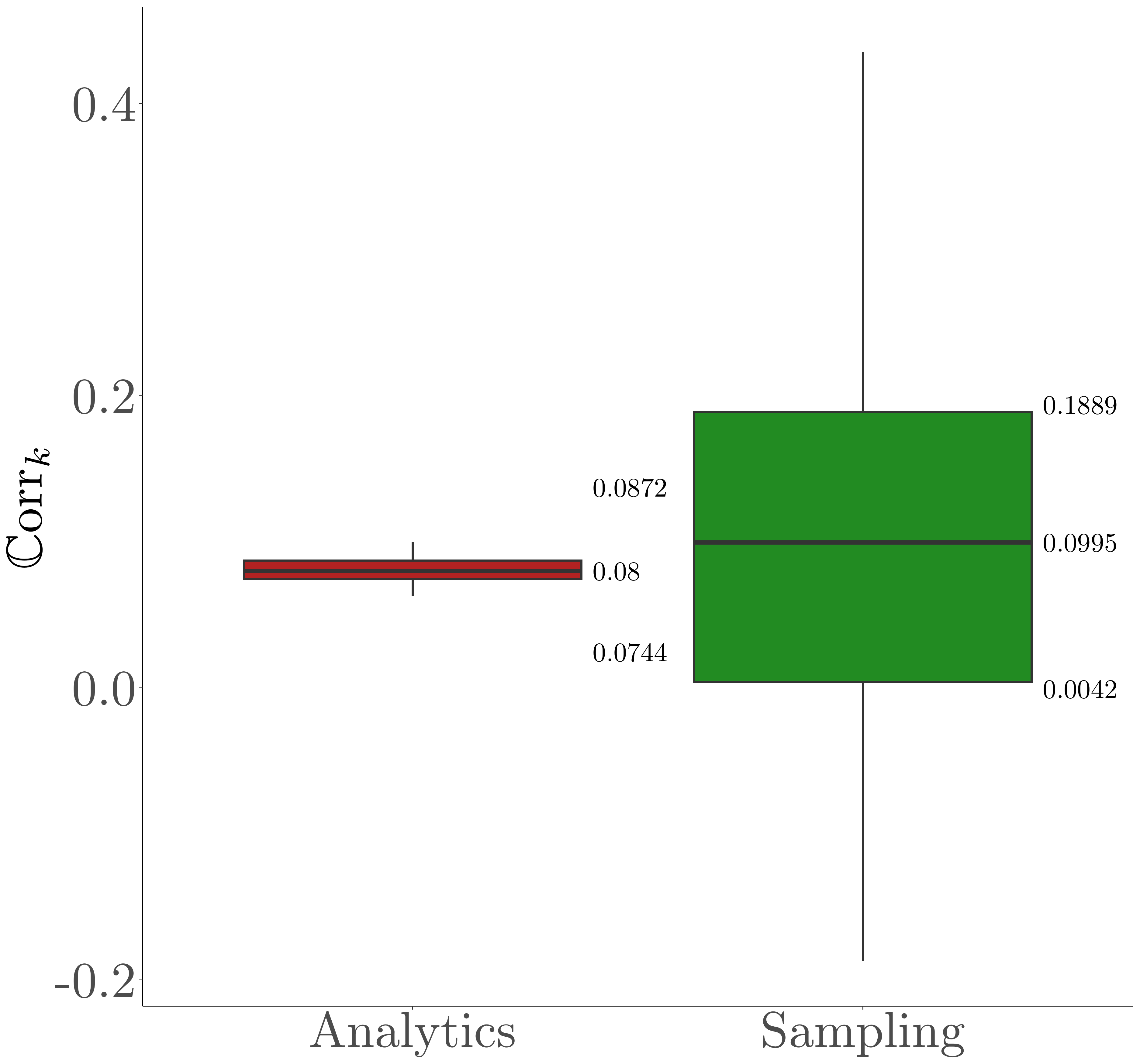} 
    \includegraphics[width=0.495\linewidth]{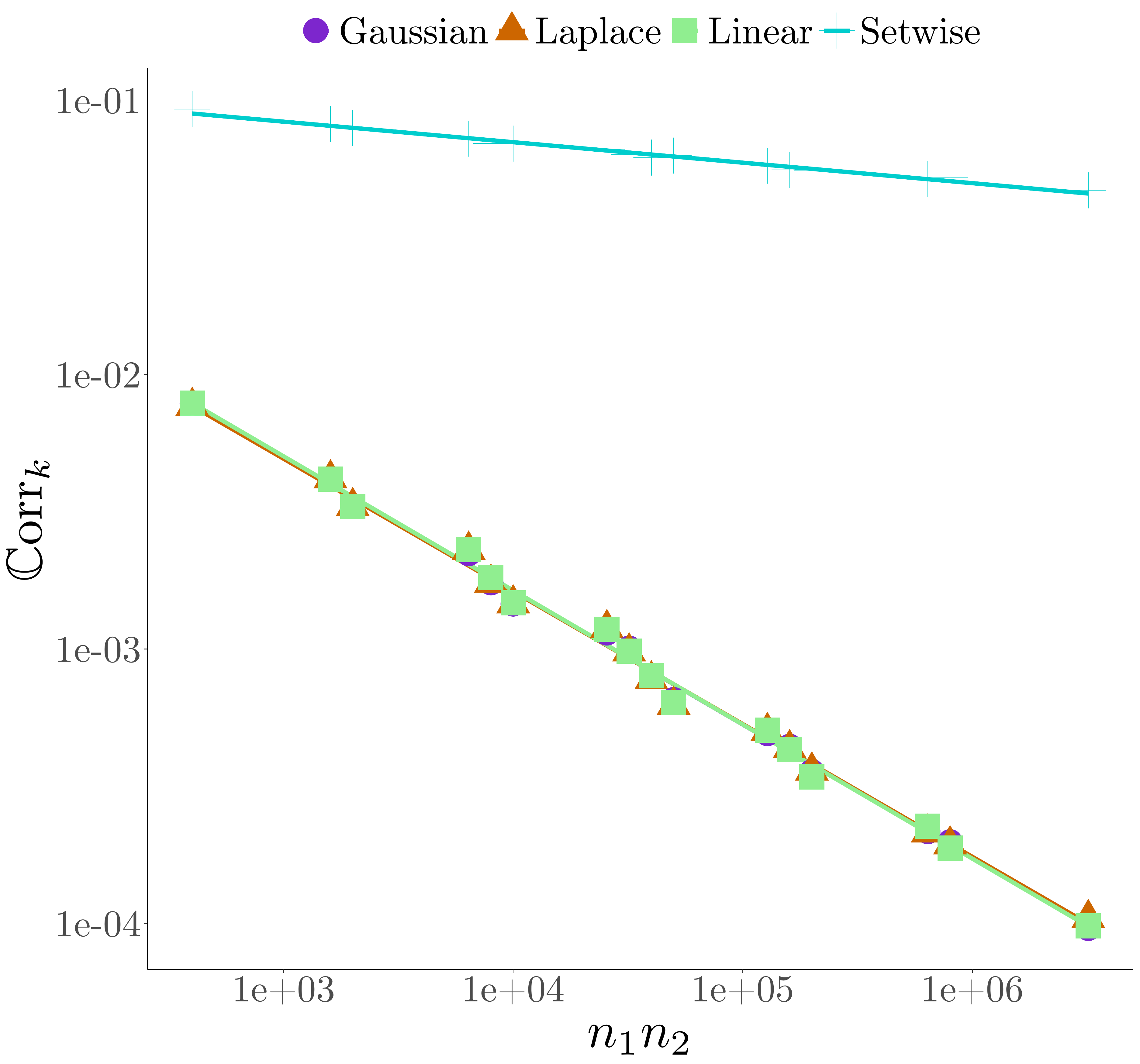}
    \caption{\textbf{Left:} Box plots for \emph{Analytics-based} and \emph{Sampling-based} algorithms for $100$ realizations of the estimate of the posterior kernel correlation for $\bigl(\tilde{P}_{1},\tilde{P}_{2}\bigr) \sim \mathrm{hDP}\bigl(c=1,c_0=1,P_{0} = \mathrm{Unif}_{[0,1]}\bigr)$ when we condition on $n_{1} = n_{2} = 10$ data points sampled from the model and use a Gaussian kernel with $\sigma = 1$. \textbf{Right:} Log-log plot of the kernel correlation a posteriori as a function of $n_{1}n_{2}$ with $n_{1} = 4^{i}$ and $n_{2} = 5^{j}$ for $i,j = 2,3,4,5$ for $\bigl(\tilde{P}_{1},\tilde{P}_{2}\bigr) \sim \mathrm{hDP}\bigl(c=1,c_0=1,P_{0} = \mathrm{Unif}_{[0,1]}\bigr)$ and data points sampled from the model. The kernel correlation is computed for different kernels: Gaussian with $\sigma = 1$, Laplace with $\beta = 1$, set-wise with $A = [0,0.95]$, and linear.}
    \label{fig:kernels}
\end{figure}

We test the two algorithms by computing the posterior kernel correlation for the Gaussian kernel with $\sigma = 1$ coming from an hDP prior \eqref{hDP} with $c_{0} = c = 1$ and $P_{0} \sim \mathrm{Unif}_{[0,1]}$ when the data are $n_{1} = n_{2} = 10$ observations from the model. We run both algorithms on the same data $100$ times. In the analytics-based algorithm, we run the Gibbs sampler $R=10$ times to approximate the expectations in the law of total covariance. In contrast, the sampling-based algorithm uses the empirical covariance estimator, as described in \cref{prop_cov_unbiased}, with $M = 10,000$ independent samples. The box plots on the left of \cref{fig:kernels} show the estimated posterior kernel covariance for the two methods, annotating the first, second, and third quartiles. We observe that the sampling-based method exhibits significantly more variability than the analytics-based method. A possible explanation is that the former, despite the vast number of generated samples, does not account for the additional knowledge about the posterior, whereas the latter is built \emph{ad hoc} for this model. There is an evident trade-off between generality and variability: when possible, analytic computations can reduce the variability.

\subsection{Convergence Rate}

Our aim is to empirically recover the convergence rate $\mathcal{O}\bigl((n_{1}n_{2})^{-1/2}\bigr)$ of \cref{convergence_corr_hDP}. We recover this for several injective kernels, while empirically showing that the set-wise kernel has a slower convergence rate. 

We compute the posterior kernel correlation with Gaussian, Laplace, linear, and set-wise kernels for $\bigl(\tilde{P}_{1},\tilde{P}_{2}\bigr) \sim \mathrm{hDP}\bigl(c=1,c_0=1,P_{0} = \mathrm{Unif}_{[0,1]}\bigr)$, when the observations are sampled from the model for $n_{1} = 4^{i}$ and $n_{2} = 5^{j}$ for $i,j = 2,3,4,5$. We choose these values to have a grid of different values for $n_{1}n_{2}$. The Gaussian and the Laplace kernel have their parameters set to $\sigma = \beta = 1$, while the set-wise kernel is taken for $A = [0,0.95]$. We use the analytics-based method because of its lower variance, with $R = 1000$. The log-log plot on the right of \cref{fig:kernels} reports the value of the estimated posterior kernel correlation for the different values of $n_{1}n_{2}$ and different kernel choices. Moreover, we report the regression line of the log-kernel correlation as a function of $\log(n_{1}n_{2})$. The estimated value of the slope is approximately $-1/2$ for Gaussian, Laplace, and linear kernels, as it is noticeable in \cref{fig:kernels}. This result, in the logarithmic domain, confirms our theoretical findings since the sequence of observables is non-degenerate for injective kernels. However, the set-wise correlation in \cref{convergence_corr_hDP} fails to capture this rate of convergence and has a much slower decrease due to the choice of $A = [0,0.95]$. Since the marginal distribution of the observables is $\mathrm{Unif}_{[0,1]}$, almost all the values will be contained in $A$, making the sequence \emph{almost} degenerate for the chosen set.

\subsection{Stability of the Kernel}

We now investigate the stability of the posterior kernel correlation with respect to the hyperparameters of the kernel. As expected, the choice of hyperparameters does not have a significant impact on the value of the kernel correlation for Gaussian or Laplace kernels, while it dramatically changes the value of the set-wise correlation. 

As shown in \cref{post_ex}, the set-wise kernel can depend heavily on the choice of the set $A$, and we expect this to be the case for the posterior since it has a similar structure. We investigate the same question for the parameters of Gaussian and Laplace kernels by considering a similar setting to the one in the previous section. We compute the posterior kernel correlation with Gaussian, Laplace, linear, and set-wise kernels for $\bigl(\tilde{P}_{1},\tilde{P}_{2}\bigr) \sim \mathrm{hDP}\bigl(c=1,c_0=1,P_{0} = \mathrm{Unif}_{[0,1]}\bigr)$, when the observations are sampled from the model for $n_{1} = n_{2} \in \{0,10,100,1000\}$. We choose the parameters for Gaussian, Laplace, and set-wise kernels to take one of three values as follows: $\sigma, \beta \in \{10^{-3},1,10^{3}\}$, and $A \in \{[0,.1],[0,.5],[0,.9]\}$. We use the analytics-based algorithm with $R = 1000$. \cref{tab:stability} shows the values of the kernel correlation for different kernels and parameters as the sizes of the observables, $n_{1} $ and $ n_{2}$, increase. Both Gaussian and Laplace kernels show little variability for any sample size. In contrast, for different choices of $A$, the variability of the set-wise correlation rises dramatically as $n_{1},n_{2}$ increase, suggesting a lack of robustness of the index a posteriori. This is especially problematic as it is challenging to elicit $A$ without prior knowledge of the data since, for a sensible measurement, $A$ must not contain too many or too few observations.

\begin{table}[b!]
\centering
\caption{Value of kernel correlation for Gaussian, Laplace, and set-wise kernels for different choices of parameters and sample sizes $n_{1},n_{2}$.} \label{tab:stability}
\resizebox{\textwidth}{!}{
\begin{tabular}{l c c c c c c c c c}
    \toprule
    & \multicolumn{3}{c}{\textbf{Gaussian}, $\sigma > 0$} & \multicolumn{3}{c}{\textbf{Laplace}, $\beta > 0$} & \multicolumn{3}{c}{\textbf{Set-wise}, $A = [0,b]$} \\
    \cmidrule(rl){2-4} \cmidrule(rl){5-7} \cmidrule(rl){8-10} 
    $n_{1},n_{2}$ & $\sigma = 10^{-3}$ & $\sigma =  1$ & $\sigma = 10^{3}$ & $\beta = 10^{-3}$ & $\beta = 1$ & $\beta = 10^{3}$ & $b = 0.1$ & $b = 0.5$ & $b = 0.9$ \\
    \midrule
    $0$ & $6.7 \cdot 10^{-1}$ & $6.7 \cdot 10^{-1}$ & $6.7 \cdot 10^{-1}$ & $6.7 \cdot 10^{-1}$ & $6.7 \cdot 10^{-1}$ & $6.7 \cdot 10^{-1}$ & $6.7 \cdot 10^{-1}$ & $6.7 \cdot 10^{-1}$ & $6.7 \cdot 10^{-1}$ \\
    $10$ & $1.8 \cdot 10^{-2}$ & $1.7 \cdot 10^{-2}$ & $1.7 \cdot 10^{-2}$ & $1.8 \cdot 10^{-2}$ & $1.7 \cdot 10^{-2}$ & $1.7 \cdot 10^{-2}$ & $1.2 \cdot 10^{-1}$ & $1.6 \cdot 10^{-2}$ & $1.2 \cdot 10^{-1}$ \\
    $100$ & $1.8 \cdot 10^{-3}$ & $1.7 \cdot 10^{-3}$ & $1.7 \cdot 10^{-3}$ & $1.8 \cdot 10^{-3}$ & $1.6 \cdot 10^{-3}$ & $1.7 \cdot 10^{-3}$ & $6.9 \cdot 10^{-2}$ & $1.6 \cdot 10^{-3}$ & $6.9 \cdot 10^{-2}$ \\
    $1000$ & $1.9 \cdot 10^{-4}$ & $1.7 \cdot 10^{-4}$ & $1.7 \cdot 10^{-4}$ & $1.7 \cdot 10^{-4}$ & $1.8 \cdot 10^{-4}$ & $1.7 \cdot 10^{-4}$ & $5.0 \cdot 10^{-2}$ & $1.8 \cdot 10^{-4}$ & $5.0 \cdot 10^{-2}$ \\
    \bottomrule
\end{tabular}}
\end{table}

\section{Model Comparison} \label{sec:models}

Having established an RKHS-based index for measuring partial exchangeability for both parametric and nonparametric models, we can now use it to calibrate prior parameters across different models to match the same index value. We consider the Gaussian parametric model in \cref{par_ex} and the hierarchical Dirichlet Process in \eqref{hDP} with $P_{0} = \mathcal{N}(0,t^{2})$ with $t^{2} = s^{2} + \tau^{2}$, so that both models share the same prior predictive distribution. We set the remaining parameters to ensure the same value of marginal kernel correlation a priori and study the implications for posterior inference.

\begin{figure*}[t]
    \includegraphics[width=0.320\linewidth]{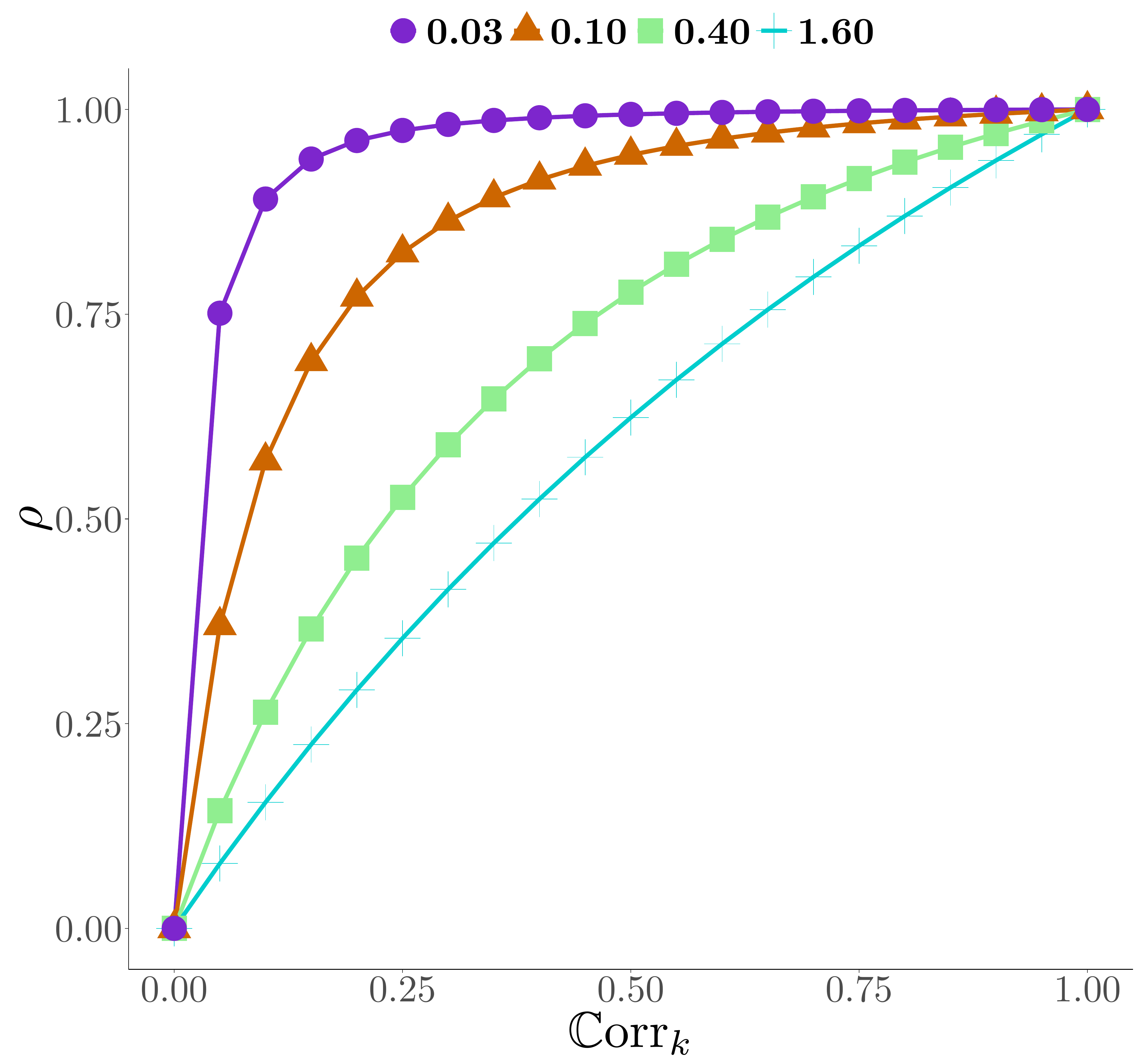}
    \includegraphics[width=0.320\linewidth]{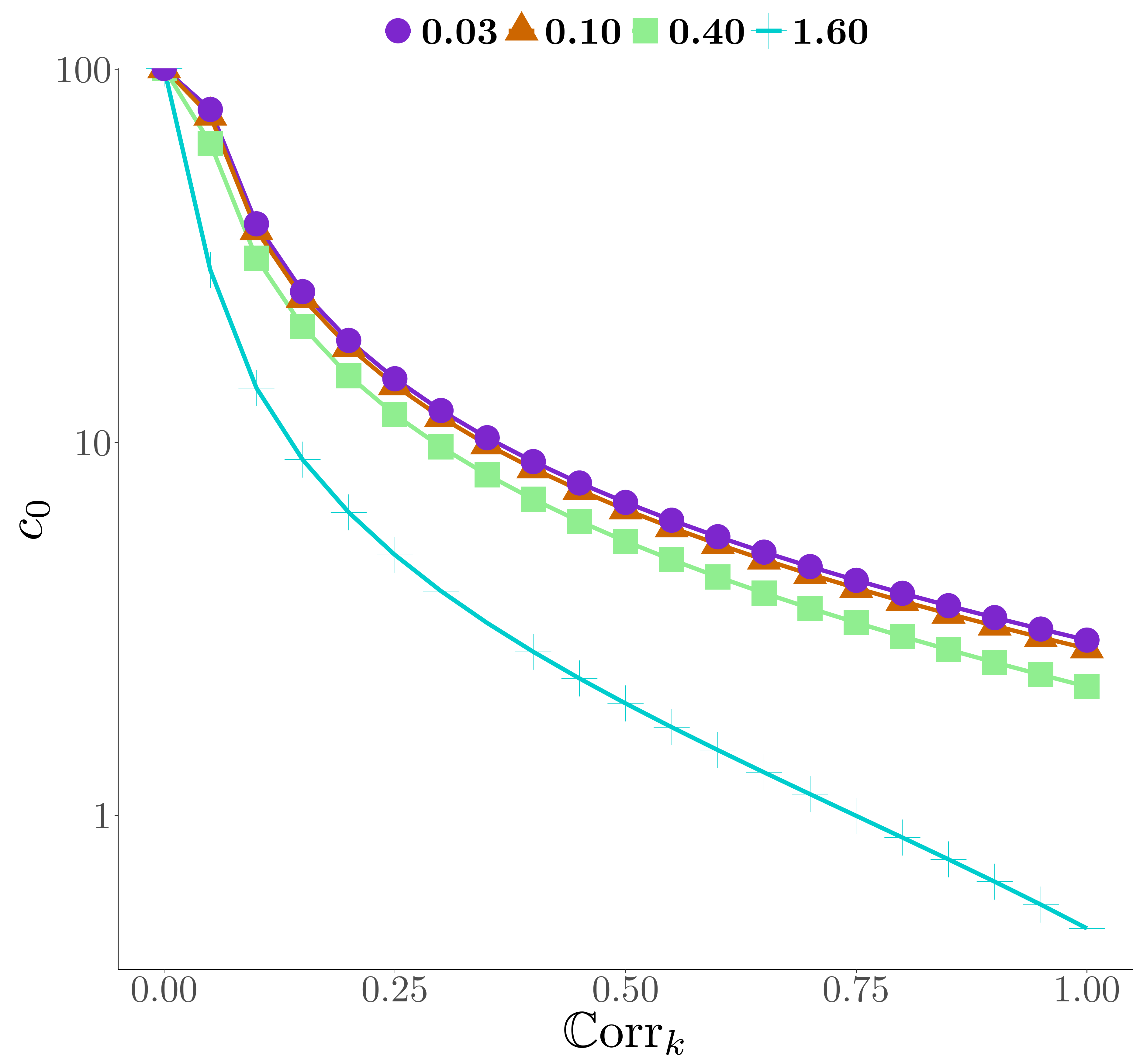}
    \includegraphics[width=0.320\linewidth]{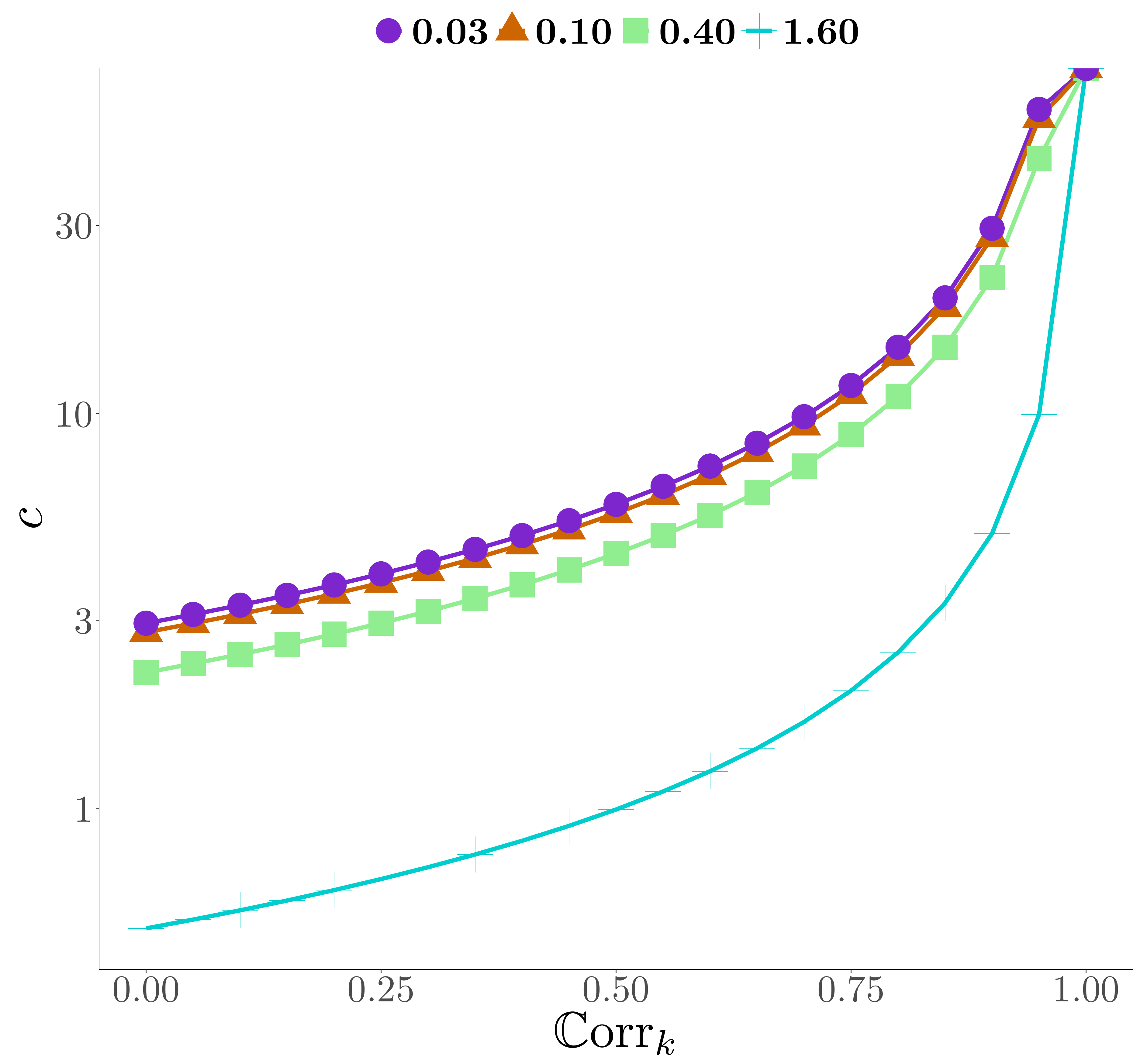}
    \caption{Values of different parameters for the Gaussian model and the hDP for $v = 1/4$, $t^{2} = 2$, and different values of $\sigma$ as the kernel correlation varies. \textbf{Left:} Value of $\rho$ for the Gaussian model. Centre: Value of $c_{0}$ for the hDP. \textbf{Right:} Value of $c$ for the hDP.}
    \label{fig:corr_params}
\end{figure*}

As shown in Section SM5 in the Supplementary Material, given a fixed marginal variance of the observables ($\var(X_i) = t^{2}$) and fixed kernel variances a priori ($\var_{k}(\tilde P_i) = v$), we can compute the model parameters ($s,\tau,\rho$ for the Gaussian model, $c_{0},c$ for the hDP) for any value of the kernel correlation $\corr_k(\tilde{P}_1, \tilde{P}_2)$, provided that the parameter $\sigma$ of the Gaussian kernel satisfies $\sigma < \sigma^{*} := \sqrt{2}t/\sqrt{1/(1-v)^{2} - 1}$. In \cref{fig:corr_params}, we study how the values of the parameters change as the kernel correlation varies for $t^{2} = 2$, $v = 1/4$, and $\sigma^{2} = (\sigma^{*}/4^{i})^{2}/2$ for $i = 0,1,2,3$. For the Gaussian model, the value of $\rho$ (on the left) increases as the kernel correlation increases, with $\rho = 0$ corresponding to independence between the two groups and a null kernel correlation, and $\rho = 1$ corresponding to full exchangeability between the two groups and a kernel correlation equal to $1$. For the hDP, we notice that $c_{0}$ (in the middle) diverges as the kernel correlation approaches $0$; this can be explained by the fact that we obtain almost sure equality of the two random probability measures whenever $\tilde{P}_{0}$ converges to $P_{0}$ a.s.. Conversely, $c$ (on the right) diverges as the kernel correlation approaches $1$ since both $\tilde{P}_{1}$ and $\tilde{P}_{2}$ converge to $\tilde{P}_{0}$ a.s.. We observe that a high value of $\sigma$ yields the same correlation value for smaller values of $c$ and $c_{0}$, resulting in greater numerical stability. Consequently, we fix $\sigma = \sigma^{*}/\sqrt{2}$ for the rest of the section.

\begin{table}[t]
\centering
\caption{Predictive distributions for different values of the kernel correlation for the Gaussian case and the hDP case, conditionally on $(X_{1,j})_{j=1}^{n_{1}}$ and $(X_{2,j})_{j=1}^{n_{2}}$ as in \cref{eq:post_smpl} for $n_{1} = 200$ and $n_{2} = 5$.} \label{tab:models}
\begin{tabular}{c c c}
    \toprule
    $\corr_{k}$ & Parametric (Gaussian) & Nonparametric (hDP) \\
    \midrule
    \raisebox{0.25\textwidth}{0.01} & \includegraphics[width=0.25\textwidth]{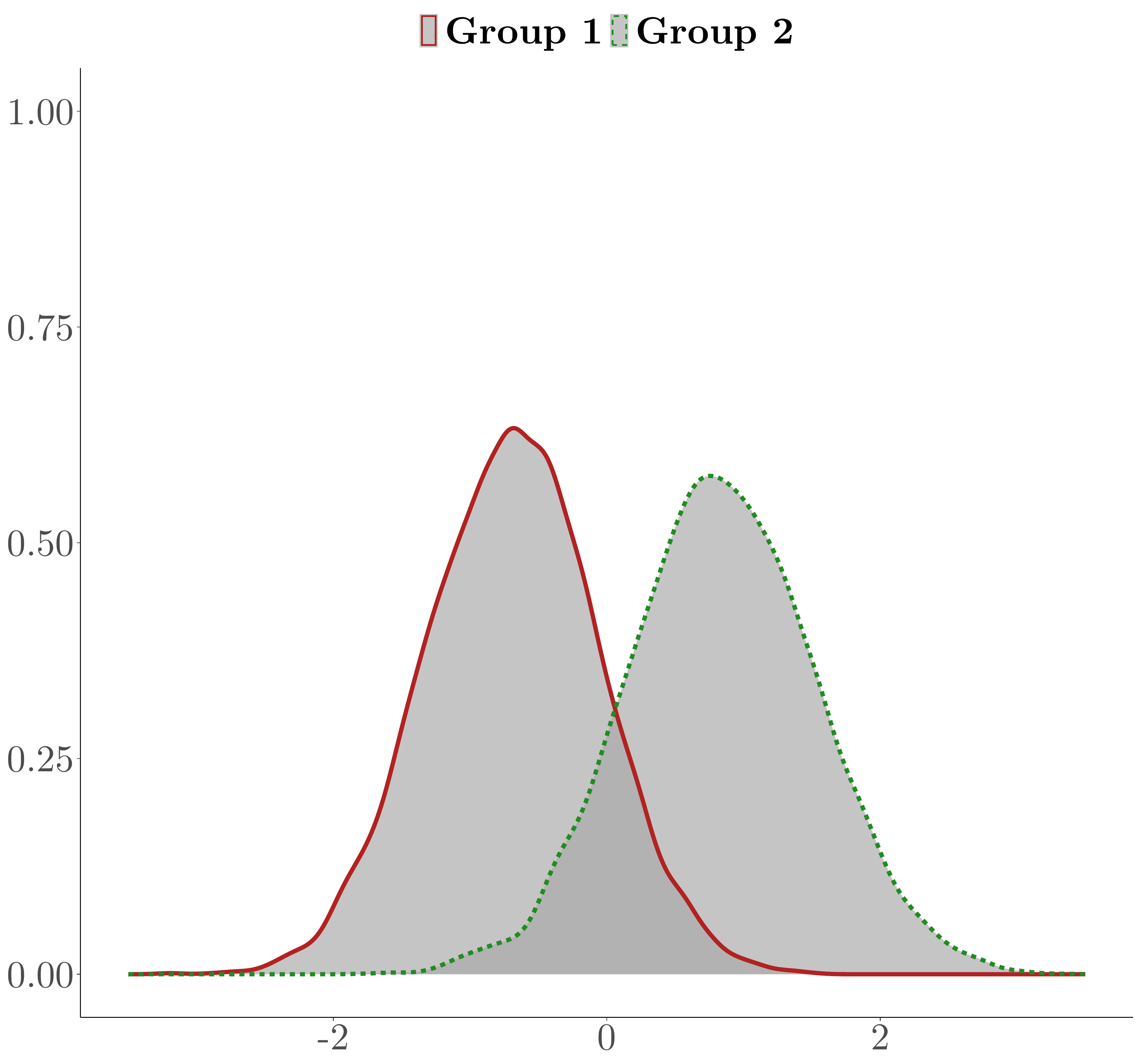} &
    \includegraphics[width=0.25\textwidth]{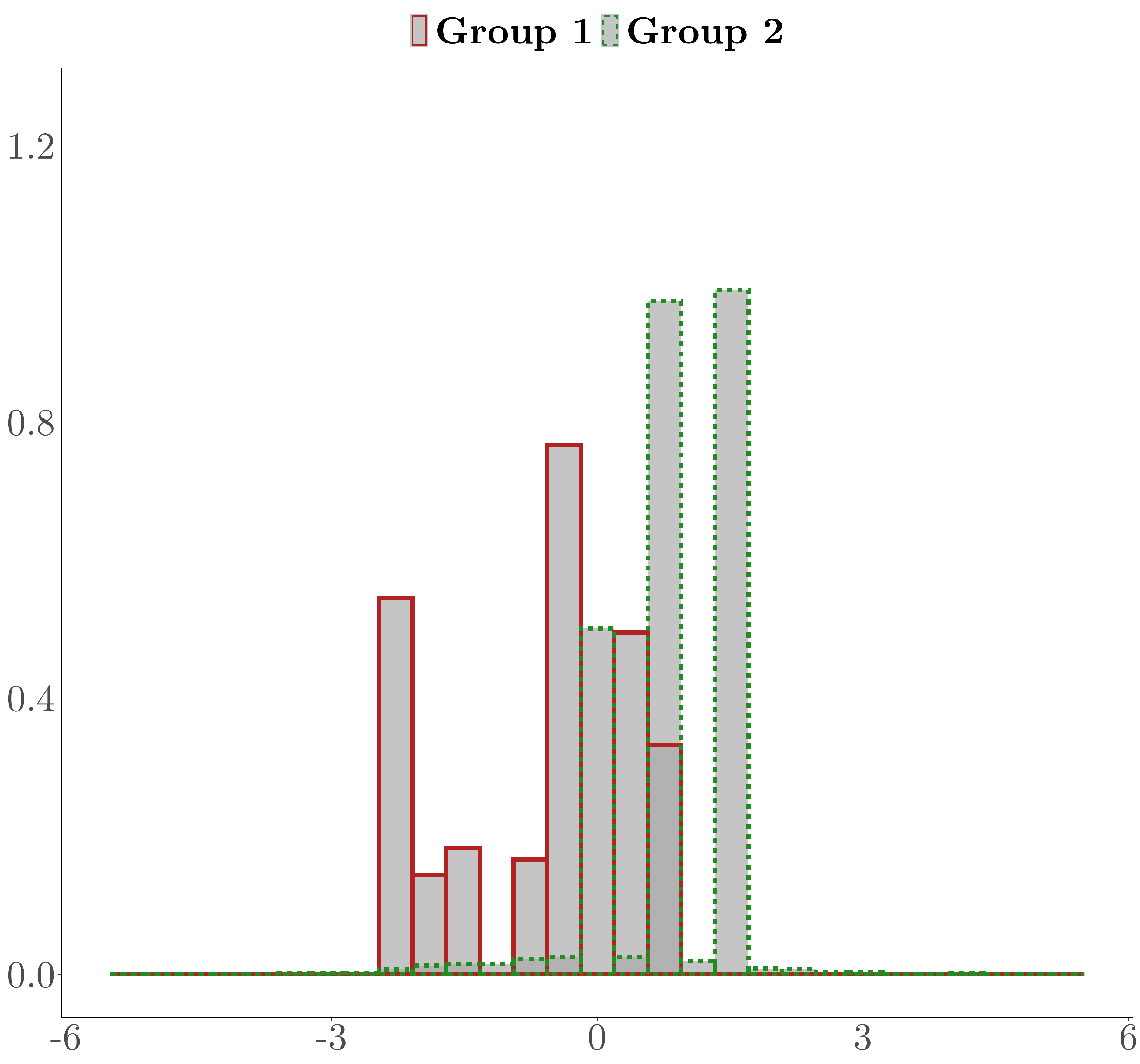} \\
    \raisebox{0.25\textwidth}{0.50} & 
    \includegraphics[width=0.25\textwidth]{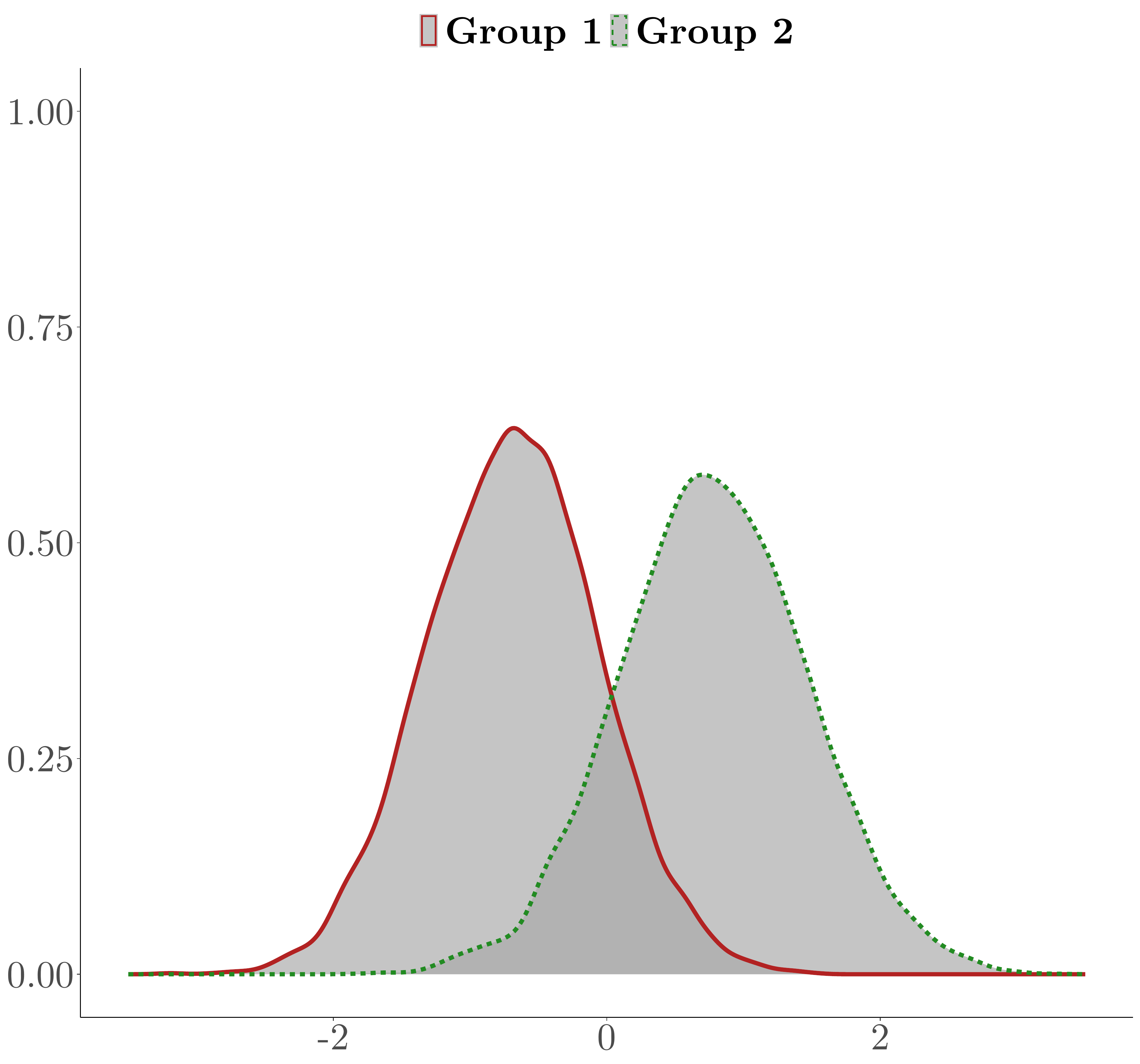} &
    \includegraphics[width=0.25\textwidth]{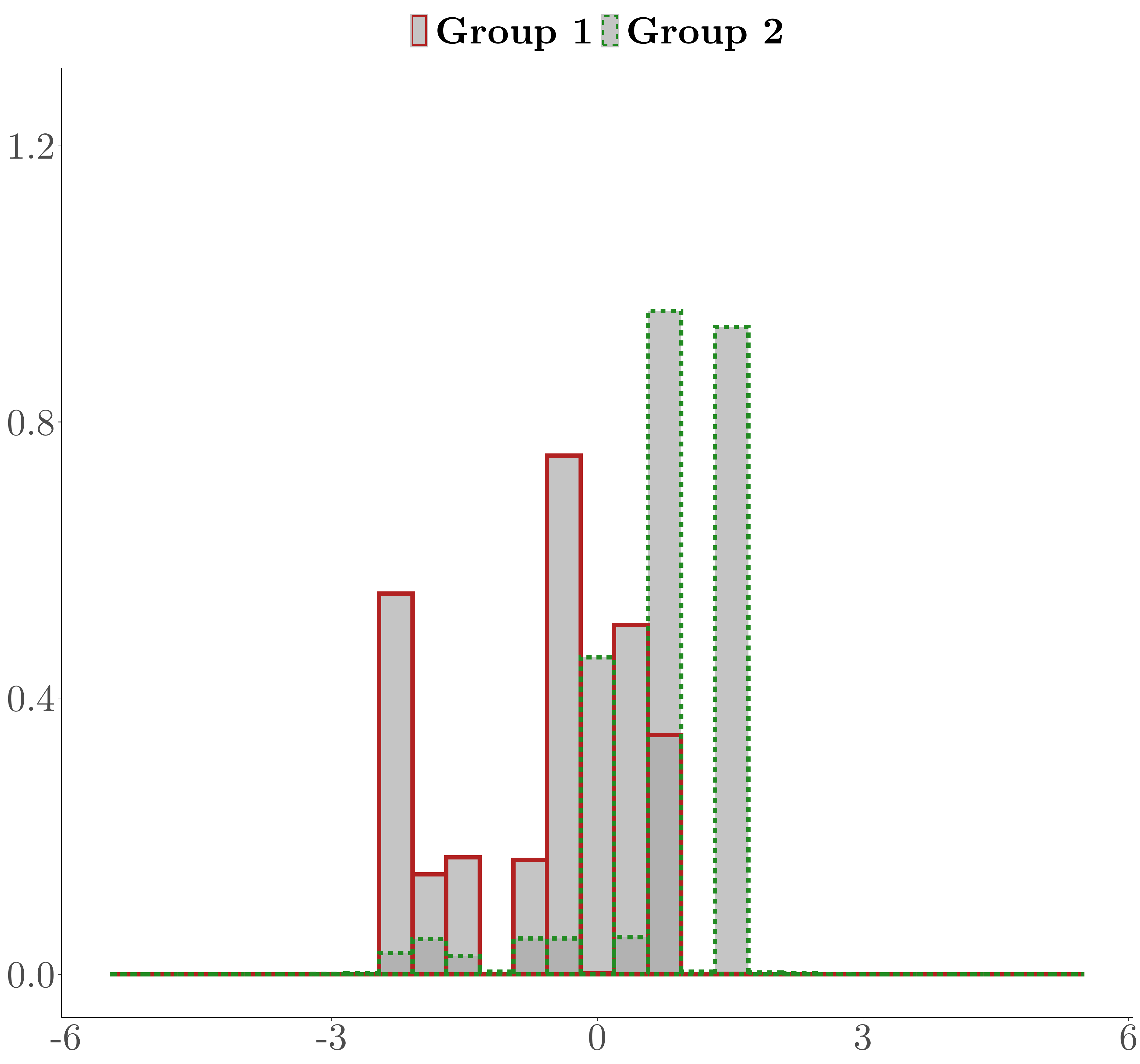} \\
    \raisebox{0.25\textwidth}{0.99} & \includegraphics[width=0.25\textwidth]{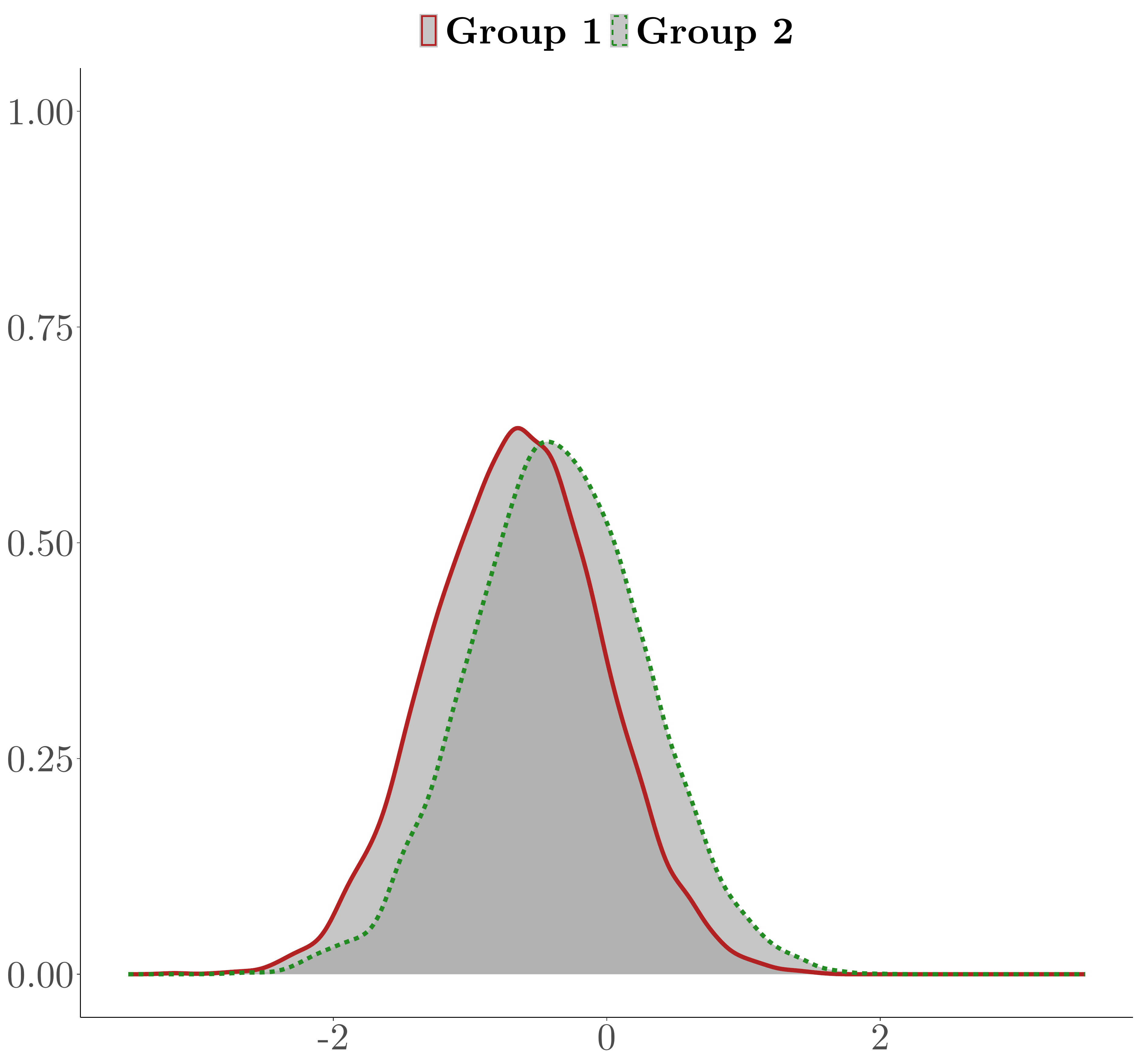} &
    \includegraphics[width=0.25\textwidth]{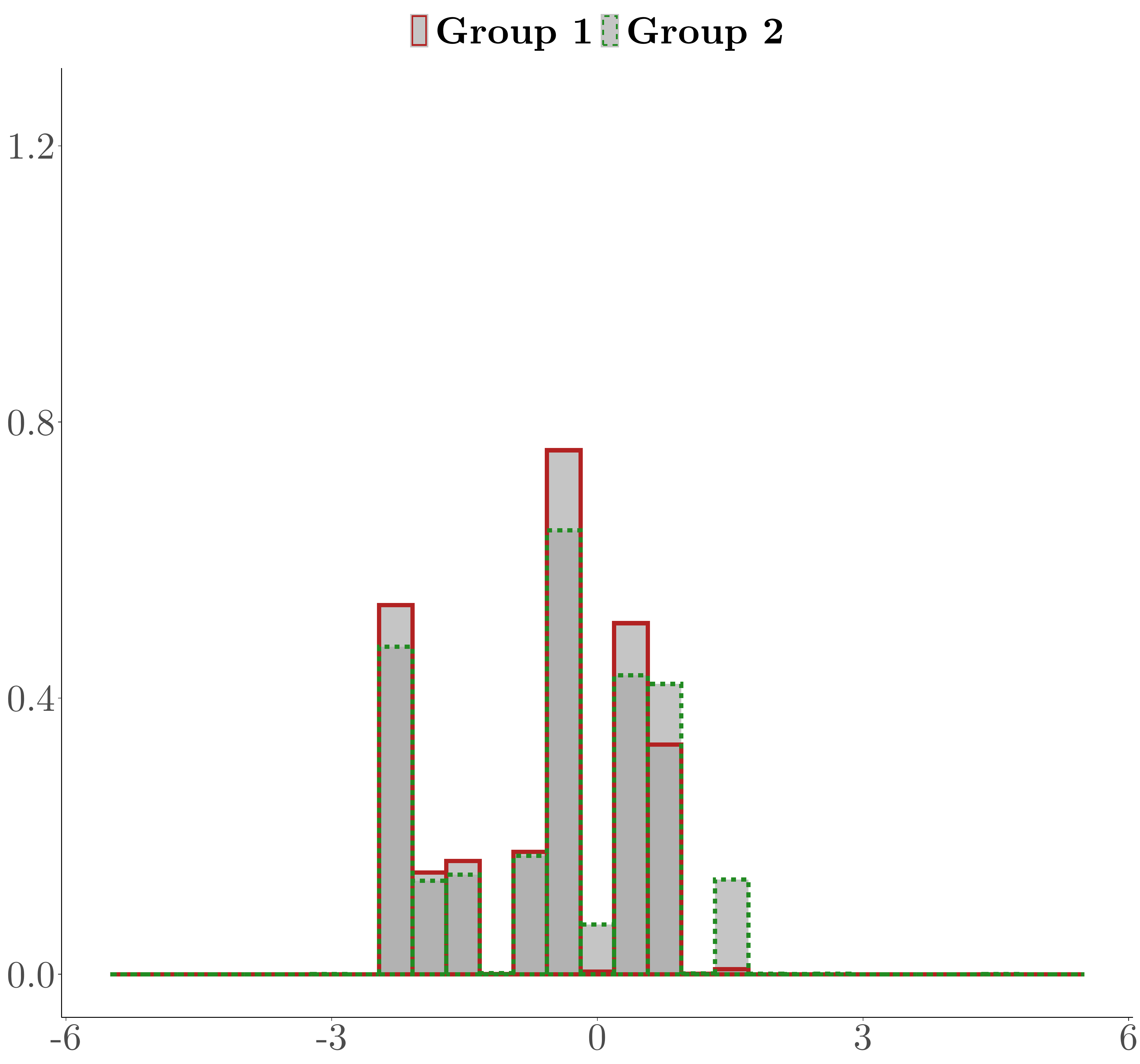} \\
    \bottomrule
\end{tabular}
\end{table}

\begin{figure*}[t]
    \includegraphics[width=0.320\linewidth]{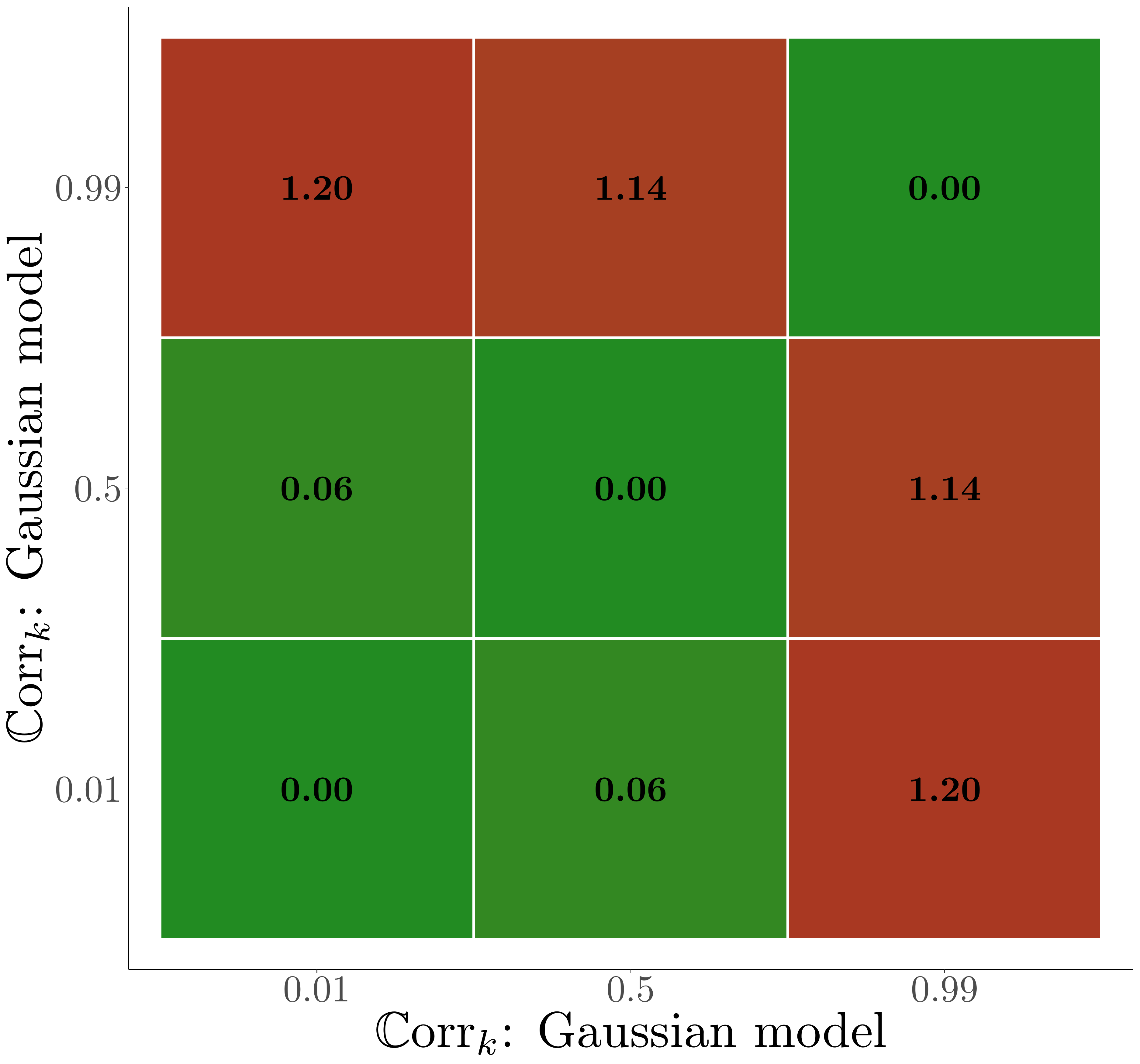}
    \includegraphics[width=0.320\linewidth]{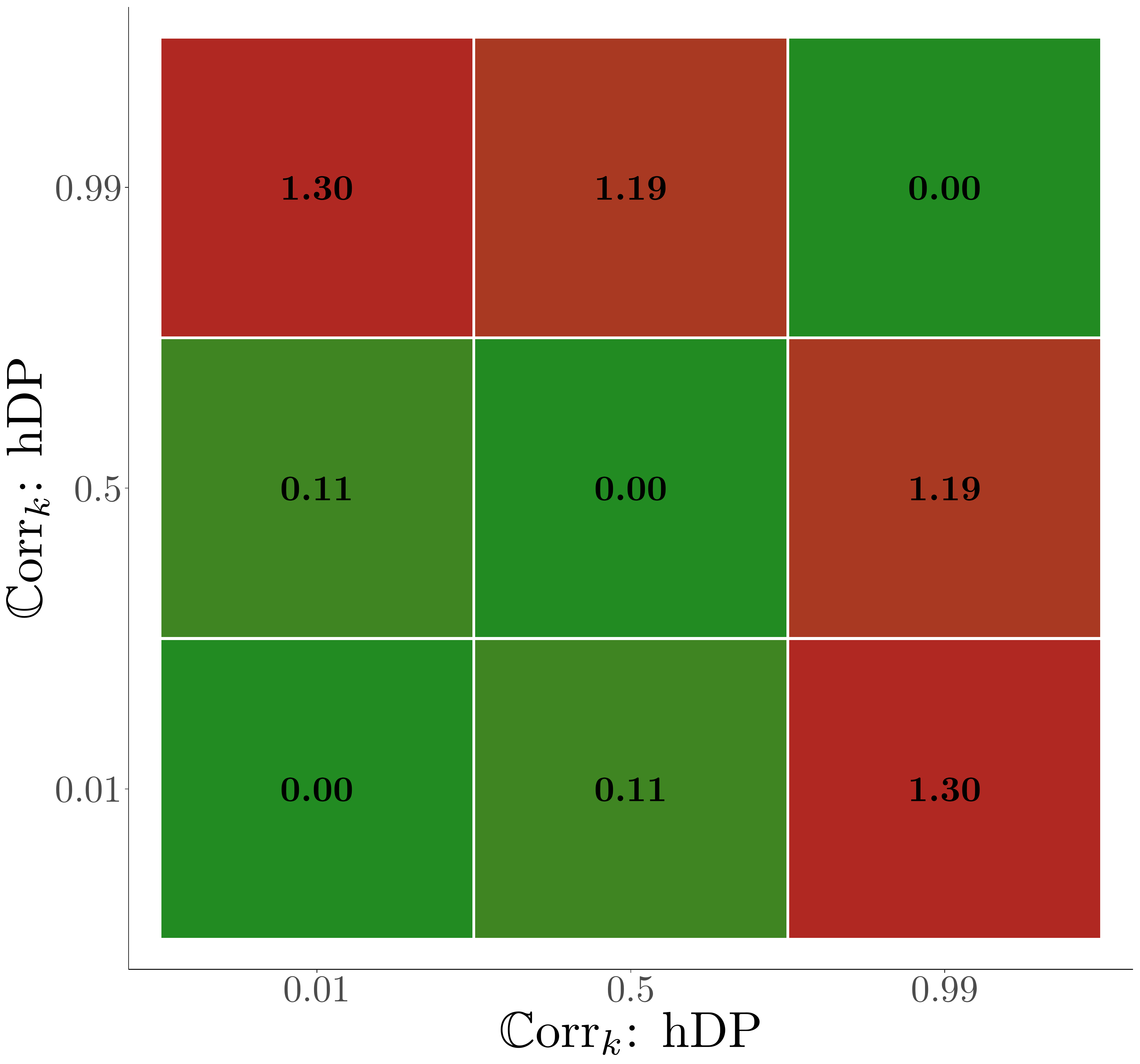}
    \includegraphics[width=0.320\linewidth]{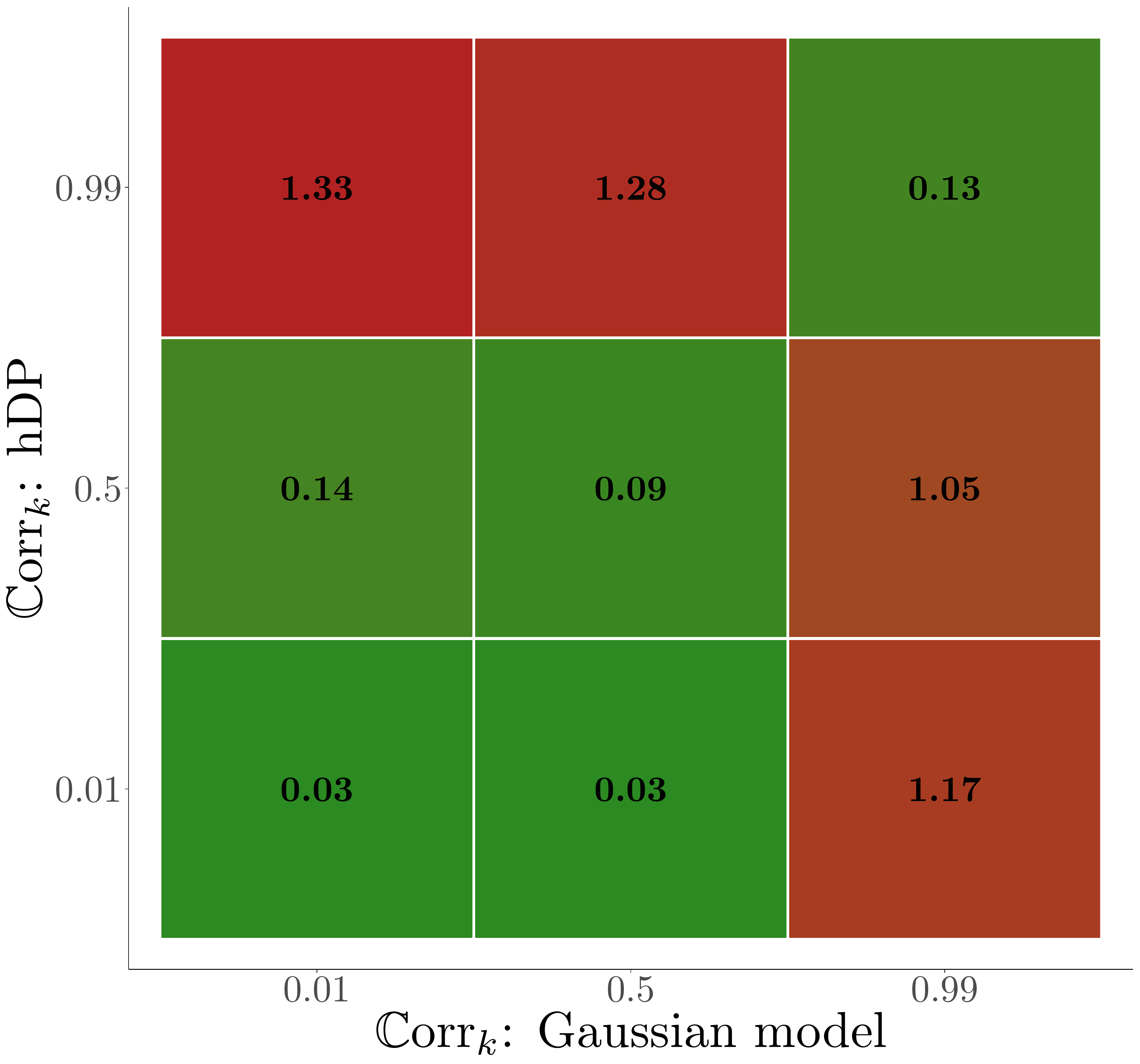} \\
    \includegraphics[width=0.320\linewidth]{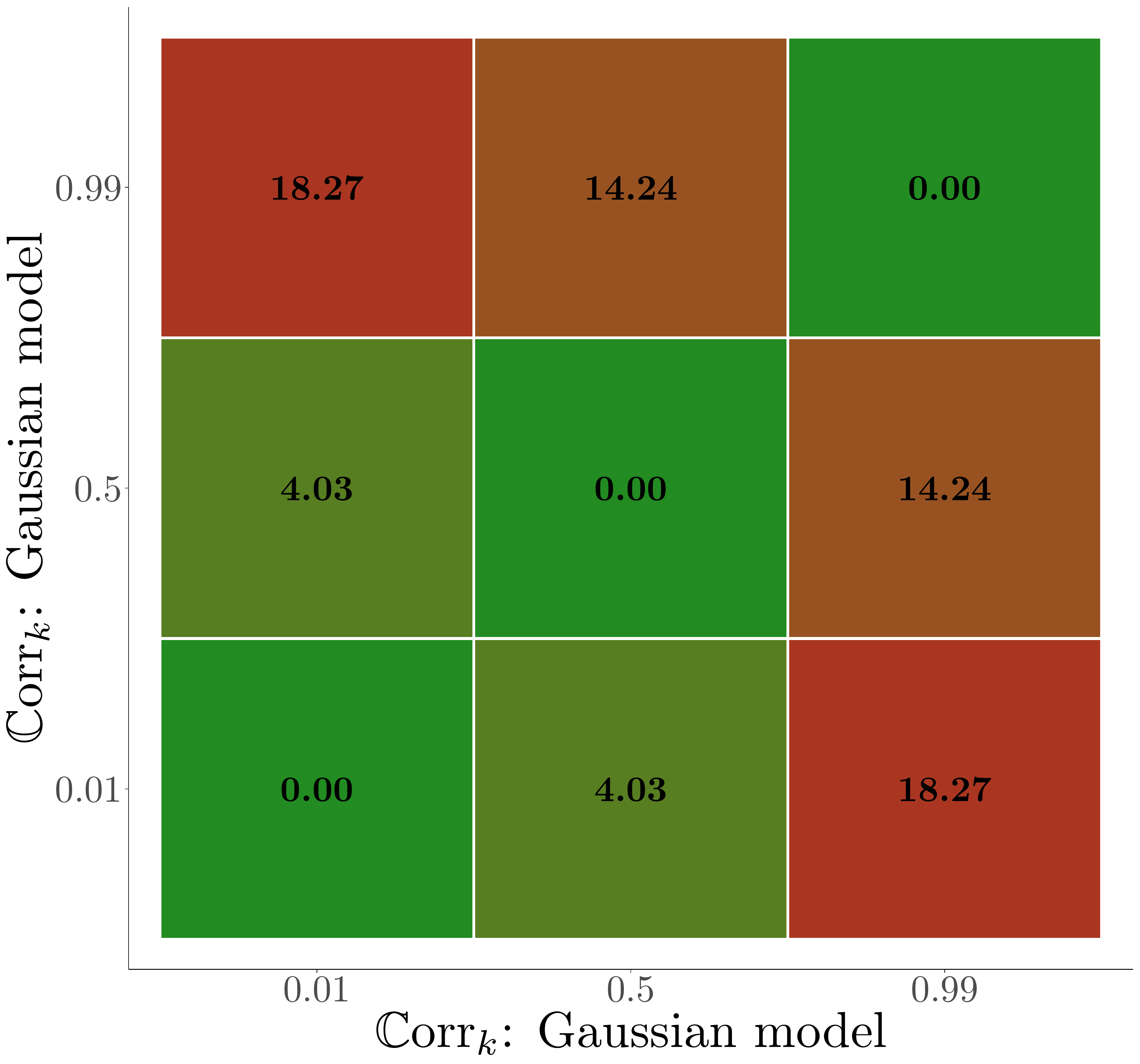}
    \includegraphics[width=0.320\linewidth]{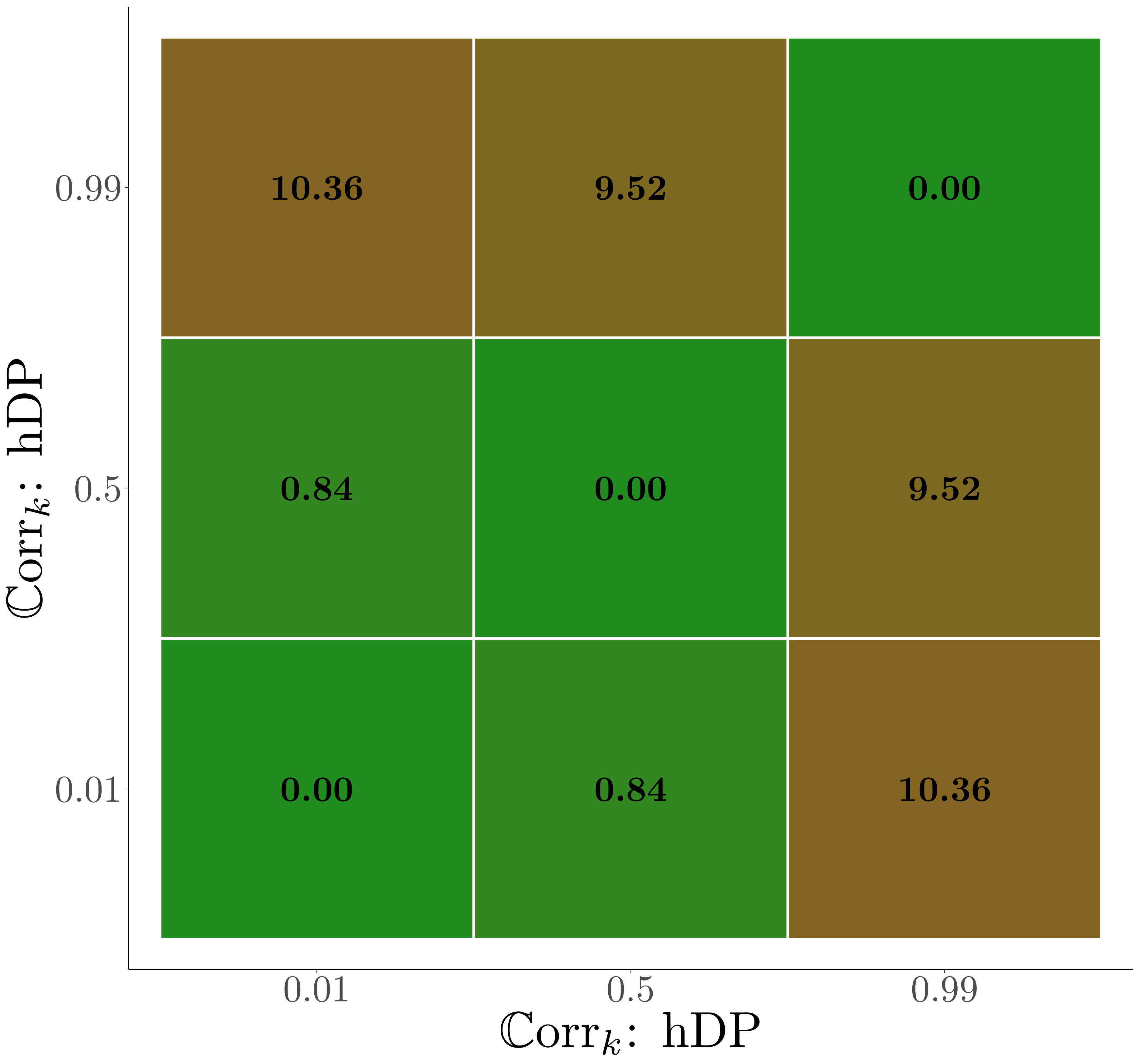}
    \includegraphics[width=0.320\linewidth]{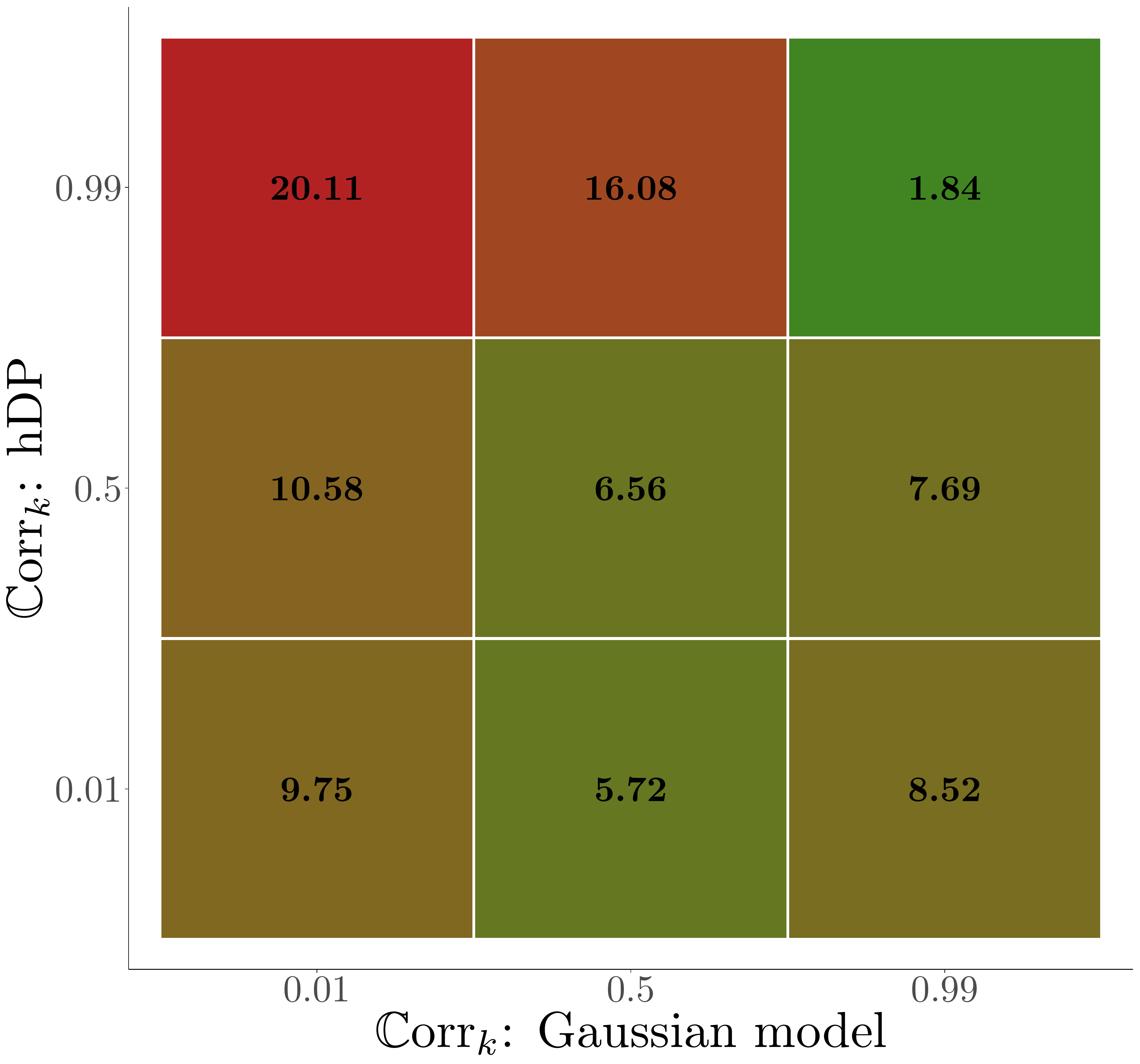}
    \caption{Absolute difference between the empirical averages of two samples of size $M = 10,000$ from the predictive distribution of Group 2 for different values of kernel correlation for Gaussian vs. Gaussian (\textbf{left}), hDP vs. hDP (\textbf{centre}), and Gaussian vs. hDP (\textbf{right}) for $v = 1/4$, $t^{2} = 2$, and $\sigma = t/\sqrt{1/(1-v)^{2} - 1}$. \textbf{Top row}: the data are simulated from \eqref{eq:post_smpl}; \textbf{Bottom row}: the data comes from the Palmer Penguins dataset \cite{horst2020palmerpenguins}.}
    \label{fig:corr_comparison}
\end{figure*}

We provide two illustrations of how the dependence a priori impacts posterior inference, on simulated and real data, highlighting the necessity of matching dependencies for accurate model comparisons. First, we generate the observations in each group as independent sequences $(X_{1,j})_{j=1}^{n_{1}}$ and $(X_{2,j})_{j=1}^{n_{2}}$ according to the following sampling scheme,
\begin{equation} \label{eq:post_smpl}
\begin{split}
	X_{1,1},\dots,X_{1,n_{1}} | \tilde{P}_{1} \stackrel{\iid}{\sim} \tilde{P}_{1} \qquad \tilde{P}_{1} &\sim \mathrm{hDP}\bigl(c = 10,c_{0} = 10, P_{0} = \mathcal{N}(-1,2)\bigr), \\
   X_{2,1},\dots,X_{2,n_{2}} | \tilde{P}_{2} \stackrel{\iid}{\sim} \tilde{P}_{2} \qquad  \tilde{P}_{2} &\sim \mathrm{hDP}\bigl(c = 10,c_{0} = 10, P_{0} = \mathcal{N}(1,2)\bigr).
\end{split}
\end{equation}
We consider unbalanced groups, as borrowing information is particularly useful in this setting, with $n_1 = 200$ and $n_2 = 5$. We consider three different values of the index, namely $\corr_{k} \in \{0.01, 0.50, 0.99\}$, corresponding to the situation of almost independence, intermediate dependence, and almost exchangeability, respectively. In \cref{tab:models}, we show the empirical distribution of a sample of size $M = 10,000$ from the posterior mean measure, which coincides with the posterior one-step-ahead predictive distribution, for the two groups for both models. We set $v = 1/4$, $t^{2} = 2$, $\sigma = \sigma^{*}/\sqrt{2}$, and the values of $\rho$, $c_{0}$, $c$ determined by the value of the kernel correlation, as shown in \cref{fig:corr_comparison}. We notice that the two distributions are more similar for values of the kernel correlation close to $1$, as expected, and the second group is heavily affected by the larger number of observations in the first group. This intuition can be quantified through the absolute difference of their means, which we estimate through the sample means. The results are shown in the top row of \cref{fig:corr_comparison}. On the left, we plot the absolute mean differences for different values of the kernel correlations for the Gaussian model; at the centre, we reproduce the same analysis for the hDP. As expected, the estimates tend to be closer for similar values of the index; in other words, the values near the main diagonal tend to be smaller. On the right, we compare the absolute mean differences between the Gaussian model and the hDP. We notice that this distance tends to be bigger between two instances of the same model (both Gaussian or hDP) with different kernel correlations, rather than between a Gaussian and an hDP model with the same kernel correlation. 

We conclude with a similar type of analysis on the Palmer Archipelago (Antarctica) penguin data \citep{horst2020palmerpenguins} with the goal of making predictions on the flipper length for male and female penguins. We consider all male penguins in the dataset (168) and a subsample of 5 female penguins to benefit from the borrowing of information across groups. The presence of ties could be ascribed to rounding error, leading to a dominated hierarchical model such as the parametric Gaussian hierarchical model in Example~\ref{par_ex}, or to the presence of latent subpopulations, leading to an a.s. discrete hierarchical model such as the hierarchical Dirichlet process in \eqref{hDP}. We show that setting the same kernel correlation a priori leads to a more meaningful model comparison: indeed, the predictions using the same model with different kernel correlations can lead to differences between mean predictions that are larger than those obtained by using different models with the same value of kernel correlation. This is pictured in the bottom row of \cref{fig:corr_comparison}, where, e.g., the difference between the hierarchical Gaussian model with kernel correlation 0.5 and the hDP model with kernel correlation 0.99 dominates the difference between the two models with the same kernel correlation 0.5. This analysis provides evidence of the importance of fixing the same prior kernel correlation, when possible, to mitigate the effect of the choice of the hyperparameters in the model comparison.

\section{Discussion}

In this work, we have introduced a measure of partial exchangeability by quantifying the dependence of random probability measures through reproducing kernel Hilbert spaces. A distinctive feature of our index, termed kernel correlation, is to detect exchangeability for a broad class of models, including almost surely discrete random probabilities, their posterior updates, mixture models, and standard parametric models. We have identified some mild conditions on the marginal distributions that are remarkably agnostic of the dependence structure and that we expect to hold in many other settings not considered in this work.

The kernel correlation extends the widely used set-wise correlation and coincides with the former for multivariate species sampling models \cite{franzolini_multivariate_2025}, which encompass most discrete priors in the Bayesian nonparametric literature and is easily computable in such settings. We have shown that these computations easily extend to mixture models as well. For other random probabilities, such as those arising from posterior and parametric models, we have provided a simple and efficient estimator that only uses four observables of the partially exchangeable sequence. This makes it possible to perform a fair model comparison between parametric and nonparametric models by fixing the same amount of prior dependence. We were also able to investigate the behaviour of the dependence structure a posteriori. Remarkably, we have found that the kernel correlation for the hierarchical Dirichlet process \cite{teh_hierarchical_2006} goes to zero at a parametric rate of convergence. As the dependence structure drives the borrowing of information, our results show that as the sample size grows, each additional datapoint contributes progressively less to the borrowing of information across groups. This work lays the groundwork for analyzing the dependence a posteriori for other hierarchical models \citep{camerlenghi_distribution_2019} or more general dependent priors \cite{quintana_dependent_2022, wade2025bayesian}, complementing well-established frequentist asymptotic analyses on the recovery of the true distribution in partially exchangeable settings \citep{nguyen_borrowing_2016, catalano_posterior_2022}.

\bibliography{rkhs}

\begin{thebibliography}{7}
\providecommand{\natexlab}[1]{#1}
\providecommand{\url}[1]{\texttt{#1}}
\expandafter\ifx\csname urlstyle\endcsname\relax
  \providecommand{\doi}[1]{doi: #1}\else
  \providecommand{\doi}{doi: \begingroup \urlstyle{rm}\Url}\fi

\bibitem[Camerlenghi et~al.(2018)Camerlenghi, Lijoi, and
  Prünster]{camerlenghi_bayesian_2018}
F.~Camerlenghi, A.~Lijoi, and I.~Prünster.
\newblock {Bayesian Nonparametric Inference Beyond the Gibbs-Type Framework}.
\newblock \emph{Scand. J. Statist.}, 45\penalty0 (4):\penalty0 1062--1091,
  2018.

\bibitem[Camerlenghi et~al.(2019)Camerlenghi, Lijoi, Orbanz, and
  Prünster]{camerlenghi_distribution_2019}
F.~Camerlenghi, A.~Lijoi, P.~Orbanz, and I.~Prünster.
\newblock {Distribution Theory for Hierarchical Processes}.
\newblock \emph{Ann. Statist.}, 47\penalty0 (1):\penalty0 67--92, 2019.

\bibitem[Catalano et~al.(2023)Catalano, Del~Sole, Lijoi, and
  Prünster]{catalano_unified_2023}
M.~Catalano, C.~Del~Sole, A.~Lijoi, and I.~Prünster.
\newblock {A Unified Approach to Hierarchical Random Measures}.
\newblock \emph{Sankhya A}, 86:\penalty0 255--287, 2023.

\bibitem[Hoffmann-Jorgensen and Pisier(1976)]{hoffmann_law_1976}
J.~Hoffmann-Jorgensen and G.~Pisier.
\newblock {The Law of Large Numbers and the Central Limit Theorem in Banach
  Spaces}.
\newblock \emph{Ann. Probab.}, 4\penalty0 (4):\penalty0 587--599, 1976.

\bibitem[Sriperumbudur et~al.(2011)Sriperumbudur, Fukumizu, and
  Lanckriet]{sriperumbudur_universality_2011}
B.~K. Sriperumbudur, K.~Fukumizu, and G.~R.~G. Lanckriet.
\newblock {Universality, Characteristic Kernels and RKHS Embedding of
  Measures}.
\newblock \emph{J. Mach. Learn. Res.}, 12\penalty0 (70):\penalty0 2389--2410,
  2011.

\bibitem[Teh et~al.(2006)Teh, Jordan, Beal, and Blei]{teh_hierarchical_2006}
Y.~W. Teh, M.~I. Jordan, M.~J. Beal, and D.~M. Blei.
\newblock {Hierarchical Dirichlet Processes}.
\newblock \emph{J. Amer. Statist. Assoc.}, 101\penalty0 (476):\penalty0
  1566--1581, 2006.

\bibitem[van~der Vaart and Wellner(1996)]{vandervaart_weak_1996}
A.~W. van~der Vaart and J.~Wellner.
\newblock \emph{{Weak Convergence and Empirical Processes with Applications to
  Statistics}}.
\newblock Springer, 1996.

\end{thebibliography}


\begin{thebibliography}{65}
\providecommand{\natexlab}[1]{#1}
\providecommand{\url}[1]{\texttt{#1}}
\expandafter\ifx\csname urlstyle\endcsname\relax
  \providecommand{\doi}[1]{doi: #1}\else
  \providecommand{\doi}{doi: \begingroup \urlstyle{rm}\Url}\fi

\bibitem[Aronszajn(1950)]{aronszajn_theory_1950}
N.~Aronszajn.
\newblock {Theory of Reproducing Kernels}.
\newblock \emph{Trans. Amer. Math. Soc.}, 68\penalty0 (3):\penalty0 337--404,
  1950.

\bibitem[Ascolani et~al.(2023)Ascolani, Franzolini, Lijoi, and
  Prünster]{ascolani_nonparametric_2023}
F.~Ascolani, B.~Franzolini, A.~Lijoi, and I.~Prünster.
\newblock {Nonparametric Priors with Full-Range Borrowing of Information}.
\newblock \emph{Biometrika}, 111\penalty0 (3):\penalty0 945--969, 2023.

\bibitem[Bach and Jordan(2002)]{bach_kernel_2002}
F.~R. Bach and M.~I. Jordan.
\newblock {Kernel Independent Component Analysis}.
\newblock \emph{J. Mach. Learn. Res.}, 3:\penalty0 1--48, 2002.

\bibitem[Beraha et~al.(2021)Beraha, Guglielmi, and
  Quintana]{beraha_semihierarchical_2021}
M.~Beraha, A.~Guglielmi, and F.~A. Quintana.
\newblock {The Semi-Hierarchical Dirichlet Process and Its Application to
  Clustering Homogeneous Distributions}.
\newblock \emph{Bayesian Anal.}, 16\penalty0 (4):\penalty0 1187--1219, 2021.

\bibitem[Berlinet and Thomas-Agnan(2004)]{berlinet_reproducing_2004}
A.~Berlinet and C.~Thomas-Agnan.
\newblock \emph{{Reproducing Kernel Hilbert Spaces in Probability and
  Statistics}}.
\newblock Springer, 2004.

\bibitem[Blackwell and MacQueen(1973)]{blackwell_ferguson_1973}
D.~Blackwell and J.~B. MacQueen.
\newblock Ferguson distributions via {P}ólya urn schemes.
\newblock \emph{Ann. Statist.}, 1\penalty0 (2):\penalty0 353--355, 1973.

\bibitem[Bochner(1959)]{bochner_lectures_1959}
S.~Bochner.
\newblock \emph{Lectures on Fourier Integrals; with an Author's Supplement on
  Monotonic Functions, Stieltjes Integrals, and Harmonic Analysis}.
\newblock Princeton University Press, 1959.

\bibitem[Camerlenghi et~al.(2019)Camerlenghi, Lijoi, Orbanz, and
  Prünster]{camerlenghi_distribution_2019}
F.~Camerlenghi, A.~Lijoi, P.~Orbanz, and I.~Prünster.
\newblock {Distribution Theory for Hierarchical Processes}.
\newblock \emph{Ann. Statist.}, 47\penalty0 (1):\penalty0 67--92, 2019.

\bibitem[Catalano et~al.(2021)Catalano, Lijoi, and
  Prünster]{catalano_measuring_2021}
M.~Catalano, A.~Lijoi, and I.~Prünster.
\newblock {Measuring dependence in the Wasserstein distance for Bayesian
  nonparametric models}.
\newblock \emph{Ann. Statist.}, 49\penalty0 (5):\penalty0 2916--2947, 2021.

\bibitem[Catalano et~al.(2022)Catalano, De~Blasi, Lijoi, and
  Prünster]{catalano_posterior_2022}
M.~Catalano, P.~De~Blasi, A.~Lijoi, and I.~Prünster.
\newblock {Posterior Asymptotics for Boosted Hierarchical Dirichlet Process
  Mixtures}.
\newblock \emph{J. Mach. Learn. Res.}, 23\penalty0 (80):\penalty0 1--23, 2022.

\bibitem[Catalano et~al.(2023)Catalano, Del~Sole, Lijoi, and
  Prünster]{catalano_unified_2023}
M.~Catalano, C.~Del~Sole, A.~Lijoi, and I.~Prünster.
\newblock {A Unified Approach to Hierarchical Random Measures}.
\newblock \emph{Sankhya A}, 86:\penalty0 255--287, 2023.

\bibitem[Catalano et~al.(2024)Catalano, Lavenant, Lijoi, and
  Prünster]{catalano_wasserstein_2024}
M.~Catalano, H.~Lavenant, A.~Lijoi, and I.~Prünster.
\newblock {A Wasserstein Index of Dependence for Random Measures}.
\newblock \emph{J. Amer. Statist. Assoc.}, 119\penalty0 (547):\penalty0
  2396--2406, 2024.

\bibitem[Chakrborty et~al.(2012)Chakrborty, Ghosh, and
  Mallick]{chakraborty_bayesian_2012}
S.~Chakrborty, M.~Ghosh, and B.~K. Mallick.
\newblock {Bayesian Nonlinear Regression for Large $p$ small $n$ problems}.
\newblock \emph{J. Multivariate Anal.}, 108:\penalty0 28--40, 2012.

\bibitem[Chen et~al.(2019)Chen, Barp, Briol, Gorham, Girolami, Mackey, and
  Oates]{chen2019stein}
W.~Y. Chen, A.~Barp, F.-X. Briol, J.~Gorham, M.~Girolami, L.~Mackey, and
  C.~Oates.
\newblock Stein point {M}arkov chain {M}onte {C}arlo.
\newblock In \emph{Proceedings of the 36th International Conference on Machine
  Learning}, volume~97, pages 1011--1021, 2019.

\bibitem[Cifarelli and Regazzini(1978)]{cifarelli_problemi_1978}
M.~D. Cifarelli and E.~Regazzini.
\newblock {Problemi Statistici non Parametrici in Condizioni di Scambiabilità
  Parziale: Impiego di Medie Associative}.
\newblock Technical report, Istituto di Matematica Finanziaria
  dell’Università di Torino, 1978.

\bibitem[Colombi et~al.(2025)Colombi, Argiento, Camerlenghi, and
  Paci]{colombi_hierarchical_2025}
A.~Colombi, R.~Argiento, F.~Camerlenghi, and L.~Paci.
\newblock {Hierarchical Mixture of Finite Mixtures}.
\newblock \emph{Bayesian Anal.}, Advance Publication:\penalty0 1--29, 2025.

\bibitem[Cortes et~al.(2012)Cortes, Mohri, and
  Rostamizadeh]{cortes_algorithms_2012}
C.~Cortes, M.~Mohri, and A.~Rostamizadeh.
\newblock Algorithms for learning kernels based on centered alignment.
\newblock \emph{J. Mach. Learn. Res.}, 13\penalty0 (28):\penalty0 795--828,
  2012.

\bibitem[Cristianini et~al.(2001)Cristianini, Shawe-Taylor, Elisseeff, and
  Kandola]{cristianini_kernel_2001}
N.~Cristianini, J.~Shawe-Taylor, A.~Elisseeff, and J.~Kandola.
\newblock On kernel-target alignment.
\newblock In \emph{Advances in Neural Information Processing Systems},
  volume~14, pages 367--373. MIT Press, 2001.

\bibitem[de~Finetti(1937)]{definetti_prevision_1937}
B.~de~Finetti.
\newblock {La Prévision : ses Lois Logiques, ses Sources Subjectives}.
\newblock \emph{Ann. Inst. Henri Poincaré}, 7\penalty0 (1):\penalty0 1--68,
  1937.

\bibitem[de~Finetti(1938)]{definetti_condition_1938}
B.~de~Finetti.
\newblock {Sur la Condition d'Equivalence Partielle}.
\newblock \emph{Act. Sci. Ind.}, 739:\penalty0 5--18, 1938.

\bibitem[Denti et~al.(2023)Denti, Camerlenghi, Guindani, and
  Mira]{denti_common_2023}
F.~Denti, F.~Camerlenghi, M.~Guindani, and A.~Mira.
\newblock {A Common Atoms Model for the Bayesian Nonparametric Analysis of
  Nested Data}.
\newblock \emph{J. Amer. Statist. Assoc.}, 118\penalty0 (541):\penalty0
  405--416, 2023.

\bibitem[Efron and Morris(1973)]{efron_steins_1973}
B.~Efron and C.~Morris.
\newblock {Stein's Estimation Rule and Its Competitors -- An Empirical Bayes
  Approach}.
\newblock \emph{J. Amer. Statist. Assoc.}, 68\penalty0 (341):\penalty0
  117--130, 1973.

\bibitem[Escobar and West(1995)]{escobar_bayesian_1995}
M.~D. Escobar and M.~West.
\newblock {Bayesian Density Estimation and Inference Using Mixtures}.
\newblock \emph{J. Amer. Statist. Assoc.}, 90\penalty0 (430):\penalty0
  577--588, 1995.

\bibitem[Ferguson(1983)]{ferguson_bayesian_1983}
T.~S. Ferguson.
\newblock {Bayesian Density Estimation by Mixtures of Normal Distributions}.
\newblock In M.~H. Rizvi, J.~S. Rustagi, and D.~Siegmund, editors, \emph{Recent
  Advances in Statistics}, pages 287--302. Academic Press, 1983.

\bibitem[Franzolini et~al.(2025)Franzolini, Lijoi, Prünster, and
  Rebaudo]{franzolini_multivariate_2025}
B.~Franzolini, A.~Lijoi, I.~Prünster, and G.~Rebaudo.
\newblock {Multivariate Species Sampling Models}.
\newblock \emph{arXiv:2503.24004}, 2025.

\bibitem[Fukumizu et~al.(2013)Fukumizu, Song, and Gretton]{fukumizu2013}
K.~Fukumizu, L.~Song, and A.~Gretton.
\newblock Kernel {B}ayes' rule: {B}ayesian inference with positive definite
  kernels.
\newblock \emph{J. Mach. Learn. Res.}, 14\penalty0 (118):\penalty0 3753--3783,
  2013.

\bibitem[Gelman and Hill(2007)]{gelman_data_2007}
A.~Gelman and J.~Hill.
\newblock \emph{{Data Analysis Using Regression and Multilevel/Hierarchical
  Models}}.
\newblock Cambridge University Press, 2007.

\bibitem[Gelman et~al.(2003)Gelman, Carlin, Stern, and
  Rubin]{gelman_bayesian_2003}
A.~Gelman, J.~B. Carlin, H.~S. Stern, and D.~B. Rubin.
\newblock \emph{{Bayesian Data Analysis. Second Edition}}.
\newblock Chapman \& Hall, 2003.

\bibitem[Gretton et~al.(2005{\natexlab{a}})Gretton, Bousquet, Smola, and
  Schölkopf]{gretton_measuring_2005}
A.~Gretton, O.~Bousquet, A.~Smola, and B.~Schölkopf.
\newblock {Measuring Statistical Dependence with Hilbert-Schmidt norms}.
\newblock In \emph{Proceedings of the 16th International Conference on
  Algorithmic Learning Theory}, pages 63--77. Springer, 2005{\natexlab{a}}.

\bibitem[Gretton et~al.(2005{\natexlab{b}})Gretton, Herbrich, Bousquet, Smola,
  and Schölkopf]{gretton_kernel_2005}
A.~Gretton, R.~Herbrich, O.~Bousquet, A.~Smola, and B.~Schölkopf.
\newblock {Kernel Methods for Measuring Independence}.
\newblock \emph{J. Mach. Learn. Res.}, 6\penalty0 (70):\penalty0 2075--2129, b
  2005{\natexlab{b}}.

\bibitem[Griffin and Leisen(2017)]{griffin_compound_2017}
J.~E. Griffin and F.~Leisen.
\newblock {Compound Random Measures and Their Use in Bayesian Non-Parametrics}.
\newblock \emph{J. R. Stat. Soc. Ser. B}, 79\penalty0 (2):\penalty0 525--545,
  2017.

\bibitem[Guella(2022)]{guella2022gaussian}
J.~C. Guella.
\newblock {On Gaussian kernels on Hilbert spaces and kernels on hyperbolic
  spaces}.
\newblock \emph{J. Approx. Theory}, 279:\penalty0 105765, 2022.

\bibitem[Horiguchi et~al.(2024)Horiguchi, Chan, and Ma]{horiguchi_tree_2024}
A.~Horiguchi, C.~Chan, and L.~Ma.
\newblock {A Tree Perspective on Stick-Breaking Models in Covariate-Dependent
  Mixtures}.
\newblock \emph{Bayesian Anal.}, Advance Publication:\penalty0 1--28, 2024.

\bibitem[Horst et~al.(2020)Horst, Hill, and Gorman]{horst2020palmerpenguins}
A.~M. Horst, A.~P. Hill, and K.~B. Gorman.
\newblock \emph{palmerpenguins: Palmer Archipelago (Antarctica) penguin data},
  2020.
\newblock R package version 0.1.0.

\bibitem[James et~al.(2009)James, Lijoi, and Prünster]{james_posterior_2009}
L.~F. James, A.~Lijoi, and I.~Prünster.
\newblock {Posterior Analysis for Normalized Random Measures with Independent
  Increments}.
\newblock \emph{Scand. J. Statist.}, 36\penalty0 (1):\penalty0 76--97, 2009.

\bibitem[Kingman(1967)]{kingman_completely_1967}
J.~F.~C. Kingman.
\newblock {Completely Random Measures}.
\newblock \emph{Pac. J. Math.}, 21\penalty0 (1):\penalty0 59--78, 1967.

\bibitem[Kornblith et~al.(2019)Kornblith, Norouzi, Lee, and
  Hinton]{kornblith2019similarity}
S.~Kornblith, M.~Norouzi, H.~Lee, and G.~Hinton.
\newblock Similarity of neural network representations revisited.
\newblock In \emph{International conference on machine learning}, pages
  3519--3529. PMlR, 2019.

\bibitem[Legramanti et~al.(2025)Legramanti, Durante, and
  Alquier]{legramanti2025concentration}
S.~Legramanti, D.~Durante, and P.~Alquier.
\newblock {Concentration of discrepancy-based approximate Bayesian computation
  via Rademacher complexity}.
\newblock \emph{Ann. Statist.}, 53\penalty0 (1):\penalty0 37 -- 60, 2025.

\bibitem[Leisen et~al.(2013)Leisen, Lijoi, and Spanò]{leisen_vector_2013}
F.~Leisen, A.~Lijoi, and D.~Spanò.
\newblock {A Vector of Dirichlet Processes}.
\newblock \emph{Electron. J. Stat.}, 7:\penalty0 62--90, 2013.

\bibitem[Lijoi and Prünster(2010)]{lijoi_models_2010}
A.~Lijoi and I.~Prünster.
\newblock {Models Beyond the Dirichlet Process}.
\newblock In N.~L. Hjort, C.~Holmes, P.~Müller, and S.~G. Walker, editors,
  \emph{Bayesian Nonparametrics}, pages 80--136. Cambridge University Press,
  2010.

\bibitem[Lijoi et~al.(2023)Lijoi, Prünster, and Rebaudo]{lijoi_flexible_2023}
A.~Lijoi, I.~Prünster, and G.~Rebaudo.
\newblock {Flexible Clustering via Hidden Hierarchical Dirichlet Priors}.
\newblock \emph{Scand. J. Statist.}, 50\penalty0 (1):\penalty0 213--234, 2023.

\bibitem[Lindley and Smith(1972)]{lindley_bayes_1972}
D.~V. Lindley and A.~F.~M. Smith.
\newblock {Bayes Estimates for the Linear Model}.
\newblock \emph{J. R. Stat. Soc. Ser. B}, 34\penalty0 (1):\penalty0 1--41,
  1972.

\bibitem[Liu and Wang(2016)]{qiang2016stein}
Q.~Liu and D.~Wang.
\newblock Stein variational gradient descent: a general purpose {B}ayesian
  inference algorithm.
\newblock In \emph{Proceedings of the 30th International Conference on Neural
  Information Processing Systems}, page 2378–2386, Red Hook, NY, USA, 2016.
  Curran Associates Inc.
\newblock ISBN 9781510838819.

\bibitem[Lo(1984)]{lo_class_1984}
A.~Y. Lo.
\newblock {On a Class of Bayesian Nonparametric Estimates: I. Density
  Estimates}.
\newblock \emph{Ann. Statist.}, 12\penalty0 (1):\penalty0 351--357, 1984.

\bibitem[MacEachern(1999)]{maceachern_dependent_1999}
S.~N. MacEachern.
\newblock {Dependent Nonparametric Processes}.
\newblock In \emph{ASA Proceedings of the Section on Bayesian Statistical
  Science}, 1999.

\bibitem[MacEachern(2000)]{maceachern_dependent_2000}
S.~N. MacEachern.
\newblock {Dependent Dirichlet Processes}.
\newblock Technical report, Ohio State University, 2000.

\bibitem[MacLehose and Dunson(2009)]{maclehose_nonparametric_2009}
R.~F. MacLehose and D.~B. Dunson.
\newblock {Nonparametric Bayes Kernel-Based Priors for Functional Data
  Analysis}.
\newblock \emph{Statistica Sinica}, 19\penalty0 (2):\penalty0 611--629, 2009.

\bibitem[Muandet et~al.(2017)Muandet, Fukumizu, Sriperumbudur, and
  Schölkopf]{muandet_kernel_2017}
K.~Muandet, K.~Fukumizu, B.~Sriperumbudur, and B.~Schölkopf.
\newblock {Kernel Mean Embedding of Distributions: A Review and Beyond}.
\newblock \emph{Found. Trends Mach. Learn.}, 10\penalty0 (1--2):\penalty0
  1--141, 2017.

\bibitem[Nguyen(2013)]{nguyen_convergence_2013}
X.~L. Nguyen.
\newblock {Convergence of Latent Mixing Measures in Finite and Infinite Mixture
  Models}.
\newblock \emph{Ann. Statist.}, 41\penalty0 (1):\penalty0 370--400, 2013.

\bibitem[Nguyen(2016)]{nguyen_borrowing_2016}
X.~L. Nguyen.
\newblock {Borrowing Strengh in Hierarchical Bayes: Posterior Concentration of
  the Dirichlet Base Measure}.
\newblock \emph{Bernoulli}, 22\penalty0 (3):\penalty0 1535--1571, 2016.

\bibitem[Park et~al.(2016)Park, Jitkrittum, and
  Sejdinovic]{park2016approximate}
M.~Park, W.~Jitkrittum, and D.~Sejdinovic.
\newblock K2-abc: Approximate {B}ayesian computation with kernel embeddings.
\newblock In \emph{Proceedings of the 19th International Conference on
  Artificial Intelligence and Statistics}, volume~51, pages 398--407, Cadiz,
  Spain, 09--11 May 2016. PMLR.

\bibitem[Pillai et~al.(2007)Pillai, Wu, Liang, Mukherjee, and
  Wolpert]{pillai_characterizing_2007}
N.~S. Pillai, Q.~Wu, F.~Liang, S.~Mukherjee, and R.~L. Wolpert.
\newblock {Characterizing the Function Space for Bayesian Kernel Models}.
\newblock \emph{J. Mach. Learn. Res.}, 8\penalty0 (62):\penalty0 1769--1797,
  2007.

\bibitem[Pitman(1996)]{pitman_developments_1996}
J.~Pitman.
\newblock {Some Developments of the Blackwell-MacQueen Urn Scheme}.
\newblock In T.~S. Ferguson, L.~S. Shapley, and J.~B. MacQueen, editors,
  \emph{Statistics, Probability and Game Theory: Papers in Honor of David
  Blackwell}, pages 245--268. Institute of Mathematical Statistics, 1996.

\bibitem[Quintana et~al.(2022)Quintana, Müller, Jara, and
  MacEachern]{quintana_dependent_2022}
F.~A. Quintana, P.~Müller, A.~Jara, and S.~N. MacEachern.
\newblock {The Dependent Dirichlet Process and Related Models}.
\newblock \emph{Stat. Sci.}, 37\penalty0 (1):\penalty0 24--41, 2022.

\bibitem[Rasmussen and Williams(2006)]{RasmussenWilliams2006}
C.~E. Rasmussen and C.~K.~I. Williams.
\newblock \emph{Gaussian Processes for Machine Learning}.
\newblock MIT Press, 2006.

\bibitem[Regazzini et~al.(2003)Regazzini, Lijoi, and
  Prünster]{regazzini_distributional_2003}
E.~Regazzini, A.~Lijoi, and I.~Prünster.
\newblock {Distributional Results for Means of Normalized Random Measures with
  Independent Increments}.
\newblock \emph{Ann. Statist.}, 31\penalty0 (2):\penalty0 560--585, 2003.

\bibitem[Rodríguez et~al.(2008)Rodríguez, Dunson, and
  Gelfand]{rodriguez_nested_2008}
A.~Rodríguez, D.~B. Dunson, and A.~E. Gelfand.
\newblock {The Nested Dirichlet Process}.
\newblock \emph{J. Amer. Statist. Assoc.}, 103\penalty0 (483):\penalty0
  1131--1154, 2008.

\bibitem[Schölkopf and Smola(2001)]{scholkopf_learning_2001}
B.~Schölkopf and A.~J. Smola.
\newblock \emph{Learning with Kernels: Support Vector Machines, Regularization,
  Optimization, and Beyond}.
\newblock The MIT Press, 2001.

\bibitem[Sollich(2002)]{sollich_bayesian_2002}
P.~Sollich.
\newblock {Bayesian Methods for Support Vector Machines: Evidence and
  Predictive Class Probabilities}.
\newblock \emph{Mach. Learn.}, 46:\penalty0 21--52, 2002.

\bibitem[Sriperumbudur et~al.(2011)Sriperumbudur, Fukumizu, and
  Lanckriet]{sriperumbudur_universality_2011}
B.~K. Sriperumbudur, K.~Fukumizu, and G.~R.~G. Lanckriet.
\newblock {Universality, Characteristic Kernels and RKHS Embedding of
  Measures}.
\newblock \emph{J. Mach. Learn. Res.}, 12\penalty0 (70):\penalty0 2389--2410,
  2011.

\bibitem[Teh et~al.(2006)Teh, Jordan, Beal, and Blei]{teh_hierarchical_2006}
Y.~W. Teh, M.~I. Jordan, M.~J. Beal, and D.~M. Blei.
\newblock {Hierarchical Dirichlet Processes}.
\newblock \emph{J. Amer. Statist. Assoc.}, 101\penalty0 (476):\penalty0
  1566--1581, 2006.

\bibitem[Teicher(1961)]{teicher_identifiability_1961}
H.~Teicher.
\newblock {Identifiability of Mixtures}.
\newblock \emph{Ann. Math. Statist.}, 32\penalty0 (1):\penalty0 244--248, 1961.

\bibitem[Tipping(2001)]{tipping_sparse_2001}
M.~E. Tipping.
\newblock {Sparse Bayesian Learning and the Relevance Vector Machine}.
\newblock \emph{J. Mach. Learn. Res.}, 1:\penalty0 211--244, 2001.

\bibitem[Tukey(1972)]{tukey_data_1972}
J.~W. Tukey.
\newblock {Data Analysis, Computation and Mathematics}.
\newblock \emph{Q. Appl. Math.}, 30\penalty0 (1):\penalty0 51--65, 1972.

\bibitem[Wade and In{\'a}cio(2025)]{wade2025bayesian}
S.~Wade and V.~In{\'a}cio.
\newblock {Bayesian Dependent Mixture Models: A Predictive Comparison and
  Survey}.
\newblock \emph{Stat. Sci.}, 40\penalty0 (1):\penalty0 81 -- 108, 2025.

\end{thebibliography}

\end{document}


\maketitle

{\em Organization of the supplementary material.}
The Supplementary Material contains further details on relevant examples, the proofs of our statements, both the theoretical background and the derivations of the algorithms for our numerical simulations. To ease cross-reading between the main manuscript and the supplement, here we use the prefix SM for the numbering of results, sections, and definitions (e.g., Proposition SM1, Section SM1, Equation (SM1)).

\bigskip

We first recap some useful notions on RKHS that we will use repeatedly in our proofs. A measurable kernel $k: \XX \times \XX \rightarrow \RR$ is bounded, symmetric, and positive-definite if for some $K > 0$, we have $k(x,y) = k(y,x) \leq K$ for any $x,y \in \XX$, and, moreover, 
\begin{equation*} 
	\sum_{i=1}^{m}{\sum_{j=1}^{m}{a_{i}a_{j}k(x_{i},x_{j})}} \geq 0,
\end{equation*}
for any $m \in \NN$, any $x_{1},\dots,x_{m} \in \XX$ and $a_{1},\dots,a_{m} \in \RR$. The reproducing property of a RKHS guarantees that for any $h \in \HH_{k}$,  $h(x) = \langle h(\cdot), k(x,\cdot) \rangle _{\HH_{k}}$.

Any kernel $k$ induces a squared pseudo-metric
\begin{equation} \label{eq:d_k}
    d_{k}^{2}(x,y) := k(x,x) - 2k(x,y) + k(y,y),    
\end{equation}
which is a squared distance whenever the feature map $x \mapsto k(x,\cdot)$ is injective, that is, when $k$ is injective. It is often useful to rewrite the integral expressions appearing in the kernel covariance as
\begin{equation} \label{eq:variance_distance}
    \int{k(x,x)\diff P_{0}(x)} - \iint{k(x,y)\diff P_{0}(x)\diff P_{0}(y)} = \frac{1}{2} \iint d_k^2(x,y) \diff P_0(x) \diff P_0(y).  
\end{equation}

\section{Examples} \label{sec:examples}

\cref{par_ex} Consider $X_{i,j} | \theta_{1},\theta_{2} \stackrel{\iid}{\sim} \mathcal{N}\left(\theta_{i},s^{2}\right)$ for $j \in \NN$ and $i = 1,2$, with $(\theta_1,\theta_2) \sim  \mathcal{N}\left(\boldsymbol{0},\tau^{2}\Sigma\right)$ for $s,\tau > 0$, where $\Sigma_{11} = \Sigma_{22} = 1$ and $\Sigma_{12} = \Sigma_{21} = \rho \in [-1,1]$.

The negative log-likelihood can be expressed, up to an additive constant, as
\[
    -\log \mathbb{P}\bigl(\boldsymbol{X}^{(n_1,n_2)} \big| \theta_1, \theta_2\bigr) = \frac{1}{2s^{2}} \sum_{j=1}^{n_1} (X_{1,j} - \theta_1)^2 + \frac{1}{2s^{2}} \sum_{j=1}^{n_2} (X_{2,j} - \theta_2)^2.
\]
Thus, by Bayes' rule, $\theta_{1},\theta_{2} \big| \boldsymbol{X}^{(n_1,n_2)}$ has density proportional to 
\begin{equation*}
    \mathbb{P}\bigl(\theta_{1},\theta_{2} \big| \boldsymbol{X}^{(n_1,n_2)}\bigr) \propto \exp \left( - \frac{1}{2s^{2}} \sum_{j=1}^{n_1} (X_{1,j} - \theta_1)^2 - \frac{1}{2s^{2}} \sum_{j=1}^{n_2} (X_{2,j} - \theta_2)^2 - \frac{1}{2} \boldsymbol{\theta}^\top \Sigma^{-1} \boldsymbol{\theta}  \right).
\end{equation*}

By standard algebraic manipulations, one can prove that $\theta_{1},\theta_{1} \big| \boldsymbol{X}^{(n_1,n_2)}$ follows a Gaussian distribution $\mathcal{N}(\boldsymbol{\theta}^*, \Sigma^*)$ with 
\begin{equation} \label{eq:sigma_star}
    \Sigma^* =
    \frac{s^{2}\tau^{2}}{s^{4} + (n_{1} + n_{2})s^{2}\tau^{2} + n_{1}n_{2}\tau^{4}(1 - \rho^{2})}
	\begin{pmatrix}
		s^{2} + n_{2}\tau^{2}(1 - \rho^{2}) & s^{2}\rho \\
		s^{2}\rho & s^{2} + n_{1}\tau^{2}(1 - \rho^{2})
	\end{pmatrix}
\end{equation}
and
\begin{equation} \label{eq:theta_star}
    \boldsymbol{\theta}^* =
    \frac{\tau^{2}}{s^{4} + (n_{1} + n_{2})s^{2}\tau^{2} + n_{1}n_{2}\tau^{4}(1 - \rho^{2})}
    \begin{pmatrix}
        (s^{2} + n_{2}\tau^{2}(1 - \rho^{2}))n_{1}\overline{X}_{1} + s^{2}\rho n_{2}\overline{X}_{2} \\
        (s^{2} + n_{1}\tau^{2}(1 - \rho^{2}))n_{2}\overline{X}_{2} + s^{2}\rho n_{1}\overline{X}_{1}  
    \end{pmatrix},
\end{equation}
where $\overline{X}_{i} = \sum_{j=1}^{n_{i}}{X_{i,j}}/n_{i}$ for $i = 1,2$ is the empirical mean.

The correlation between $\theta_1$ and $\theta_2$ a posteriori depends only on $\Sigma^*$ but not on $\boldsymbol{X}^{(n_1,n_2)}$. A direct computation gives
\begin{equation*}
     \corr\bigl(\theta_{1},\theta_{1} \big| \boldsymbol{X}^{(n_1,n_2)} = \boldsymbol{x}^{(n_{1},n_{2})}\bigr) = \frac{\Sigma^*_{12}}{\sqrt{\Sigma^*_{11}} \sqrt{\Sigma^*_{22}}} =  \frac{\rho}{\sqrt{1 + n_{1}\frac{\tau^{2}}{s^{2}}(1 - \rho^{2})}\sqrt{1 + n_{2}\frac{\tau^{2}}{s^{2}}(1 - \rho^{2})}}.  
\end{equation*}

\begin{example} \label{countex_linear_th1}
    Consider the linear kernel $k(x,y) = xy$ on $\XX = [0,2]$. For $X \sim \mathrm{Unif}_{[0,1]}$, let $\tilde{P}_{1} := \delta_{X}$ and $\tilde{P}_{2} := \delta_{X+1}$. It follows that $\tilde{P}_{1}$, $\tilde{P}_{2}$ are a.s. discrete with atomless mean measures $\E\bigl[\tilde{P}_{1}\bigr] = \mathrm{Unif}_{[0,1]}$ and $\E\bigl[\tilde{P}_{2}\bigr] = \mathrm{Unif}_{[1,2]}$. However, $\corr_{k}\bigl(\tilde{P}_{1},\tilde{P}_{2}\bigr) = \corr(X,X+1) = 1$.
\end{example}

\begin{example}
    Let $\tilde{P}_{i} = \omega_{i}\delta_{x_{i}} + (1 - \omega_{i})\tilde{P}$, where $\tilde{P}$ is an a.s. discrete random probability, $\omega_{i} \in (0,1)$, and $x_{1} \neq x_{2} \in \XX$, for $i = 1,2$. Note that $\tilde{P}_{1}$ and $\tilde{P}_{2}$ have fixed jumps at deterministic points and, thus, do not satisfy the assumptions of \cref{index_equal1_post}. Since $\var_{k}\bigl(\tilde{P}_{i}\bigr) = (1 - \omega_{i})^{2}\var_{k}\bigl(\tilde{P}\bigr)$ and $\cov_{k}(\tilde{P}_{1},\tilde{P}_{2}) = (1 - \omega_{1})(1 - \omega_{2})\var_{k}\bigl(\tilde{P}\bigr)$, it follows that $\corr_{k}(\tilde{P}_{1},\tilde{P}_{2}) = 1$.
\end{example}

\begin{example} \label{countex_th1}
    For $X \sim \mathrm{Unif}_{[0,1]}$, let $\tilde{P}_{1} := \delta_{X}$ and $\tilde{P}_{2} := \delta_{1-X}$, which have atomless mean measure. If we take $A = [1/4,3/4]$, then $\tilde{P}_{1}(A) = \tilde{P}_{2}(A)$ a.s.. Thus, $\corr\bigl(\tilde{P}_{1}(A),\tilde{P}_{2}(A)\bigr) = 1$.
\end{example}

\begin{example} \label{cont_ex}
    Let $W \sim \mathrm{Unif}_{[0,1]}$ and let $P \in \mathcal{P}(\XX)$. For any closed set $A$ there exists a sequence $(x_{n})_{n \in \NN} \subset A^c$ such that $x_{n} \rightarrow x_{\infty} \in A$ as $n \rightarrow \infty$. The sequence of random probabilities $\tilde{P}_{n} := W\delta_{x_{n}} + (1 - W)P$ converges weakly to $\tilde{P}_{\infty} := W\delta_{x_{\infty}} + (1 - W)P$ a.s.. If $k$ is continuous and bounded, the Dominated Convergence Theorem easily implies that $\corr_{k}\bigl(\tilde{P}_{n},\tilde{P}_{\infty}\bigr) \rightarrow 1$ as $n \rightarrow +\infty$, as one would expect. However, by taking the set-wise kernel for $A$, \cref{post_ex} shows that $\corr\bigl(\tilde{P}_{n}(A),\tilde{P}_{\infty}(A)\bigr) = -1$ for every $n\in \NN$.
\end{example}

\begin{example} \label{sec:corr_hdp}
    We first focus on the evaluation of the kernel correlation a priori. To compute the variance we observe that since $\tilde P_0 \sim \mathrm{DP}(c_0, P_0)$, $\tilde{P}_0(A)\sim \mathrm{Beta}(cP_0(A), c(1-P_0(A)))$. This implies that $\var\bigl(\tilde{P}_0(A)\bigr) = P_0(A)(1-P_0(A))/(1+c_0)$. To find the variance of $\tilde P_i(A)$, we apply the law of total variance conditioning on $\tilde{P}_0$. With some algebraic manipulations, we obtain
    \[
        \var\bigl(\tilde P_i(A)\bigr) = \frac{1+c+c_0}{(1+c)(1+c_0)}P_0(A)(1-P_0(A)).
    \]
    Thus the hDP a priori satisfies \cref{var_set_ker} with $\lambda_{1} = \lambda_{2} = (1 + c +c_0)(1+c)^{-1}(1+c_0)^{-1}$. With similar calculations based on the law of total covariance and conditional independence,
    \[
        \cov\bigl(\tilde P_1(A), \tilde P_2(A)\bigr) = \var \bigl(\tilde P_0(A)\bigr) =  \frac{P_0(A)(1-P_0(A))}{1+c_0}.
    \]
    This proves that the hDP also satisfies \cref{cov_set_ker} with $\eta = (1 + c_0)^{-1}$. Thus by \cref{corr_set_ker},
    \[
        \corr_{k}\bigl(\tilde{P}_{1},\tilde{P}_{2}\bigr) = \corr\bigl(\tilde{P}_{1}(A),\tilde{P}_{2}(A)\bigr) = \frac{1+c}{1+c+c_{0}},
    \]
    for any injective kernel $k$ and any measurable set $A$ such that $P_{0}(A) \notin \{0,1\}$. 
\end{example}

\section{Proofs of the Statements}

\subsection*{Proof of \cref{int_ker_cov_var}}

Since $\E\bigl[\mu_{k}(\tilde{P}_{i})\bigr] = \mu_{k}(P_{0,i})$ for $i = 1,2$, we have 
\[
    \cov_{k}\bigl(\tilde{P}_{1},\tilde{P}_{2}\bigr) = \E\left[\bigl\langle \mu_{k}(\tilde{P}_{1}), \mu_{k}(\tilde{P}_{2}) \bigr\rangle_{\HH_{k}}\right] - \bigl\langle\mu_{k}(P_{0,1}), \mu_{k}(P_{0,2}) \bigr\rangle_{\HH_{k}}.
\]
Thus, the result derives from the application of \cref{reproduce_kernel_int}.

\subsection*{Proof of \cref{lemma:var_equal_0}}

Note that 
\begin{equation*}
    \var_k(\tilde{P}) = \E\left[\bigl\|\mu_k(\tilde{P}) - \E\bigl[\mu_k(\tilde{P})\bigr]\bigr\|_{\HH_k}^2\right].
\end{equation*}
Thus $\var_k(\tilde{P}) \geq 0$ with equality if and only if $\bigl\|\mu_k(\tilde{P}) - \E\bigl[\mu_k(\tilde{P})\bigr]\bigr\|_{\HH_k}^2 = 0$ almost surely, that is, $\mu_k(\tilde{P}) = \E\bigl[\mu_k(\tilde{P})\bigr]$ almost surely. If $k$ is characteristic, and as $\mu_k$ is linear, it happens if and only if $\tilde{P} = \E[\tilde{P}]$ almost surely, which means $\tilde{P}$ is deterministic.

\subsection*{Proof of \cref{extreme_case_indep}}

The correlation belongs to $[-1,1]$ from an application of Cauchy-Schwarz: 
\begin{align*}
    \cov_{k}^{2}\bigl(\tilde{P}_{1},\tilde{P}_{2}\bigr) & = \E\left[\bigl\langle\mu_{k}(\tilde{P}_{1}) - \E\bigl[\mu_{k}(\tilde{P}_{1})\bigr], \mu_{k}(\tilde{P}_{2}) - \E\bigl[\mu_{k}(\tilde{P}_{2})\bigr] \bigr\rangle_{\HH_{k}}\right]^{2} \\
    & \leq \E\left[\bigl\|\mu_{k}(\tilde{P}_{1}) - \E\bigl[\mu_{k}(\tilde{P}_{1})\bigr]\bigr\|_{\HH_k}^2\right] \E\left[\bigl\|\mu_{k}(\tilde{P}_{2}) - \E\bigl[\mu_{k}(\tilde{P}_{2})\bigr]\bigr\|_{\HH_k}^2\right] \\
    & = \var_k\bigl(\tilde{P}_1\bigr) \var_k\bigl(\tilde{P}_2\bigr).
\end{align*}
If $\tilde{P}_1$ and $\tilde{P}_2$ are independent, so are $\mu_{k}(\tilde{P}_{1})$ and $\mu_{k}(\tilde{P}_{2})$, and thus $\E\bigl[\langle\mu_{k}(\tilde{P}_{1}),\mu_{k}(\tilde{P}_{2})\rangle_{\HH_k}\bigr] = \bigl\langle \E\bigl[\mu_{k}(\tilde{P}_{1})\bigr],\E\bigl[\mu_{k}(\tilde{P}_{2})\bigr] \bigr\rangle_{\HH_k}$.

\subsection*{Proof of \cref{linear_corr}}

Let $(\HH, \langle \cdot, \cdot \rangle_{\HH})$ be a Hilbert space and let $X$, $Y$ be $\HH$-valued random variables with finite first and second moments. Without loss of generality, we assume that $X$ and $Y$ are centred. Then $\corr_{\HH}(X,Y) = \pm1$ if and only if $\E[\langle X,Y \rangle_{\HH}]^{2} = \E\bigl[\|X\|_{\HH}^{2}\bigr]\E\bigl[\|Y\|_{\HH}^{2}\bigr]$, which means that there is equality in Cauchy-Schwarz. Thus $X$ and $Y$ are collinear, in the sense that $X = \alpha Y$ a.s. for $\alpha = \E[\langle X, Y \rangle]/\E\bigl[\|Y\|_{\HH}^{2}\bigr]$. In particular, $\alpha$ has the same sign as $\corr_{\HH}(X,Y)$.  

Applying this reasoning to $X = \mu_k(\tilde{P}_1)$ and $Y = \mu_k(\tilde{P}_2)$, we obtain
\[
\mu_{k}(\tilde{P}_{1}) - \E\bigl[\mu_{k}(\tilde{P}_{1})\bigr] = \alpha\bigl(\mu_{k}(\tilde{P}_{2}) - \E\bigl[\mu_{k}(\tilde{P}_{2})\bigr]\bigr) \qquad \as
\]
for some $\alpha \in \RR \setminus \{0\}$ having the same sign as $\corr_{\HH_k}\bigl(\mu_{k}(\tilde{P}_{1}),\mu_{k}(\tilde{P}_{2})\bigr)$. Since the kernel mean embedding is linear and injective for bounded signed measures as $k$ is $c_{0}$-universal, the equality reads $\tilde{P}_{1} - \E\bigl[\tilde{P}_{1}\bigr] = \alpha\bigl(\tilde{P}_{2} - \E\bigl[\tilde{P}_{2}\bigr]\bigr)$ a.s..

\subsection*{Proof of \cref{index_equal1_prior}}

It is trivial to prove that $\corr_{k}\bigl(\tilde{P}_{1},\tilde{P}_{2}\bigr) = 1$ whenever $\tilde{P}_{1} = \tilde{P}_{2}$ a.s..

Let us prove the converse implication. By \cref{linear_corr}, there exists some $\alpha > 0$ such that $\tilde{P}_{1} - \E\bigl[\tilde{P}_{1}\bigr] = \alpha \bigl(\tilde{P}_{2} - \E\bigl[\tilde{P}_{2}\bigr]\bigr)$ a.s., which we rewrite as
\[
	\tilde{P}_{1} - \alpha\tilde{P}_{2} = \E\bigl[\tilde{P}_{1}\bigr] - \alpha\E\bigl[\tilde{P}_{2}\bigr] \qquad \as.
\]
Since the random measure on the left is a.s. discrete, while the measure on the right is atomless, both sides must be a.s. null. In other words, $\tilde{P}_{1} = \alpha\tilde{P}_{2}$ a.s.. Evaluating both sides of this equality at $\XX$, we obtain that $\alpha = 1$ as $\tilde{P}_{1}(\XX) = \tilde{P}_{2}(\XX) = 1$ a.s..

\subsection*{Proof of \cref{index_equal1_post}}

It is trivial to prove that $\corr_{k}\bigl(\tilde{P}_{1},\tilde{P}_{2}\bigr) = 1$ whenever $\tilde{P}_{1} = \tilde{P}_{2}$ a.s..

Let us now prove the converse implication. Firstly, let us notice that the mean measures of $\tilde{P}_{1}$ and $\tilde{P}_{2}$, as any other probability measure, can be decomposed as
\[
	\E\bigl[\tilde{P}_{i}\bigr] = \omega_{i}Q_{i}^{\rm d} + (1 - \omega_{i})Q_{i}^{\rm a},
\]
where $\omega_{i} \in [0,1]$, $Q_{i}^{\rm d}$ is a discrete probability measure on $\XX$, and $Q_{i}^{\rm a}$ is an atomless probability measure on $\XX$. By \cref{linear_corr}, there exist $\alpha > 0$ such that $\tilde{P}_{1} - \E\bigl[\tilde{P}_{1}\bigr] = \alpha \bigl(\tilde{P}_{2} - \E\bigl[\tilde{P}_{2}\bigr]\bigr)$ a.s., which can be rewritten as 
\[
    \tilde{P}_{1} - \omega_{1}Q_{1}^{\rm d} - \alpha\bigl(\tilde{P}_{2} - \omega_{2}Q_{2}^{\rm d}\bigr) = (1 - \omega_{1})Q_{1}^{\rm a} - \alpha (1 - \omega_{2})Q_{2}^{\rm a} \qquad \as.
\]
Since the random measure on the left is a.s. discrete, while the measure on the right is atomless, both sides must be a.s. null. In particular, focusing on the left-hand side, we have 
\[ 
\tilde{P}_{1} = \alpha\tilde{P}_{2} + \xi \qquad \as,
\]
where $\xi :=  \omega_{1}Q_{1}^{\rm d} - \alpha\omega_{2}Q_{2}^{\rm d}$ is a discrete bounded signed measure. Let $z$ be any atom of $\xi$. If $\xi(\{ z \}) > 0$ then $\tilde{P}_{1}(\{ z \}) \geq \xi(\{ z \}) > 0$ a.s.: it implies that $z$ is a fixed atom of $\tilde{P}_1$, but it contradicts that $\mathbb{P}\bigl(\tilde{P}_{1}(\{  z \}) < \varepsilon\bigr) > 0$ for any $\varepsilon > 0$. On the other hand if $\xi(\{ z \}) < 0$ then $\tilde{P}_{2}(\{ z \}) \geq -  \xi(\{ z \})/\alpha > 0$, and again it implies that $z$ is a fixed atom of $\tilde{P}_2$, but it contradicts that $\mathbb{P}\bigl(\tilde{P}_{2}(\{  z \}) < \varepsilon\bigr) > 0$ for any $\varepsilon > 0$. Thus $\xi(\{ z \}) = 0$ for any $z$, which implies that $\xi = 0$ as it is a purely discrete measure. From $\xi = 0$ it is easy to conclude that $\tilde{P}_1 = \alpha \tilde{P}_2$ a.s., thus $\tilde{P}_1 = \tilde{P}_2$ as they are both probability measures a.s..

\subsection*{Proof of \cref{prop_cov_obs}}

By linearity of the mean kernel embedding, $\mu_k(\tilde{P}_i) = \E\bigl[k(X_{i,1},\cdot) \big| \tilde{P}_i\bigr]$. Moreover $X_{1,1}$ and $X_{2,1}$ are independent given $\tilde{P}_1, \tilde{P}_2$. Using the law of total covariance
\begin{align*}
    \cov_{\HH_k}(k(X_{1,1},\cdot), k(X_{2,1},\cdot)) & = \cov_{\HH_k}\left(\E\bigl[ k(X_{1,1},\cdot) \big| \tilde{P}_1, \tilde{P}_2\bigr], \E\bigl[ k(X_{2,1},\cdot) \big| \tilde{P}_1, \tilde{P}_2\bigr]\right) \\
    & = \cov_{\HH_k}\bigl(\mu_k(\tilde{P}_1), \mu_k(\tilde{P}_2)\bigr) \\
    & = \cov_{k}\bigl(\tilde{P}_{1},\tilde{P}_{2}\bigr).
\end{align*}
For the second term, we use in a similar way that $\mu_k(\tilde{P}_i) = \E\bigl[k(X_{i,1},\cdot) \big| \tilde{P}_i\bigr] = \E\bigl[k(X_{i,2},\cdot) \big| \tilde{P}_i\bigr]$ and $X_{i,1}$ and $X_{i,2}$ are independent given $\tilde{P}_i$. Thus
\begin{align*}
    \cov_{\HH_k}(k(X_{i,1},\cdot), k(X_{i,2},\cdot)) & = \cov_{\HH_k}\left(\E\bigl[ k(X_{i,1},\cdot) \big| \tilde{P}_i\bigr], \E\bigl[ k(X_{i,2},\cdot) | \tilde{P}_i\bigr]\right) \\
    & = \cov_{\HH_k}\bigl(\mu_k(\tilde{P}_i), \mu_k(\tilde{P}_i)\bigr) \\
    & = \var_{k}\bigl(\tilde{P}_{i}\bigr).
\end{align*}

\subsection*{Proof of \cref{prop_cov_unbiased}}

If $Y, Z$ are two random variables valued in a Hilbert space $\HH$ and $\bigl(Y^{(t)}, Z^{(t)}\bigr)_{t=1}^M$ are i.i.d. samples from them, then an unbiased estimator of the covariance between $Y$ and $Z$ is 
\begin{equation*}
    \widehat{\cov}_{\HH,M}(Y,Z) = \frac{1}{M-1} \sum_{t=1}^M \bigl\langle Y^{(t)} - \bar{Y}, Z^{(t)} - \bar{Z} \bigr\rangle_\HH, \qquad \bar{Y} = \frac{1}{M} \sum_{m=1}^M Y^{(t)} \bar{Z} = \frac{1}{M} \sum_{m=1}^M Z^{(t)}.
\end{equation*}
Expanding the formula, this expression can be rewritten as
\begin{align}
    \widehat{\cov}_{\HH,M}(Y,Z) &=  \frac{1}{M-1} \sum_{t=1}^M \bigl\langle Y^{(t)},Z^{(t)} \bigr\rangle_\HH - \frac{M}{M-1} \bigl\langle \bar{Y},\bar{Z} \bigr\rangle_{\HH}
    \label{eq:aux_inner_prod_diff} \\
    &= \frac{1}{M-1} \sum_{t=1}^M \bigl\langle Y^{(t)}, Z^{(t)}\bigr\rangle_\HH - \frac{1}{M(M-1)}  \sum_{t=1}^M \sum_{s=1}^M \bigl\langle Y^{(t)},Z^{(s)} \bigr\rangle_\HH.
    \label{eq:aux_inner_prod_U_stats}
\end{align}

For the unbiased estimator of $\cov_{k}\bigl(\tilde{P}_{1},\tilde{P}_{2}\bigr)$, we apply \cref{eq:aux_inner_prod_U_stats} with $Y = k(X_{1,1},\cdot)$ and $Z = k(X_{2,1},\cdot)$, as the covariance between $Y$ and $Z$ is the kernel covariance (see \cref{prop_cov_obs}). We recover the statement thanks to the reproducing property $\langle k(x,\cdot),k(y,\cdot) \rangle_{\HH_k} = k(x,y)$ for any $x,y$. 

For the unbiased estimator of $\var_{k}\bigl(\tilde{P}_{i}\bigr)$ we use \cref{eq:aux_inner_prod_U_stats} with $Y = k(X_{i,1},\cdot)$ and $Z = k(X_{i,2},\cdot)$.

\subsection*{Proof of \cref{prop_cor_asymp_normal}}

The proof relies on the delta method. We only sketch it, as we do not aim to find the exact formula for the asymptotic covariance. Let us introduce the random vector $\mathbf{Z}$ in $\RR^3 \times \HH_k^4$: 
\begin{equation*}
    \mathbf{Z} = 
    \frac{1}{M} \sum_{t=1}^M \mathbf{Z}^{(t)}, \qquad \text{with} \quad \mathbf{Z}^{(t)} =
    \begin{pmatrix}
        k\bigl(X^{(t)}_{1,1},X^{(t)}_{2,1}\bigr) \\
        k\bigl(X^{(t)}_{1,1},X^{(t)}_{1,2}\bigr) \\
        k\bigl(X^{(t)}_{2,1},X^{(t)}_{2,2}\bigr) \\
        k\bigl(X^{(t)}_{1,1}, \cdot \bigr) \\
        k\bigl(X^{(t)}_{1,2}, \cdot \bigr) \\
        k\bigl(X^{(t)}_{2,1}, \cdot \bigr) \\
        k\bigl(X^{(t)}_{2,2}, \cdot \bigr)
    \end{pmatrix}.
\end{equation*}
By the central limit theorem for Hilbert spaces (see e.g. \cite{hoffmann_law_1976}), and as $k$ is bounded, this vector is asymptotically normal, in the sense that $\sqrt{M}(\mathbf{Z} - \E[\mathbf{Z}])$ converges to a normal distribution. On the other hand, from the formula \cref{eq:aux_inner_prod_diff} in the proof of \cref{prop_cov_unbiased}, 
\begin{equation*}
    \widehat{\corr}_{k,M}\bigl(\tilde{P}_1,\tilde{P}_2\bigr) = \frac{\mathbf{Z}_1 - \langle \mathbf{Z}_4, \mathbf{Z}_6 \rangle_{\HH_k}}{\sqrt{ \mathbf{Z}_2 - \langle\mathbf{Z}_4, \mathbf{Z}_5\rangle_{\HH_k}} \sqrt{\mathbf{Z}_3 - \langle \mathbf{Z}_6, \mathbf{Z}_7\rangle_{\HH_k}}}.
\end{equation*}
As the inner product is Fréchet differentiable in $\HH_k$, thus Hadamard differentiable, we can apply the delta method \cite[Theorem 3.9.4]{vandervaart_weak_1996} to the differentiable function $\phi(\mathbf{z}) = (\mathbf{z}_1 - \langle\mathbf{z}_4, \mathbf{z}_6 \rangle_{\HH_k}) (\mathbf{z}_2 - \langle \mathbf{z}_4, \mathbf{z}_5\rangle_{\HH_k})^{-1/2} (\mathbf{z}_3 - \langle \mathbf{z}_6, \mathbf{z}_7\rangle_{\HH_k})^{-1/2}$, and conclude that $\sqrt M(\phi(\mathbf{Z}) - \phi(\E[\mathbf{Z}]))$ converges to a Gaussian distribution, which is the conclusion of the theorem.

\subsection*{Proof of \cref{setwise_kernel_corr}}

The kernel $k$ is symmetric and bounded by definition. To prove that it is positive semi-definite consider $x_{1},\dots,x_{m} \in \XX$ and $a_{1},\dots,a_{m} \in \RR$. Firstly, we notice that for $i,j \in \{1,\dots,m\}$, $k_{A}(x_{i},x_{j}) = 1$ if $x_i,x_j \in A$ and $k_{A}(x_{i},x_{j}) = 0$ otherwise. Hence, if we denote $I_{A} = \{i \in \{1,\dots,m\} \text{ s.t. } x_{i} \in A\}$, 
\[
    \sum_{i=1}^{m}{\sum_{j=1}^{m}{a_{i}a_{j}k(x_{i},x_{j})}} = \sum_{i \in I_{A}}{\sum_{j \in I_{A}}{a_{i}a_{j}}} = \bigg(\sum_{i \in I_{A}}{a_{i}} \bigg)^{2} \geq 0.
\]
Let us write $P_{0,1} := \E\bigl[\tilde{P}_{1}\bigr]$ and $P_{0,2} := \E\bigl[\tilde{P}_{2}\bigr]$. Then, by \cref{int_ker_cov_var} and Fubini's Theorem, it holds
\begin{align*}
\cov_{k_A}\bigl(\tilde{P}_{1},\tilde{P}_{2}\bigr) &= \E \bigg[ \iint{\mathbbm{1}_{A}(x)\mathbbm{1}_{A}(y)\diff\tilde{P}_{1}(x)\diff\tilde{P}_{2}(y)} \bigg] - \iint{\mathbbm{1}_{A}(x)\mathbbm{1}_{A}(y)\diff P_{0,1}(x)\diff P_{0,2}(y)} \\
	&= \E\bigl[\tilde{P}_{1}(A)\tilde{P}_{2}(A)\bigr] - P_{0,1}(A)P_{0,2}(A) = \cov\bigl(\tilde{P}_{1}(A),\tilde{P}_{2}(A)\bigr).
\end{align*}

\subsection*{Proof of \cref{cov_set_ker}}

By taking $k(x,y) = k_{A}(x,y) = \mathbbm{1}_{A}(x)\mathbbm{1}_{A}(y)$, \cref{setwise_kernel_corr} guarantees that $(ii)$ implies $(i)$. To show the converse, we prove that $(ii)$ holds for an increasingly larger class of kernels.

\emph{Step 1.} We observe that by linearity, $(ii)$ holds for any kernel function $k$ that can be written as a linear combination of kernel functions $k_A$.

\emph{Step 2.} For a measurable function $f : \XX \times \XX \to \RR$, let us define its ``symmetrized'' version $\mathcal{S}[f]$ by $\mathcal{S}[f](x,y) = (f(x,y) + f(y,x))/2$. If $A,B$ measurable then $(ii)$ holds for $k = \mathcal{S}[\mathbbm{1}_{A \times B}]$. This comes from linearity as
\[  \mathcal{S}[\mathbbm{1}_{A  \times B}] = \frac{1}{2} ( k_{A \cup B} + k_{A \cap B} - k_{A \setminus B} - k_{B \setminus A}). \]

\emph{Step 3.} Thanks to the monotone class lemma, the property $(ii)$ holds for any $\mathcal{S}[f]$ where $f$ is the indicator function of a measurable set in $\mathcal{X} \otimes \mathcal{X}$. Then from stability by linear combination and the monotone convergence theorem, $(ii)$ holds for any $\mathcal{S}[f]$, where $f$ is a measurable and bounded function over $\XX \times \XX$. The conclusion follows as $k = \mathcal{S}[k]$ since $k$ is symmetric.
 
Lastly, we use $(i)$ to prove $\eta \in [-1,1]$. Indeed, we have 
\[
    \eta = \frac{\cov\bigl(\tilde{P}_{1}(A),\tilde{P}_{2}(A)\bigr)}{P_{0}(A) - P_{0}(A)^{2}} = \frac{\E\bigl[\tilde{P}_{1}(A)\tilde{P}_{2}(A)\bigr] - P_{0}(A)^{2}}{P_{0}(A) - P_{0}(A)^{2}}.
\]
On the one hand, from $\tilde{P}_{2}(A) \leq 1$ we see that $\eta \leq \bigl(\E\bigl[\tilde{P}_{1}(A)\bigr] - P_{0}(A)^{2}\bigr)/(P_{0}(A) - P_{0}(A)^{2}) = 1$. On the other hand,
\[
    \eta \geq - \frac{P_{0}(A)^{2}}{P_{0}(A) - P_{0}(A)^{2}} = - \frac{P_0(A)}{P_0(A^c)},
\]
being $A^c$ the complement of $A$. Up to exchanging the role of $A$ and $A^c$, we can always find a set $A$ such that $P_0(A) \leq P_0(A^c)$, hence $\eta \geq -1$.

\subsection*{Proof of \cref{var_set_ker}}

We apply \cref{cov_set_ker} to $\tilde{P}_1 = \tilde{P}_2 = \tilde{P}$. The only new element is that $\lambda \geq 0$, which is easily deduced from the fact that the variance is always non-negative, see \cref{lemma:var_equal_0}.

\subsection*{Proof of \cref{corr_set_ker}}

The result comes from a direct application of \cref{cov_set_ker} and \cref{var_set_ker}, provided that both variances are non-null for both correlations.

Firstly, $\var\bigl(\tilde{P}_{1}(A)\bigr),\var\bigl(\tilde{P}_{2}(A)\bigr) \neq 0$ since $\tilde{P}_{i}(A)$ is not deterministic for $i = 1,2$.

Secondly, $\var_{k}(\tilde{P}_{1}),\var_{k}(\tilde{P}_{2}) \neq 0$. If by contradiction,
\[
    \int{k(x,x) \diff P_{0}(x)} - \iint{k(x,y) \diff P_{0}(x) \diff P_{0}(x)} = \frac{1}{2}\iint{d_{k}^{2}(x,y) \diff P_{0}(x) \diff P_{0}(x)} = 0,
\]
implies that $d_{k}^{2}(x,y)$ for $P_{0} \times P_{0}$-almost every $x,y \in \XX$. Now, $d_{k}^{2}(x,y)$ is a distance since $k$ is injective. Hence, $x = y$ for $P_{0} \times P_{0}$-almost every $x,y \in \XX$. This implies that there exists $z \in \XX$ such that $P_{0} = \delta_{z}$ almost surely. Hence, $\tilde{P}_{1} = \tilde{P}_{2} = \delta_{z}$, which contradicts the assumption that $\tilde{P}_{1}$ and $\tilde{P}_{2}$ are non-determinist almost surely.

\subsection*{Proof of \cref{mix_corr}}

Firstly, $k_f$ is symmetric due to the symmetry of $k$. With Fubini's theorem and given the definition of $k_f$, for any signed bounded measures $\xi_1$, $\xi_2$ on $\XX$, 
\begin{align}
    \iint_{\YY \times \YY} &{k(y_1,y_2) f_{\xi_1}(y_1)  f_{\xi_2}(y_2)} \diff y_1 \diff y_2
    \notag \\
    &=
    \iint_{\XX \times \XX}{\bigg(\iint_{\YY \times \YY}{k(y_1,y_2) f(y_1;x_1) f(y_2;x_2)} \diff y_1 \diff y_2 \bigg)\diff \xi_1(x_1)\diff \xi_2(x_2)}
    \notag \\
    & = \iint_{\XX \times \XX} k_f(x_1,x_2) \, \diff \xi_1(x_1) \diff \xi_2(x_2). 
    \label{eq:mixture_kernel}
\end{align}
Applying this formula with $\xi_1 = \xi_2$, we see that the last expression in \cref{eq:mixture_kernel} is non-negative (because $k$ is positive definite), and thus $k_f$ is positive definite.

Then, let us write $P_{0,1} := \E\bigl[\tilde{P}_{1}\bigr]$ and $P_{0,2} := \E\bigl[\tilde{P}_{2}\bigr]$. Thus, $\E\bigl[f_{\tilde{P}_i}\bigr] = f_{P_{0,i}}$. Applying \cref{eq:mixture_kernel} with $(\xi_1,\xi_2) = (\tilde{P}_1,\tilde{P}_2)$ and then $(P_{0,1},P_{0,2})$, we have 
\begin{align*}
    \E \left[ \iint_{\YY \times \YY} {k(y_1,y_2) f_{\tilde{P}_1}(y_1)  f_{\tilde{P}_2}(y_2)} \diff y_1 \diff y_2 \right] & = \E \left[ \iint_{\XX \times \XX} k_f(x_1,x_2) \, \diff \tilde{P}_1(x_1) \diff \tilde{P}_2(x_2) \right], \\
    \iint_{\YY \times \YY} {k(y_1,y_2) f_{P_{0,1}}(y_1)  f_{P_{0,2}}(y_2)} \diff y_1 \diff y_2  & = \iint_{\XX \times \XX} k_f(x_1,x_2) \, \diff P_{0,1}(x_1) \diff P_{0,2}(x_2).
\end{align*}
Hence by subtracting these two equalities we obtain $\cov_{k}\bigl(f_{\tilde{P}_1},f_{\tilde{P}_2}\bigr) = \cov_{k_f}\bigl(\tilde{P}_1,\tilde{P}_2\bigr)$.

For the second part of the statement, we assume identifiability of the family $\{ f(\cdot;x) \ : \ x \in \XX \}$ and that $k$ is characteristic, and we want to show that $k_f$ is characteristic. By \cref{lm:equi_char} below, which is an easy variant of \citet[Proposition 4]{sriperumbudur_universality_2011}, we need to show that $k_f$ is conditionally integrally strictly positive definite. Take $\xi \in \mathcal{M}_b(\mathbb{X})$ with $\xi(\mathbb{X}) = 0$, assume $\iint k_f(x_1,x_2) \diff \xi(x_1) \diff \xi(x_2) = 0$ and we want to show $\xi = 0$. With \cref{eq:mixture_kernel} 
\begin{equation*}
    \iint_{\XX \times \XX} k_f(x_1,x_2) \, \diff \xi(x_1) \diff \xi(x_2) = \iint_{\YY \times \YY} {k(y_1,y_2) f_{\xi}(y_1)  f_{\xi}(y_2)} \diff y_1 \diff y_2,
\end{equation*}
and, thus, if the left-hand side vanishes, so does the right-hand side. By Fubini, we can check that $f_\xi(y) \diff y$ is a measure with zero total mass, and as $k$ itself is characteristic, \cref{lm:equi_char} implies that, as a measure, $f_{\xi}(y) \diff y = 0$. We write $\xi = \xi_+ - \xi_-$ for the Hahn-Jordan decomposition of $\xi$ with $\xi_\pm$ non-negative measures.  Calling $ a = \xi_+(\XX) = \xi_-(\XX)$ and $P_1 = \xi_+/a$, $P_2 = \xi_- / a$, we write $\xi = a(P_1 - P_2)$ with $P_1, P_2$ probability distributions over $\XX$. By linearity, $f_\xi = a (f_{P_1} - f_{P_2})$, but we already know that $f_\xi = 0$ a.e., thus, $a = 0$ or $f_{P_1} = f_{P_2}$ a.e. which implies $P_1 = P_2$ by identifiability. It means, in any case, that $\xi = 0$, and the proof is concluded.

\begin{lemma} \label{lm:equi_char}
    Let $k$ be a bounded kernel on a space $\mathbb{X}$. Then $k$ is characteristic if and only if it is conditionally integrally strictly positive definite, that is, if and only if, for any $\xi \in \mathcal{M}_b(\mathbb{X})$ with $\xi(\mathbb{X}) = 0$ and $\xi \neq 0$, we have
    \begin{equation*}
        \iint k(x_1,x_2) \, \diff \xi(x_1) \diff \xi(x_2) > 0.
    \end{equation*}
\end{lemma}

\begin{proof}
    From the reproducing property \cref{reproduce_kernel_int}, we have 
    \begin{equation*}
        \iint k(x_1,x_2) \, \diff \xi(x_1) \diff \xi(x_2) = \| \mu_k(\xi) \|_{\HH_k}^2,
    \end{equation*}
    and thus the left-hand side vanishes if and only if $\mu_k(\xi) = 0$.
    
    First, assume that $k$ is characteristic. By the Hahn-Jordan decomposition any $\xi$ with $\xi(\mathbb{X}) = 0$ but $\xi \neq 0$ can be written as $\xi = a(P_1-P_2)$ with $a > 0$ for $P_1$, $P_2$ two distinct probability distributions. Thus, $\mu_k(P_1) \neq \mu_k(P_2)$ and by linearity $\mu_k(\xi) \neq 0$, which implies $\| \mu_k(\xi) \|_{\HH_k}^2 > 0$.
    
    Conversely, assume that $k$ is not characteristic, and let $P_1$ and $P_2$ be two distinct probability distributions with $\mu_k(P_1) = \mu_k(P_2)$. With $\xi = P_1 -P_2$, we have $\xi \neq 0$ and $\xi(\mathbb{X}) = 0$. But $\mu_k(\xi) = 0$, so $\| \mu_k(\xi) \|_{\HH_k}^2 =0$. It shows that $k$ is not conditionally integrally strictly positive definite.
\end{proof}

\subsection*{Proof of \cref{corr_set_ker_phi}}
From \cref{mix_corr}, it holds that $\corr_{k}\bigl(f_{\tilde{P}_1},f_{\tilde{P}_2}\bigr) = \corr_{k_f}\bigl(\tilde{P}_{1},\tilde{P}_{2}\bigr)$. Now, since $k_f$ is an injective kernel, the result follows by \cref{corr_set_ker}.

\subsection*{Proof of \cref{mix_trans_corr}}

The result comes directly by Fubini's Theorem, since $k_{f}(x_1,x_2)$ is equal to
\begin{align*}
    &\iint \bigg(\int  e^{- i\langle y_1 - y_2, z \rangle} \diff \nu(z) \bigg) f(y_1; x_1) f(y_2; x_2) \diff y_1 \diff y_2  \\
    &= \int \bigg(\int e^{-i\langle y_1, z \rangle} f(y_1; x_1) \diff y_1  \int e^{i\langle y_2, z \rangle} f(y_2; x_2) \diff y_2 \bigg) \diff\nu(z) = \int \hat{f}(z; x_1)\overline{\hat{f}(z; x_2)}\diff \nu(z).
\end{align*}

\subsection*{Proof of \cref{exchangeability_mixtures}}

We use $\corr_{k}\bigr(f_{\tilde{P}_{1}},f_{\tilde{P}_{2}}\bigl) = \cov_{k_f}\bigl(\tilde{P}_{1},\tilde{P}_{2}\bigr)$ from \cref{mix_corr} and then apply \cref{index_equal1_prior} to $\tilde{P}_1, \tilde{P}_2$ with kernel $k_f$.

\subsection*{Proof of \cref{convergence_corr_hDP}}

For this proof, we rely on the quasi-conjugacy property of the augmented hDP model, which is recalled below, in \cref{sec:hdp_prop}. This augmented model relies on the introduction of a specific sequence of latent random variables $\boldsymbol{T}^{(n_{1},n_{2})}$ \citep{teh_hierarchical_2006,camerlenghi_distribution_2019,catalano_unified_2023}, commonly referred to as ``tables'' in the restaurant franchise metaphor. For the sake of compactness, we write $\boldsymbol{X}$ and $\boldsymbol{T}$ instead of $\boldsymbol{X}^{(n_{1},n_{2})}$ and $\boldsymbol{T}^{(n_{1},n_{2})}$. The key result is that $(\tilde{P}_1, \tilde{P}_2)$, given $\boldsymbol{X}$ and $\boldsymbol{T}$, follows a hierarchical Dirichlet process with updated parameters, see \cref{post_hDP} below. We start with an auxiliary lemma before moving to the core of the proof.

\begin{lemma} \label{var_cov_post_hDP_tabs}
    In the augmented hDP model in \cref{post_hDP}, we have
    \[
        \cov_{k}\bigl(\tilde{P}_{1},\tilde{P}_{2} \big| \boldsymbol{X},\boldsymbol{T},\tilde{P}_{0}\bigr) = 0,    
    \]
    and, for $i=1,2$, with $\hat{P}_{i} = n_i^{-1}\sum_{j=1}^{n_i} \delta_{X_{i,j}}$,
    \begin{multline*}
        \var_{k}\bigl(\tilde{P}_{i} \big| \boldsymbol{X}, \boldsymbol{T},\tilde{P}_{0}\bigr) = \frac{1}{c + n_{i} + 1}\Bigg(\frac{1}{2} \frac{c^{2}}{(c + n_{i})^{2}}\iint{d_{k}^{2}(x,y)\diff \tilde{P}_{0}(x)\diff \tilde{P}_{0}(y)} \\  
    	+ \frac{cn_{i}}{(c + n_{i})^{2}}\iint{d_{k}^{2}(x,y)\diff\tilde{P}_{0}(x)\diff\hat P_{i}(y)} + \frac{1}{2}\frac{n_{i}^{2}}{(c + n_{i})^{2}}\iint{d_{k}^{2}(x,y)\diff\hat P_{i}(x)\diff\hat P_{i}(y)} \Bigg).
    \end{multline*}
\end{lemma}
\begin{proof}
    The covariance result follows from observing that in the augmented model \cref{post_hDP}, $\tilde{P}_1$ and $\tilde{P}_2$ are independent given $\boldsymbol{X}$, $\boldsymbol{T}$ and $\tilde{P}_0$.

    The variance result comes from the fact that, conditionally on $\boldsymbol{X},\boldsymbol{T},\tilde{P}_{0}$, as in \cref{post_hDP}, $\tilde{P}_{i}$ follows a Dirichlet Process prior with baseline measure $\tilde{P}_i^* = c/(c+n_i) \tilde{P}_0 + n_i/(c+n_i) \hat{P}_i$ and concentration parameter $c + n_{i}$. Hence,
    \[
        \var_{k}\bigl(\tilde{P}_{i} \big| \boldsymbol{X}, \boldsymbol{T},\tilde{P}_{0}\bigr) = \frac{1}{2}\frac{1}{c + n_{i} + 1}\iint{d_{k}^{2}(x,y) \diff \tilde{P}^{*}_i(x) \diff \tilde{P}^{*}_i(y)}.
    \]
    The result comes from the expansion of the integral on the right-hand side.
\end{proof}

We now move to the proof of \cref{convergence_corr_hDP}. To bound the correlation from above, we need to bound the covariance from above and bound the variance from below. We condition on the tables $\boldsymbol{T}$ and use \cref{var_cov_post_hDP_tabs} that we just proved. By the variance decomposition, for $i = 1,2$,
\begin{align*}
    \var_{k}\bigl(\tilde{P}_{i} \big| \boldsymbol{X}\bigr) &\geq \E\bigl[\var_{k}\bigl(\tilde{P}_{i} \big|\boldsymbol{X}, \boldsymbol{T}, \tilde{P}_{0}\bigr)\big| \boldsymbol{X}\bigr] \\
    &\geq  \frac{1}{2}\frac{n_{i}^{2}}{(c + n_{i} + 1)(c + n_{i})^{2}}\iint{d_{k}^{2}(x,y)\diff\hat{P}_{i}(x)\diff\hat{P}_{i}(y)}, 
\end{align*}
where, crucially, we use that $\hat{P}_i = n_i^{-1}\sum_{j=1}^{n_i} \delta_{X_{i,j}}$ is deterministic conditionally on $\boldsymbol{X}$ and does not depend on $\boldsymbol{T}$.

From the covariance decomposition, 
\begin{align*}
    0 \leq \cov_{k}\bigl(\tilde{P}_{1},\tilde{P}_{2}|\boldsymbol{X}\bigr) &= \cov_{k}\bigl(\E\bigl[\tilde{P}_{1} \big| \boldsymbol{X},\boldsymbol{T},\tilde{P}_{0}\bigr],\E\bigl[\tilde{P}_{2} \big| \boldsymbol{X},\boldsymbol{T},\tilde{P}_{0}\bigr]\big|\boldsymbol{X}\bigr).
\end{align*}
Now, recalling that for $i = 1,2$,
\begin{equation} \label{eq:expression_mean_Ptildei}
    \E\bigl[\tilde{P}_{i}|\boldsymbol{X},\boldsymbol{T},\tilde{P}_{0}\bigr] = \frac{c}{c+n_{i}}\tilde{P}_{0} + \frac{n_{i}}{c+n_{i}}\hat{P}_{i},
\end{equation}
and as $\hat{P}_i$ is deterministic conditionally on $\boldsymbol{X}$, we have by bilinearity
\[
    \cov_{k}\bigl(\tilde{P}_{1},\tilde{P}_{2} |\boldsymbol{X}\bigr) = \frac{c^{2}}{(c + n_{1})(c + n_{2})}\var_{k}\bigl(\tilde{P}_{0}|\boldsymbol{X}\bigr).
\]
With $K>0$ an upper bound on the kernel, we have $\var_{k}\bigl(\tilde{P}_{0}|\boldsymbol{X}\bigr) \leq K$ thus
\[
    \cov_{k}\bigl(\tilde{P}_{1},\tilde{P}_{2}|\boldsymbol{X}\bigr) \leq \frac{c^{2}K}{(c + n_{1})(c + n_{2})}.
\]

Putting everything together, we deduce, 
\[
    0 \leq \corr_{k}\bigl(\tilde{P}_{1},\tilde{P}_{2} | \boldsymbol{X}\bigr) \leq 2 c^2 K \prod_{i=1}^{2}\bigg({\frac{\sqrt{c + n_{i} + 1}}{n_{i}}\bigg( \iint d_{k}^{2}(x,y)\diff\hat{P}_{i}(x)\diff \hat{P}_{i}(y)\bigg)^{-1/2}}\bigg).
\]
By the assumption of non-degeneracy, the last terms in the product are larger than a strictly positive constant in the limit. The conclusion follows. 

\begin{remark}
    As can be seen in the proof, the result relies only on the conditional distribution of $(\tilde{P}_1,\tilde{P}_2)$ given $\boldsymbol{X}^{(n_{1},n_{2})}$, $\boldsymbol{T}^{(n_{1},n_{2})}$ and $\tilde{P}_0$, and not on the marginal distribution of $\boldsymbol{T}^{(n_{1},n_{2})}$ given $\boldsymbol{X}^{(n_{1},n_{2})}$.
\end{remark}

\section{Distributional Properties of the Hierarchical Dirichlet Process for Proofs and Simulations} \label{sec:hdp_prop}

The partially exchangeable model with an hDP prior \cref{hDP} is quasi-conjugate a posteriori, as shown in \cite{teh_hierarchical_2006, camerlenghi_distribution_2019}. To explain this structure, we introduce an augmented model with additional latent variables, the \emph{tables} $T_{i,j}$ \citep{teh_hierarchical_2006}, using the formulation a priori in \cite{catalano_unified_2023}.
\[
    (X_{i,j},T_{i,j}) \big| \tilde{P}_{1, \textsc{xt}},\tilde{P}_{2,\textsc{xt}} \stackrel{\iid}{\sim} \tilde{P}_{i, \textsc{xt}}, \qquad \tilde{P}_{1,\textsc{xt}},\tilde{P}_{2,\textsc{xt}} \big| \tilde{P}_{0} \stackrel{\iid}{\sim} \mathrm{DP}\bigl(c,\tilde{P}_{0} \times H\bigr), \qquad \tilde{P}_{0} \sim \mathrm{DP}(c_{0},P_{0}),
\]
for $j \in \NN$ and $i=1,2$, where $H$ is an atomless probability measure. Clearly, the sequence $\{X_{i,j}\}_{i,j}$ generated by this augmented model has the same marginal law as the partially exchangeable model with hDP prior as in  \cref{hDP}. Moreover, we observe that, marginally, $\{T_{1,j}\}_j$,$\{T_{2,j}\}_j$ are two independent exchangeable sequences directed by the same Dirichlet Process prior with concentration parameter $c$ and baseline measure $H$.

From this augmented model, following \cite{camerlenghi_bayesian_2018, camerlenghi_distribution_2019}, we can specify the joint one-step-ahead predictive distribution for the observations and the tables, which can be used to sample from the model both a priori and a posteriori. Let us denote with $X_{1}^{*},\dots,X_{K}^{*}$ the unique values in $\{X_{i,j}\}_{i,j}$. For $h \in \{ 1, \ldots, K \}$ let $\ell_{h}$ be the number of unique values in $\{ T_{i,j}: X_{i,j}=X_h^{*} \}_{i,j}$ and $|\mathbf{\ell}| := \ell_1 + \ldots + \ell_K$. Then, for $i=1,2$,
\begin{align} \label{sampling_hDP}
    \nonumber
    &\mathbb{P}\bigl(X_{i,n_i+1} \in  A, T_{i,n_i+1} \in B \big| \boldsymbol{X}^{(n_{1},n_{2})}, \boldsymbol{T}^{(n_{1},n_{2})}\bigr) = \\
    &= \sum_{j=1}^{n_{i}}{\frac{1}{c + n_{i}}\delta_{X_{i,j}}(A)\delta_{T_{i,j}}(B)} + \frac{c}{c + n_{i}}\left(\frac{c_{0}}{c_{0} + |\ell|}P_{0}(A) + \sum_{h=1}^{K}{\frac{\ell_{h}}{c_{0} + |\ell|}\delta_{X^{*}_{h}}(A)}\right)H(B),
\end{align}
where $\boldsymbol{X}^{(n_{1},n_{2})} = \bigl((X_{1,j})_{j=1}^{n_{1}},(X_{2,j})_{j=1}^{n_{2}}\bigr)$ and $\boldsymbol{T}^{(n_{1},n_{2})} = \bigl((T_{1,j})_{j = 1}^{n_{1}},(T_{2,j})_{j = 1}^{n_{2}}\bigr)$.

This formula is usually interpreted through a \textit{restaurant franchise} metaphor \citep{teh_hierarchical_2006}: two restaurants of a franchise share the same menu, made of infinitely many dishes sampled from the atomless baseline measure $P_{0}$. Each restaurant has potentially infinitely many tables and serves only one dish per table. The first customer enters one of the two restaurants, say restaurant $i$, and sits at a new table, whose label is randomly generated from $H$, and eats the unique dish served at that table. Each of the next customers entering restaurant $i$ sits at the same table as one of the other $n_{i}$ customers with probability $1/(c + n_{i})$ and, thus, eats their same dish. Alternatively, they sit at a new table, randomly generated from $H$, with probability $c/(c + n_{i})$. There, they choose one of the other $|\ell|$ tables across the franchise with probability $1/(c_{0} + |\ell|)$  and eat the same dish that is being eaten at that table. Alternatively, they eat a new dish from the menu with probability $c_{0}/(c_{0} + |\ell|)$. In this metaphor, $X_{i,j}$ represents the dish eaten by the $j$-th customer in the $i$-th restaurant, while $T_{i,j}$ is the label of the table where they sit.

\begin{remark}
    The total number of tables is $|\ell| = \ell_{1} + \dots + \ell_{K}$, and $K \leq |\ell| \leq n_{1} + n_{2}$. It is $K$ when there is a unique table for each dish. It is $n_{1} + n_{2}$ when there is one table per customer.
\end{remark}

\begin{remark} \label{compatible_TX}
    The predictive distribution above forces the sequences $\{X_{i,j}\}_{i,j}$ and $\{T_{i,j}\}_{i,j}$ to have some compatibility properties, which illustrate their dependence. Firstly, $T_{1,j_{1}} \neq T_{2,j_{2}}$ a.s. for every $j_{i} = 1,\dots,n_{i}$ for $i = 1,2$. Secondly, if $T_{i,j_{1}} = T_{i,j_{2}}$ a.s. for some $j_{i} \in \{1,\dots,n_{i}\}$,
    then $X_{i,j_{1}} = X_{i,j_{2}}$ a.s.. Both these conditions can be interpreted in light of the restaurant franchise metaphor. The former states that a table cannot be shared across different restaurants. The latter says that if two customers are seated at the same table, they must eat the same dish.
\end{remark}

For a posterior characterization of the process, we report the quasi-conjugacy result of \cite{camerlenghi_distribution_2019}. Conditionally on the latent tables,
\begin{equation} 
\label{post_hDP}
\begin{split}
    \tilde{P}_{i} | \boldsymbol{X}^{(n_{1},n_{2})}, \boldsymbol{T}^{(n_{1},n_{2})},\tilde{P}_{0}   &\stackrel{\text{ind.}}{\sim} \mathrm{DP}\bigg(c + n_{i}, \frac{c}{c + n_{i}}\tilde{P}_{0} + \frac{n_{i}}{c + n_{i}}\hat P_{i}\bigg), \\
    \tilde{P}_{0} | \boldsymbol{X}^{(n_{1},n_{2})}, \boldsymbol{T}^{(n_{1},n_{2})} &\sim \mathrm{DP}\bigg(c_{0} + |\ell|, \frac{c_{0}}{c_{0} + |\ell|} P_{0} + \frac{|\ell|}{c_{0} + |\ell|}\hat P_{0} \bigg), 
\end{split}
\end{equation}
where $\hat P_{0} = |\ell|^{-1} \sum_{h=1}^{K} \ell_{h}\delta_{X^{*}_{h}}$ and $\hat P_{i} = n_i^{-1}\sum_{j=1}^{n_i} \delta_{X_{i,j}} = n_i^{-1} \sum_{h=1}^{K}n_{i,h}\delta_{X^{*}_{h}}$, independently for $i = 1,2$. In particular, this tells us that conditionally on the tables, the posterior can be interpreted as an hDP with unequal marginals. It follows that we can reproduce the type of calculations in \cref{sec:corr_hdp} to find explicit expressions of $\cov_{k}\bigl(\tilde{P}_{1},\tilde{P}_{2}\big|\boldsymbol{X}^{(n_{1},n_{2})}, \boldsymbol{T}^{(n_{1},n_{2})},\tilde{P}_{0}\bigr)$ and $\var_{k}\bigl(\tilde{P}_{i}\big|\boldsymbol{X}^{(n_{1},n_{2})}, \boldsymbol{T}^{(n_{1},n_{2})},\tilde{P}_{0}\bigr)$. These are key to the proof of \cref{convergence_corr_hDP} and can be found in \cref{var_cov_post_hDP_tabs}.

To conclude the posterior characterization, we need to provide the conditional law of $\boldsymbol{T}^{(n_{1},n_{2})}$ given $\boldsymbol{X}^{(n_{1},n_{2})}$. However, in practice, we only need the conditional law of $(\ell_1,\dots,\ell_K)$ given $\boldsymbol{X}^{(n_{1},n_{2})}$. We observe that if we define $\ell_{i,h}$ the number of unique values in $\{ T_{i,j}: X_{i,j}=X_h^{*} \}_{j}$ for fixed $i = 1,2$, since there is no intersection between the tables at different restaurants, $\ell_{h} := \ell_{1,h} + \ell_{2,h}$ for $h = 1,\dots,K$. If we denote as $n_{i,h}$ the number of observations equal to $X_{h}^*$ in group $i$, then the joint law of $\{\ell_{i,h}\}_{i,h}$ conditionally on $\boldsymbol{X}^{(n_{1},n_{2})}$ is proportional to 
\begin{equation} \label{eq:latent_posterior}
    \frac{c_{0}^{k}}{(c)_{n_{1}}(c)_{n_{2}}}\frac{c^{|\ell|}}{(c_{0})_{|\ell|}}\prod_{h=1}^{K}{(\ell_{h} - 1)!|\mathfrak{s}(n_{1,h},\ell_{1,h})||\mathfrak{s}(n_{2,h},\ell_{2,h})|\mathbbm{1}_{\{1,\dots,n_{1,h}\}}(\ell_{1,h})\mathbbm{1}_{\{1,\dots,n_{2,h}\}}(\ell_{2,h})},
\end{equation}
where $|\mathfrak{s}(n,\ell)|$ are the signless Stirling number of the first kind and $(a)_{q} = \tau(a + q)/\tau(a)$ is the rising factorial. The proof follows from specializing the partially exchangeable partition probability function (pEPPF) \cite{camerlenghi_distribution_2019} for a given configuration of $(\ell_1,\dots,\ell_K)$.

In principle, we could use the expression in \cref{eq:latent_posterior} to compute law of the latent tables $\boldsymbol{T}^{(n_{1},n_{2})}$ given $\boldsymbol{X}^{(n_{1},n_{2})}$. However, in practice, it becomes rapidly prohibitive since we have an unnormalized distribution and the normalization step can be time-consuming due to the size of the support, which increases with the number of observations. The most popular workaround is to implement a Gibbs sampler \citep{camerlenghi_distribution_2019}. After fixing an initial allocation of the tables $\boldsymbol{T}^{(n_{1},n_{2})}$ that satisfies the compatibility properties mentioned in \cref{compatible_TX}, for every $j = 1,\dots, n_{i}$ and $i = 1,2$ we remove $T_{i,j}$ and sample another value for it from the following discrete distribution,
\begin{equation} \label{gibbs}
    \mathbb{P}\bigl(T_{i,j} = T_{i,j^{*}} \big|\boldsymbol{X}^{(n_{1},n_{2})},\boldsymbol{T}_{-(i,j)}\bigr) = q_{i,j^{*}}, \qquad
    \mathbb{P}\bigl(T_{i,j} = T^{\star}\big|\boldsymbol{X}^{(n_{1},n_{2})},\boldsymbol{T}_{-(i,j)}\bigr) = \frac{c}{c_{0} + |\ell|}\ell_{h},
\end{equation}
where $\boldsymbol{T}_{-(i,j)}$ is the set $\boldsymbol{T}^{(n_{1},n_{2})}$ without $T_{i,j}$,  $q_{i,j^{*}}$ is the frequency of the table $T_{i,j^{*}}$ in $\boldsymbol{T}_{-(i,j)}$, $h$ is such that $X_{i,j} = X_{h}^{*}$,  $\ell_{h}$ is the number of unique values in $\boldsymbol{T}_{-(i,j)}$ associated to $X_{h}^{*}$, and $|\ell|$ is the number of unique values in $\boldsymbol{T}_{-(i,j)}$. Finally, $T^{\star} \sim H$ is a new value for the table.

\section{Algorithms for Numerical Simulations}

Once we have set all the distributional properties, we can estimate the kernel correlation a posteriori for an hDP model in two different ways: either a \emph{sampling-based} algorithm or an \emph{analytics-based} algorithm. We write $\boldsymbol{X}$ and $\boldsymbol{T}$ instead of $\boldsymbol{X}^{(n_{1},n_{2})}$ and $\boldsymbol{T}^{(n_{1},n_{2})}$ for compactness.

\subsection*{Sampling-Based Algorithm}

The first method consists of using the estimator defined in \cref{sec:estimator}. This estimator only uses our ability to generate samples from the model a posteriori.

Given the sequence of observable $\boldsymbol{X}$, we can initialize the sequence $\boldsymbol{T}$ to be i.i.d. from $H$. Then, we use the joint one-step-ahead predictive distribution in \cref{sampling_hDP} twice to generate a $2 \times 2$ sample for the observations and the tables. By discarding the future tables, we get a sample from $\mathcal{L}\bigl(X_{1,n_1+1}, X_{2,n_2+1}, X_{1,n_1+2}, X_{2,n_2+2}\big| \boldsymbol{X}, \boldsymbol{T}\bigr)$. Once this routine is completed, we update $\boldsymbol{T}$ conditionally to $\boldsymbol{X}$ using the Gibbs sampler introduced above.

The sampling procedure is repeated $M$ times to generate as many independent and identically distributed $2 \times 2$ samples $\bigl(X_{1,n_1+1}^{(t)}, X_{2,n_2+1}^{(t)}, X_{1,n_1+2}^{(t)}, X_{2,n_2+2}^{(t)}\bigr)_{t=1}^{M}$ and compute the sampling-based estimator as in \cref{prop_cor_asymp_normal}. 

See \cref{alg:smpl} for a thorough step-by-step description of the sampling-based method to compute the kernel correlation.

\begin{algorithm}[H]
\caption{Sampling-Based Algorithm for Kernel Correlation}\label{alg:smpl}
\begin{algorithmic}
    \Require $\mathtt{X}$, $P_{0}$, $H$, $\mathtt{c_{0}} \geq 0$, $\mathtt{c} \geq 0$, $\mathtt{M} \in \NN$.
    \State Inizialize $\mathtt{T}$ as an i.i.d. sample from $H$.
    \For{$\mathtt{t = 1,\dots,M}$}
        \State Sample $\mathtt{(X_{i,n_{i}+1}^{(t)},T_{i,n_{i}+1}^{(t)}), (X_{i,n_{i}+2}^{(t)},T_{i,n_{i}+2}^{(t)})}$ for $\mathtt{i = 1,2}$ according to \cref{sampling_hDP}.
        \State Update $\mathtt{T}$ through the Gibbs updating scheme \cref{gibbs}, conditionally on $\mathtt{X}$.
    \EndFor
    \State $\mathtt{varX_{1}} \gets \mathtt{\sum_{t=1}^{M}{k(X_{1,n_{1}+1}^{(t)},X_{1,n_{1}+2}^{(t)})}/(M-1) - \sum_{t=1}^{M}{\sum_{s=1}^{M}{k(X_{1,n_{1}+1}^{(t)},X_{1,n_{1}+2}^{(s)})}}/(M(M-1))}$
    \State $\mathtt{varX_{2}} \gets \mathtt{\sum_{t=1}^{M}{k(X_{2,n_{2}+1}^{(t)},X_{2,n_{2}+2}^{(t)})}/(M-1) - \sum_{t=1}^{M}{\sum_{s=1}^{M}{k(X_{2,n_{2}+1}^{(t)},X_{2,n_{2}+2}^{(s)})}}/(M(M-1))}$
    \State $\mathtt{covX} \gets \mathtt{\sum_{t=1}^{M}{k(X_{1,n_{1}+1}^{(t)},X_{2,n_{2}+1}^{(t)})}/(M-1) - \sum_{t=1}^{M}{\sum_{s=1}^{M}{k(X_{1,n_{1}+1}^{(t)},X_{2,n_{2}+1}^{(s)})}}/(M(M-1))}$
    \State \Return $\mathtt{corrX} \gets \mathtt{covX/\sqrt{varX_{1} varX_{2}}}$
\end{algorithmic}
\end{algorithm}

\subsection*{Analytics-Based Algorithm}

The second method uses our knowledge of the distributional properties of the posterior.

\medskip

\emph{Computation of the variances.}
If we apply the law of total variance, conditionally on $\boldsymbol{T}$ and $\tilde{P}_{0}$, using in particular \cref{eq:expression_mean_Ptildei} and the bilinearity of the covariance, we may write for $i = 1,2$,
\begin{equation} \label{varX}
    \var_{k}\bigl(\tilde{P}_{i}\big|\boldsymbol{X}\bigr) = \E\bigl[\var_{k}\bigl(\tilde{P}_{i}\big|\boldsymbol{X},\boldsymbol{T},\tilde{P}_{0}\bigr)\big|\boldsymbol{X}\bigr] + \frac{c^2}{(c + n_{i})^{2}}\var_{k}\bigl(\tilde{P}_{0}\big|\boldsymbol{X}\bigr).
\end{equation}

For the first summand, we use the expression of $\var_{k}\bigl(\tilde{P}_{i}\big|\boldsymbol{X},\boldsymbol{T},\tilde{P}_{0}\bigr)$ in \cref{var_cov_post_hDP_tabs}, and obtain 
\begin{equation} \label{varXT_sum}
    \E\bigl[\var_{k}\bigl(\tilde{P}_{i}\big|\boldsymbol{X},\boldsymbol{T},\tilde{P}_{0}\bigr)\big|\boldsymbol{X}\bigr] = \E\bigl[V_{i,1}(\boldsymbol{X},\boldsymbol{T}) + V_{i,2}(\boldsymbol{X},\boldsymbol{T}) + V_{i,3}(\boldsymbol{X},\boldsymbol{T}) \big| \boldsymbol{X}\bigr],
\end{equation}
with $V_{i,1}$, $V_{i,2}$ and $V_{i,3}$ defined below. Indeed, the first term in \cref{varXT_sum} is
\begin{equation*} 
    V_{i,1}(\boldsymbol{X},\boldsymbol{T}) := \frac{1}{2}\frac{1}{(c + n_{i} + 1)(c + n_{i})^{2}}\E\left[\iint{d_{k}^{2}(x,y)\diff\tilde{P}_{0}(x)\diff\tilde{P}_{0}(y)}\Bigg|\boldsymbol{X},\boldsymbol{T}\right].
\end{equation*}
Let us introduce
\begin{equation} \label{P0_star}
    P_{0}^{*} := \E\bigl[\tilde{P}_{0}\big|\boldsymbol{X},\boldsymbol{T}\bigr] = W_{0}(\boldsymbol{T})P_{0} + \sum_{h=1}^{K}{W_{h}(\boldsymbol{T})\delta_{X^{*}_{h}}},
\end{equation}
for $W_{0}(\boldsymbol{T}) = c_{0}/(c_{0} + |\ell|)$ and $W_{h}(\boldsymbol{T}) = \ell_{h}/(c_{0} + |\ell|)$ for $h = 1,\dots, K$. By quasi-conjugacy in \cref{post_hDP} and the computations a priori, 
\begin{align*}
    \var_k( \tilde{P}_0 | \boldsymbol{X},\boldsymbol{T} )&  =  \E\left[\iint{k(x,y)\diff\tilde{P}_{0}(x)\diff\tilde{P}_{0}(y)}\Bigg|\boldsymbol{X},\boldsymbol{T}\right] - \iint k(x,y) \diff P_{0}^{*}(x) \diff P_{0}^{*}(y) \\
    & = \frac{1}{2(c_{0} + |\ell| + 1)} \iint d_{k}^{2}(x,y)\diff P_{0}^{*}(x) \diff P_{0}^{*}(y).   
\end{align*}
Thus, we see that 
\begin{equation} \label{varXT_1}
    V_{i,1}(\boldsymbol{X},\boldsymbol{T}) = \frac{1}{2}\frac{1}{(c + n_{i} + 1)(c + n_{i})^{2}}\left(1 - \frac{1}{c_{0} + |\ell| + 1}\right)\iint{d_{k}^{2}(x,y) \diff P_{0}^{*}(x) \diff P_{0}^{*}(y)}.
\end{equation}

The second term in \cref{varXT_sum} is rewritten by the linearity of the expectation:
\begin{align}
    V_{i,2}(\boldsymbol{X},\boldsymbol{T}) &:= \E\left[\frac{cn_{i}}{(c + n_{i} + 1)(c + n_{i})^{2}}\iint{d_{k}^{2}(x,y)\diff\tilde{P}_{0}(x)\diff\hat P_{i}(y)}\Bigg|\boldsymbol{X},\boldsymbol{T}\right]
    \notag \\
    & =\frac{cn_{i}}{(c + n_{i} + 1)(c + n_{i})^{2}}\iint{d_{k}^{2}(x,y)\diff P_{0}^{*}(x)\diff\hat P_{i}(y)}
    \label{varXT_2}
\end{align}
with $P_{0}^{*}$ as in \cref{P0_star}. Lastly, the third term in \cref{varXT_sum} is
\begin{equation} \label{varXT_3}
    V_{i,3}(\boldsymbol{X},\boldsymbol{T}) := \frac{1}{2}\frac{n_{i}^{2}}{(c + n_{i} + 1)(c + n_{i})^{2}}\iint{d_{k}^{2}(x,y)\diff\hat P_{i}(x)\diff\hat P_{i}(y)},
\end{equation}
as the expectation is discarded, being the integrand completely determined by $\boldsymbol{X}$.

For the second summand in \cref{varX}, we need $\var_{k}\bigl(\tilde{P}_{0}\big|\boldsymbol{X}\bigr)$. We start from the definition
\begin{equation*} 
    \var_{k}\bigl(\tilde{P}_{0}\big|\boldsymbol{X}\bigr) = \E\left[ \left. \iint k(x,y) \diff \tilde{P}_0(x) \diff \tilde{P}_0(y) \right| \boldsymbol{X} \right] - \iint{k(x,y)\diff\E\bigl[\tilde{P}_{0}\big|\boldsymbol{X}\bigr](x)\diff\E\bigl[\tilde{P}_{0}\big|\boldsymbol{X}\bigr](y)}
\end{equation*}
For both expectations, we apply the tower law, conditioning to $\boldsymbol{T}$. That yields, with the functions $V_{0,1}$ and $V_{0,2}$ defined below, 
\begin{equation} \label{varX_0}
    \var_{k}\bigl(\tilde{P}_{0}\big|\boldsymbol{X}\bigr) = \E\bigl[V_{0,1}(\boldsymbol{X},\boldsymbol{T})\big|\boldsymbol{X}\bigr] - V_{0,2}(\boldsymbol{X}).
\end{equation}
Following the computations we did for $\var_{i,1}(\boldsymbol{X},\boldsymbol{T})$, the function $V_{0,1}$ is defined and rewritten as   
\begin{align}
    V_{0,1}(\boldsymbol{X},\boldsymbol{T}) & := \E\left[\left.\iint{k(x,y)\diff \tilde{P}_{0}(x)\diff \tilde{P}_{0}(y)}\right|\boldsymbol{X}, \boldsymbol{T}\right]
    \notag \\
    & = \iint{\left(k(x,y) + \frac{1}{2(c_{0} + |\ell| + 1)}d_{k}^{2}(x,y)\right) \diff P_{0}^{*}(x) \diff P_{0}^{*}(y)}.
    \label{v01}
\end{align}
On the other hand 
\begin{equation} \label{v02}
    V_{0,2}(\boldsymbol{X}) := \iint{k(x,y)\diff\E\bigl[P_{0}^{*}\big|\boldsymbol{X}\bigr](x)\diff\E\bigl[P_{0}^{*}\big|\boldsymbol{X}\bigr](y)},
\end{equation}
where, with the notations of \cref{P0_star},
\begin{equation} \label{eq:P0starmean}
    \E\bigl[P_{0}^{*}\big|\boldsymbol{X}\bigr] = \E[W_{0}(\boldsymbol{T})|\boldsymbol{X}]P_{0} + \sum_{h=1}^{K}{E[W_{h}(\boldsymbol{T})|\boldsymbol{X}]\delta_{X^{*}_{h}}}.
\end{equation}

\medskip

\emph{Computation of the covariance.}
For the covariance between $\tilde{P}_1$ and $\tilde{P}_2$, if we apply the law of total covariance, conditionally on $\boldsymbol{T}$ and $\tilde{P}_{0}$, and \cref{var_cov_post_hDP_tabs} we may write
\begin{equation} \label{covX}
    \cov_{k}\bigl(\tilde{P}_{1},\tilde{P}_{2}\big|\boldsymbol{X}\bigr) = \frac{c^2}{(c + n_{1})(c + n_{2})}\var_{k}\bigl(\tilde{P}_{0}\big|\boldsymbol{X}\bigr).
\end{equation} 
Note that we have used \cref{eq:expression_mean_Ptildei} again and the bilinearity of the covariance. The expression for $\var_{k}\bigl(\tilde{P}_{0}\big|\boldsymbol{X}\bigr)$ has already been expanded, see \cref{varX_0} above.

\medskip

\emph{Running the Gibbs sampler.}
Given the sequence of observable $\boldsymbol{X}$, we can initialize the sequence $\boldsymbol{T}$ to be i.i.d. from $H$. Then, we can update $\boldsymbol{T}$ $R$ times conditionally on $\boldsymbol{X}$ using the Gibbs sampler introduced above. Hence, we obtain a sequence of $(\boldsymbol{T}_{1},\dots,\boldsymbol{T}_{R})$. We approximate the expression in \cref{varXT_sum} as
\[
    \E\bigl[\var_{k}\bigl(\tilde{P}_{i}\big|\boldsymbol{X},\boldsymbol{T},\tilde{P}_{0}\bigr)\big|\boldsymbol{X}\bigr] \approx \frac{1}{R}\sum_{r=1}^{R}{V_{i,1}(\boldsymbol{X},\boldsymbol{T}_{r}) + V_{i,2}(\boldsymbol{X},\boldsymbol{T}_{r}) + V_{i,3}(\boldsymbol{X},\boldsymbol{T}_{r})},
\]
using for that the expressions \cref{varXT_1}, \cref{varXT_2}, and \cref{varXT_3}. The first term appearing in the expression of $\var_{k}\bigl(\tilde{P}_{0}\big|\boldsymbol{X}\bigr)$ in \cref{varX_0} is approximated as
\[
    \E\bigl[\var_{0,1}(\boldsymbol{X},\boldsymbol{T})\big|\boldsymbol{X}\bigr] \approx \frac{1}{R}\sum_{r=1}^{R}{V_{0,1}(\boldsymbol{X},\boldsymbol{T}_{r})},
\]
using \cref{v01}. Lastly, to approximate $V_{0,2}(\boldsymbol{X})$ we use \cref{v02} together with \cref{eq:P0starmean}, where in the last formula we estimate the expectation of the weights $W_{0}(\boldsymbol{T})$, $W_{1}(\boldsymbol{T})$, \dots, $W_{K}(\boldsymbol{T})$ conditionally on $\boldsymbol{X}$ as 
\[
    \E[W_{h}(\boldsymbol{T})|\boldsymbol{X}] \approx \frac{1}{R}\sum_{r=1}^{R}{W_{h}(\boldsymbol{T_{r}})}
\]
for $h = 0,\dots,K$. 

Using these expressions, we can approximate the expression for the covariance in \cref{covX} and the variances in \cref{varX} a posteriori to obtain the corresponding correlation. See \cref{alg:anal} for a detailed, step-by-step description of the analytics-based method for computing the kernel correlation.

\medskip

\emph{Approximating integrals with respect to $P_0$.}
Notice that we need to approximate integrals with respect to $P_0$. To that end we approximate $P_{0}$ as $\sum_{t=1}^{M}{\delta_{Z_{t}}/M}$, where $Z_{1},\dots,Z_{M} \stackrel{\text{i.i.d}}{\sim} P_{0}$. Hence, all the integrals above can be rewritten as finite sums since we are integrating with respect to discrete measures. Specifically, we apply the approximation
\begin{align*}
\mathcal{I}_{P_{0}} & := \iint{k(x,y)\diff P_{0}(x) \diff P_{0}(y)} \approx \frac{1}{M^{2}}\sum_{t=1}^{M}{\sum_{r=1}^{M}{k(Z_{t},Z_{r})}}, & \\
\mathcal{I}_{P_{0},h} & := \iint{k(x,y)\diff P_{0}(x) \diff \delta_{X_h^*}(y)} \approx \frac{1}{M}\sum_{t=1}^{M}{k(Z_{t},X_h^*)} & \text{for } h=1, \ldots, K.
\end{align*}

\begin{algorithm}[H]
\caption{Analytics-Based Algorithm for Kernel Correlation}\label{alg:anal}
\begin{algorithmic}
    \Require $\mathtt{X}$, $P_{0}$, $H$, $\mathtt{c_{0}} \geq 0$, $\mathtt{c} \geq 0$, $\mathtt{R} \in \NN$.
    \State Compute $\mathtt{X^{*}_{1},\dots,X^{*}_{K}}$, $\mathtt{n_{1,1},\dots,n_{1,K}}$, and $\mathtt{n_{2,1},\dots,n_{2,K}}$.
    \State \emph{// Precompute integrals w.r.t. to $P_0$}
    \State Sample $\mathtt{Z_{1},\dots,Z_{M}}$ i.i.d. from $P_{0}$.
    \State $\mathtt{I_{P_{0}}} \gets \mathtt{\sum_{t=1}^{M}{\sum_{s=1}^{M}{k(Z_{t},Z_{s})}}/M^{2}}$
    \For{$\mathtt{h = 1,\dots,K}$}
    \State $\mathtt{I_{P_{0},h}} \gets \mathtt{\sum_{t=1}^{M}{k(Z_{t},X_{h}^*)}/M}$
    \EndFor
    \State \emph{// Run the Gibbs sampler}
    \State Inizialize $\mathtt{T}$ as an i.i.d. sample from $H$.
    \For{$\mathtt{r = 1,\dots,R}$}
        \State Compute $\mathtt{l_{1},\dots,l_{K}}$.
        \State Compute $\mathtt{W_{0}^{(r)},W_{1}^{(r)},\dots,W_{K}^{(r)}}$.
        \State $\mathtt{VarXT_{1}^{(r)}} \gets \mathtt{V_{1,1}(X,T) + V_{1,2}(X,T) + V_{1,3}(X,T) + c^{2}V_{0,1}(X,T)/(c+n_{1})^{2}}$
        \State $\mathtt{VarXT_{2}^{(r)}} \gets \mathtt{V_{2,1}(X,T) + V_{2,2}(X,T) + V_{2,3}(X,T) + c^{2}V_{0,1}(X,T)/(c+n_{2})^{2}}$
        \State $\mathtt{CovXT^{(r)}} \gets \mathtt{c^{2}V_{0,1}(X,T)/((c+n_{1})(c+n_{2}))}$
        \State Update $\mathtt{T}$ through the Gibbs updating scheme \cref{gibbs}, conditionally on $\mathtt{X}$.
    \EndFor
    \State \emph{// Evalute $\E[P_0^*|\boldsymbol{X}]$ and $V_{0,2}(\boldsymbol{X})$}
    \State Compute $\mathtt{\overline{W}_{0}},\mathtt{\overline{W}_{1}},\dots,\mathtt{\overline{W}_{K}}$.
    \State $\mathtt{V_{0,2}} \gets \mathtt{\overline{W}_{0}^{2}I_{P_{0}} + 2\overline{W}_{0}\sum_{h=1}^{K}{\overline{W}_{h} I_{P_0,h}}  + \sum_{h=1}^{K}{\sum_{j=1}^{K}{\overline{W}_{h}\overline{W}_{j}k(X^{*}_{h},X^{*}_{j})}}}$
    \State \emph{// Merge all computations together and return the correlation}
    \State $\mathtt{varX_{1}} \gets \mathtt{\overline{VarXT_{1}} - c^{2}V_{0,2}/(c+n_{1})^{2}}$
    \State $\mathtt{varX_{2}} \gets \mathtt{\overline{VarXT_{2}} - c^{2}V_{0,2}/(c+n_{2})^{2}}$
    \State $\mathtt{covX} \gets \mathtt{\overline{CovXT} - c^{2}V_{0,2}/((c+n_{1})(c+n_{2}))}$ 
    \State \Return $\mathtt{corrX} \gets \mathtt{covX/\sqrt{varX_{1} varX_{2}}}$
\end{algorithmic}
\end{algorithm}

\section{Computations for Model Comparison} \label{SM:model_comparison}

In \cref{sec:models}, we compare the Gaussian case in \cref{par_ex} and the hDP model in \cref{hDP} with $P_{0} = \mathcal{N}(0,t^{2})$ for a Gaussian kernel $k(x,y) = \exp\bigl(-(x-y)^{2}/(2\sigma^{2})\bigr)$ with parameter $\sigma > 0$.

To perform simulations for different values of the kernel correlation, we set the parameters so that the observables have the same marginal distributions. Since the marginal distribution of the observables is $\mathcal{N}(0,s^{2} + \tau^{2})$ for the Gaussian model, while it is $\mathcal{N}(0,t^{2})$ for the hDP model, we need to set the constraint $t^{2} = s^{2} + \tau^{2}$.

In a similar line of reasoning, we set the kernel variances a priori to be equal to the same value $v > 0$ for each group for each case. To compute the kernel variances, we use the identity: if $\kappa(x,y) = a \exp(-(x-y)^2/(2 b^2))$ while $P$ is $\mathcal{N}(0,c^2)$, then 
\begin{equation} \label{eq:identity_aux}
    \int \kappa(x,x) \diff P(x) - \iint \kappa(x,y) \diff P(x) \diff P(y) = a \left( 1 - \sqrt{\frac{b^2}{2 c^2 + b^2}} \right).
\end{equation}
We deduce the kernel variance for the parametric model using \cref{mix_corr} and the explicit expression of \cref{int_ker_cov_var}: we apply \cref{eq:identity_aux} with $c^{2} = \tau^{2}$; while, from \cref{mix_ex}, we have $a^{2} = \sqrt{\sigma^2 / (2 s^2 + \sigma^2)}$ and $b^{2} = 2s^2 + \sigma^2$. We also deduce the kernel variance for the nonparametric model from \cref{var_set_ker} and \cref{sec:corr_hdp}: in the case we use~\eqref{eq:identity_aux} with $a=1$, $b^{2} = \sigma^2$ and $c^2 = t^2$. By imposing that both variances are equal to $v$, we obtain the constraint
\[
    v = \sqrt{\frac{\sigma^{2}}{2s^{2} + \sigma^{2}}} - \sqrt{\frac{\sigma^{2}}{2\tau^{2} + 2s^{2} + \sigma^{2}}} = \frac{1 + c + c_{0}}{(1 + c)(1 + c_{0})}\left(1 - \sqrt{\frac{\sigma^{2}}{2t^{2} + \sigma^{2}}}\right).
\]

Now, if we set the kernel correlation to be equal to a value $\xi \in [0,1]$, we can determine the parameters for both the Gaussian case and the hDP case. In other words, once we fix $v$, $\xi$ and $t^{2}$ we can determine $s^{2}$, $\tau^{2}$, and $\rho$ for the Gaussian case, and $c_{0}$ and $c$ for the hDP case.

For the Gaussian case, we have to solve the following system.
\[
\begin{cases}
    t^{2} &= \tau^{2} + s^{2}, \\
    v &= \sqrt{\frac{\sigma^{2}}{2s^{2} + \sigma^{2}}} - \sqrt{\frac{\sigma^{2}}{2\tau^{2} + 2s^{2} + \sigma^{2}}}, \\
    v\xi &= \sqrt{\frac{\sigma^{2}}{2\tau^{2}(1 - \rho) + 2s^{2} + \sigma^{2}}} - \sqrt{\frac{\sigma^{2}}{2\tau^{2} + 2s^{2} + \sigma^{2}}},
\end{cases}
\]
which leads to
\[
\begin{cases}
    s^{2} &= \frac{\sigma^{2}}{2}\left(\left(v + \sqrt{\frac{\sigma^{2}}{2t^{2} + \sigma^{2}}}\right)^{-2} - 1\right) \\
	\tau^{2} &= t^{2} - s^{2} \\
	\rho &= \frac{t^{2} -\frac{\sigma^{2}}{2}\left(\left(v\xi + \sqrt{\frac{\sigma^{2}}{2t^{2} + \sigma^{2}}}\right)^{-2} - 1\right)}{t^{2} - \frac{\sigma^{2}}{2}\left(\left(v + \sqrt{\frac{\sigma^{2}}{2t^{2} + \sigma^{2}}}\right)^{-2} - 1\right)},
\end{cases}
\]
which is solvable with $s^{2}, \tau^{2} > 0$ and $\rho \in [0,1]$ for $v \in (0,1)$ and $t^{2} > \sigma^{2}/2\bigl(1/(1-v)^{2} - 1\bigr)$.

For the hDP case, we have to solve the system
\[
\begin{cases}
    v &= \frac{1 + c + c_{0}}{(1 + c)(1 + c_{0})}\left(1 - \sqrt{\frac{\sigma^{2}}{2t^{2} + \sigma^{2}}}\right), \\
	v\xi &= \frac{1}{1 + c_{0}}\left(1 - \sqrt{\frac{\sigma^{2}}{2t^{2} + \sigma^{2}}}\right),
\end{cases}
\]
which leads to 
\[
\begin{cases}
    c_{0} &= \frac{1}{v\xi}\left(1 - \sqrt{\frac{\sigma^{2}}{2t^{2} + \sigma^{2}}}\right) - 1, \\
    c &= \frac{1}{1 - \xi}\left(\frac{1}{v}\left(1 - \sqrt{\frac{\sigma^{2}}{2t^{2} + \sigma^{2}}}\right) - 1\right),
\end{cases}
\]
which, again, is solvable with $c_{0},c > 0$ for $v \in (0,1)$ and $t^{2} > \sigma^{2}/2\bigl(1/(1-v)^{2} - 1\bigr)$.

To complete the explanation of the model comparison in \cref{sec:models}, we need to characterize the one-step-ahead posterior predictive for both the Gaussian case in \cref{par_ex} and the hDP model in \cref{hDP} with $P_{0} = \mathcal{N}(0,t^{2})$. For the Gaussian case, we have from \cref{par_ex} in \cref{sec:examples} that the posterior distribution of $\boldsymbol{\theta} = (\theta_{1},\theta_{2})$ given $\boldsymbol{X}^{(n_{1},n_{2})}$ is $\mathcal{N}(\boldsymbol{\theta}^*, \Sigma^*)$ with $\boldsymbol{\theta}^{*}$ and $\Sigma^{*}$ are as in \cref{eq:theta_star,eq:sigma_star}, respectively. Consequently, the one-step-ahead predictive distribution for the $i$-th group is $\mathcal{N}(\theta^{*}_{i},s^{2} + \Sigma^{*}_{i,i})$ for $i = 1,2$. For the hDP model, we generate a sample from the posterior predictive distribution using a Gibbs sampler similar to the one described in \cref{sec:hdp_prop}. For each sampled value, we update the allocation of the tables $\boldsymbol{T}^{(n_{1},n_{2})}$ conditionally on $\boldsymbol{X}^{(n_{1},n_{2})}$. Then, we generate a new data point from the posterior augmented model using the one-step-ahead predictive distribution in \cref{sampling_hDP}.

\bibliography{rkhs}